\documentclass[11pt,twoside,leqno]{aomamlt2e} 
  \pageno{823}
\received{March 3, 2003}

\newtheorem{thm}{Theorem}[section]
\newtheorem{lemma}[thm]{Lemma}
\newtheorem{cor}[thm]{Corollary}
\newtheorem{prop}[thm]{Proposition}

\theorembodyfont{\upshape}
\newtheorem{defi}[thm]{{\it Definition}}

\def\CC{{\cal{C}}}

\def\CX{{\cal{X}}}

\def\bZ{{\Bbb Z}}
\def\bQ{{\Bbb Q}}
\def\bF{{\Bbb F}}

\def\Ker{\operatorname{Ker}\nolimits}

\def\Dim{\operatorname{Dim}\nolimits}
\def\Res{\operatorname{Res}\nolimits}

\def\Inf{\operatorname{Inf}\nolimits}
\def\Def{\operatorname{Def}\nolimits}
\def\Defres{\operatorname{DefRes}\nolimits}
\def\Ten{\operatorname{Ten}\nolimits}

\def\Iso{\operatorname{Iso}\nolimits}

\def\Soc{\operatorname{Soc}\nolimits}
\def\Out{\operatorname{Out}\nolimits}

\def\Rarr#1{\buildrel #1\over \longrightarrow}

\def\HHH{\operatorname{H}\nolimits}
\def\Hom{\operatorname{Hom}\nolimits}

\def\End{\operatorname{End}\nolimits}
\def\HH#1#2#3{\HHH^{#1}(#2,#3)}

\def\Ext{\operatorname{Ext}\nolimits}

\def\hgs{\operatorname{H^*\text{$(G,k)$}}\nolimits}

\def\res#1#2{\text{res}_{#1,#2}}
\def\ress#1#2{\text{res}_{#1,#2}^*}
\def\resgh{\text{res}_{G,H}}

\def\VG{{V_G(k)}}

\def\bfp{{\Bbb F}_p}

\def\bF2{{\Bbb F}_2}
\def\bF4{{\Bbb F}_4}

\begin{document}
\currannalsline{162}{2005} 

\title{The classification of torsion\\ endo-trivial modules}

 \acknowledgements{}
%\author{}
\twoauthors{Jon F. Carlson$^\ast$}{Jacques Th\'evenaz}

 \institution{University of Georgia, Athens, GA
\\
 \email{jfc@sloth.math.uga.edu}
\\
\vglue-9pt
\'Ecole Polytechnique F\'ed\'erale de lausanne, Lausanne, Switzerland
\\
\email{jacques.thevenaz@epfl.ch}}

 \shorttitle{The classification of torsion endo-trivial modules}
\shortname{Jon Carlson and Jacques Th\'evenaz}

\section{Introduction}

This paper settles a problem raised at the end of the seventies by J.L.
Alperin \cite{Alpinv}, E.C. Dade \cite{Dendo} and J.F. Carlson
\cite{Cendo}, namely the classification of torsion endo-trivial modules for
a finite $p$-group over a field of characteristic~$p$. Our results also
imply, at least when $p$ is odd, the complete classification of torsion
endo-permutation modules.

We refer to \cite{CT1} and \cite{Both} for an overview of the problem and
its importance in the representation theory of finite groups. Let us only
mention that the classification of endo-trivial modules is the crucial
step for understanding the more general class of endo-permutation modules,
and that endo-permutation modules play an important role in module theory,
in particular as source modules, in block theory where they appear in the
description of source algebras, and in both derived equivalences and
stable equivalence of block algebras, for which many new developments
have appeared recently.

Let $G$ be a finite $p$-group and $k$ be a field of characteristic $p$.
Recall that a (finitely generated) $kG$-module~$M$ is called endo-trivial
if ${\End_k(M)\cong k\oplus F}$ as $kG$-modules, where $F$ is a free module.
Typical examples of endo-trivial modules are the Heller translates
$\Omega^n(k)$ of the trivial module. Any endo-trivial $kG$-module~$M$ is a
direct sum $M=M_0\oplus L$, where $M_0$ is an indecomposable endo-trivial
$kG$-module and $L$ is free.
Conversely, by adding a free module to an endo-trivial module, we always
obtain an endo-trivial module. This defines an equivalence relation among
endo-trivial modules and each equivalence class contains exactly one
indecomposable module up to isomorphism. The set $T(G)$ of all equivalence
classes of endo-trivial $kG$-modules is a group with multiplication induced
by tensor product, called simply the group of endo-trivial $kG$-modules.
Since scalar extension of the coefficient  field induces an injective map
between the groups of endo-trivial modules, we can replace $k$ by its
algebraic closure. So we assume that $k$ is algebraically closed. We refer
to \cite{CT1} for more details about $T(G)$.
\vfill
\footnoterule
{\footnotesize$^\ast$The first author was partly supported by a grant from NSF.}\eject

Dade~\cite{Dendo} proved that if $A$ is a noncyclic abelian $p$-group then
$T(A)\cong \bZ$, generated by the class of~$\Omega^1(k)$. For any
$p$-group~$G$, Puig~\cite{Pu1} proved that the abelian group $T(G)$ is
finitely generated (but we do not use this here since it is actually a
consequence of our main results). The torsion-free rank of $T(G)$ has been
determined recently by Alperin \cite{Alpsyz} and the remaining problem lies
in the structure of the torsion subgroup $T_t(G)$.

Let us first recall some important known cases (see \cite{CT1}).
If $G=1$ or $G=C_2$, then $T(G)=0$.
If $G=C_{p^n}$ is cyclic of order~$p^n$, with $n\geq 1$ if $p$ is odd and
$n\geq 2$ if $p=2$, then $T(C_{p^n}) \cong \bZ/2\bZ$ (generated by the class
of~$\Omega^1(k)$). If $G=Q_{2^n}$ is a
quaternion group of order $2^n\geq 8$, then $T(Q_{2^n}) = T_t(Q_{2^n})
\cong \bZ/4\bZ \,\oplus\, \bZ/2\bZ$.
If $G={\rm SD}_{2^n}$ is a semi-dihedral group of order $2^n\geq 16$, then
$T({\rm SD}_{2^n}) \cong \bZ \oplus \bZ/2\bZ$ and so $T_t({\rm SD}_{2^n}) \cong
\bZ/2\bZ$. Our first main result asserts that these are the only
cases where nontrivial torsion occurs.

\begin{thm} \label{no-torsion}
Suppose that $G$ is a finite $p$\/{\rm -}\/group which is not cyclic{\rm ,} quaternion{\rm ,} or
semi\/{\rm -}\/dihedral. Then $T_t(G) = \{0\}$.
\end{thm}

As explained in \cite{CT1}, the computation of the torsion subgroup $T_t(G)$
is tightly connected to the problem of detecting nonzero elements of $T(G)$
on restriction to a suitable class of subgroups. A detection theorem was
proved in~\cite{CT1} and it was conjectured that the detecting family should
actually only consist of elementary abelian subgroups of rank at most~2 and,
in addition when $p=2$, cyclic groups of order~4 and quaternion subgroups
$Q_8$ of order~8. This conjecture is correct and the largest part of the
present paper is concerned with the proof of this conjecture.

It is in fact only for the cases of cyclic, quaternion, and semi-dihedral
groups that one needs to include cyclic groups $C_p$ or $C_4$ and quaternion
subgroups $Q_8$ in the detecting family. For all the other cases, we are
going to prove the following.

\begin{thm} \label{finaldetection}
Suppose that $G$ is a finite $p$\/{\rm -}\/group which is not cyclic{\rm ,} quaternion{\rm ,} or
semi\/{\rm -}\/dihedral. Then the restriction homomorphism
$$
\prod_E \Res_E^G \; : T(G) \Rarr{} \prod_{E} T(E) \; \cong
\; \prod_{E} \bZ
$$
is injective{\rm ,} where $E$ runs through the set of all elementary abelian
subgroups of rank~$2$.
\end{thm}

In order to explain the right-hand side isomorphism, recall that $T(E)\cong
\bZ$ by Dade's theorem \cite{Dendo}.
Notice that Theorem~\ref{no-torsion} follows immediately from
Theorem~\ref{finaldetection}.

In the case of the theorem, $T(G)$ is free abelian and the method of Alperin
\cite{Alpsyz} describes its rank by restricting drastically the list of
elementary abelian subgroups which are actually needed on the right-hand
side (see also
\cite{Both} for another approach). However, for a complete
classification of all endo-trivial modules, there is still an open
problem. Alperin's method shows that $T(G)$ is a full lattice in a free
abelian group~$A$ by showing that some explicit subgroup~$S(G)$ of the
same rank satisfies $S(G) \subseteq T(G) \subseteq A$. But there is still
the problem of describing explicitly the finite group $T(G)/S(G)
\subseteq A/S(G)$. However, this additional problem only occurs if $G$
contains maximal elementary subgroups of rank~2 (see \cite{Alpsyz} or
\cite{Both} for details). In all other cases the rank of~$T(G)$ is one and
we have the following result.

\begin{cor} \label{rank-one}
Suppose that $G$ is a finite $p$\/{\rm -}\/group for which every maximal elementary
abelian subgroup has rank at least~$3$. Then $T(G) \cong \bZ${\rm ,} 
generated by the class of the module $\Omega^1(k)$.
\end{cor}

For the proof of Theorem~\ref{finaldetection}, we first use the results of
\cite{CT1} which provide a reduction to the case of extraspecial and almost
extraspecial $p$-groups. These are the difficult cases for which we need to
prove that the groups can be eliminated from the detecting family. When $p$
is odd, this was already done in \cite{CT1} for extraspecial $p$-groups of
exponent~$p^2$ and almost extraspecial $p$-groups. So we are left with the
remaining cases and we have to prove the following theorem, which is in fact
the main result we prove in the present paper.

\begin{thm} \label{mainthm}
Suppose the following\/{\rm :}\/
\begin{itemize}
\ritem{(a)} If $p=2$\/{\rm ,}\/ $G$ is an extraspecial or almost extraspecial $2$-group and
$G$ is not isomorphic to $Q_8$.

\ritem{(b)} If $p$ is odd\/{\rm ,}\/ $G$ is an extraspecial $p$-group of exponent~$p$.
\end{itemize}
Then the restriction homomorphism
$$
\prod_H \Res_H^G \; : T(G) \Rarr{} \prod_{H} T(H)
$$
is injective{\rm ,} where $H$ runs through the set of all maximal subgroups of
$G$.
\end{thm}

As mentioned earlier, the classification of endo-trivial modules has
immediate consequences for the more general class of
endo-permutation modules. The second goal of the present paper is to
describe the consequences of the main results for the classification of
torsion endo-permutation modules. We prove a detection theorem for the Dade
group of all endo-permutation modules and also a detection
theorem for the torsion subgroup of the Dade group. For odd $p$,
this yields a complete description of this torsion subgroup, by the
results of~\cite{Both}.

\begin{thm} \label{odd-Dade-group}
If $p$ is odd and $G$ is a finite $p$-group{\rm ,} the torsion subgroup of the
Dade group of all endo-permutation $kG$-modules is isomorphic to
$(\bZ/2\bZ)^s${\rm ,} where $s$ is the number of conjugacy classes of nontrivial
cyclic subgroups of~$G$.
\end{thm}

One set of $s$ generators is described in~\cite{Both}. Since an element of
order~2 corresponds to a self-dual module, we obtain in particular the
following corollary.

\begin{cor} \label{self-dual}
If $p$ is odd and $G$ is a finite $p$-group{\rm ,} then an indecomposable
endo-permutation $kG$-module $M$ with vertex~$G$ is self-dual if and only if
the class of $M$ in the Dade group is a torsion element of this group.
\end{cor}

This is an interesting result in view of the fact that many invariants lying
in the Dade group (e.g.\ sources of simple modules) are either known or
expected to lie in the torsion subgroup, while it is not at all clear why
the modules should be self-dual.

When $p=2$, the situation is more complicated but we obtain that any torsion
element of the Dade group has order 2 or~4. Moreover, the detection result
is efficient in some cases, but examples also show that it is not always
sufficient to determine completely this torsion subgroup.

Theorem~\ref{mainthm} is the result whose proof requires most of the work.
The result has to be treated separately when $p=2$ or when $p$ is
odd. However, the strategy is similar and many of the same methods are of
use for the proof in both cases. After a preliminary Section 2 and two
sections about the cohomology of extraspecial groups, the proof of
Theorem~\ref{mainthm} occupies Sections~5--11. We use a large amount of
group cohomology, including some very recent results, as well as the
theory of support varieties of  modules. The crucial role of Serre's theorem
on products of Bocksteins appears once again and we actually need a bound
for the number of terms in this product that was recently obtained by
Yal\c{c}in~\cite{Yal} for (almost) extraspecial groups. Also, the
module-theoretic counterpart of Serre's theorem described
in~\cite{Celeab} plays a crucial role. All these results allow us to find
an upper bound for the dimension of an indecomposable endo-trivial module
which is trivial on restriction to proper subgroups. For the purposes of
the present paper, we shall call such a module a {\it critical\/} module.
The main goal is to prove that there are no nontrivial critical
modules for extraspecial and almost extraspecial 2-groups, except~$Q_8$,
and also none for extraspecial $p$-groups of exponent~$p$ (with $p$ odd).

The existence of a bound for the dimension of a critical module
had been known for more than 20 years and was  used by
Puig~\cite{Pu1} in his   proof of the finite generation of~$T(G)$. The new
aspect is that we are now able to control this bound for (almost)
extraspecial groups. One of the differences between the case where $p=2$
and the case where $p$ is odd lies in the fact that the cohomology of
extraspecial 2-groups is entirely known, so that a reasonable bound can be
computed, while for odd $p$ some more estimates are necessary. Another
difference is due to the fact that we have three families of groups to
consider when $p=2$, but only one when $p$ is odd, because the other two
were already dealt with in~\cite{CT1}.

The other main idea in the proof of Theorem~\ref{mainthm} is the
following. Under the assumption that there exists a nontrivial critical
module~$M$, we can construct many others using the action of $\Out(G)$
(which is an orthogonal or symplectic group since $G$ is (almost)
extraspecial), and then construct a very large critical module by taking
tensor products. The dimension of this large module exceeds the upper
bound mentioned above and we have a contradiction. It is this part in
which the theory of varieties associated to modules plays an essential
role.  We use it to analyze a suitable quotient module~$\overline{M}$
which turns out to be periodic as a module over the elementary abelian
group $\overline{G}=G/\Phi(G)$.

Once Theorem~\ref{mainthm} is proved, the proof of
Theorem~\ref{finaldetection} requires much less machinery and appears in
Section~12. It is very easy if $p$ is odd and, if $p=2$, it is essentially
an inductive argument using a group-theoretical lemma.
Theorem~\ref{no-torsion} also follows easily.

The paper ends with two sections about the Dade group of all
endo-permutation modules, where we prove the results mentioned above.

We wish to thank numerous people who have shared ideas and opinions in
the course of the writing of this paper. Special thanks are due to C\'edric Bonnaf\'e, Roger Carter, Ian Leary,
Gunter Malle, and Jan Saxl.  The first author also wishes to thank the Humboldt Foundation for supporting
his stay in Germany while this paper was being written.

\section{Preliminaries}

Recall that $G$ denotes a finite $p$-group, and $k$ an
algebraically closed field of characteristic $p$.  In this section we write
down some of the facts about modules and support varieties
that we will need in later developments.   All $kG$-modules are
assumed to be finitely generated.

Recall that every projective $kG$-module is free, because
$G$ is a $p$-group, and that injective and projective modules coincide.
Moreover, an indecomposable $kG$-module
$M$ is free if and only if $t_1^G\cdot M\neq 0$, where $t_1^G=\sum_{g\in
G} g$ (a generator of the socle of~$kG$). More generally, if $M$ is a
$kG$-module and if
$m_1, \dots, m_r\break \in M$ are such that $t_1^G m_1, \dots, t_1^G m_r$
are linearly independent, then
$m_1, \dots, m_r$ generate a free submodule $F$ of $M$ of rank $r$.
Moreover $F$ is a direct summand of $M$ because $F$ is also injective.

Suppose that $M$ is a $kG$-module. If
$P \Rarr{\theta} M$ is a projective cover of $M$ then
we let $\Omega(M)$ denote the kernel of~$\theta$.  We can iterate the
process and define inductively $\Omega^n(M) = \Omega(\Omega^{n-1}(M))$,
for $n > 1$.  Suppose that $M \Rarr{\mu}Q$ is an injective hull of~$M$.
Recall that $Q$ is a projective as well as injective module. Then we let
$\Omega^{-1}(M)$ be the cokernel of~$\mu$.  Again we have inductively that
$\Omega^{-n}(M) = \Omega^{-1}(\Omega^{-n+1}(M))$ for $n > 1$.  The
modules $\Omega^n(M)$ are well defined up to isomorphism and they have no
nonzero projective submodules.  In general we write $M = \Omega^0(M)
\oplus P$ where $P$ is projective and $\Omega^0(M)$ has no projective
summands.

The basic calculus of the syzygy modules $\Omega^n(M)$ is
expressed in the following.

\begin{lemma}  Suppose that $M$ and $N$ are $kG$-modules. Then
$\Omega^m(M) \otimes \Omega^n(N) \cong \Omega^{m+n}(M \otimes N) \oplus
({\rm free})$.
\end{lemma}

Here $M \otimes N$ is meant to be  the tensor product $M \otimes_k N$ over
$k$,
with the action of the group $G$ defined diagonally,
$g(m \otimes n) = gm\otimes gn$.  The proof of the lemma is a consequence
of the facts that $M\otimes_k - $ and $-\otimes_k N$ preserve exact
sequences
and that $M \otimes N$ is projective whenever either $M$ or $N$ is a
projective module.

The cohomology ring $\hgs$ is a finitely generated $k$-algebra and for
any $kG$-modules $M$ and $N$, $\Ext^*_{kG}(M,N)$ is a finitely generated
module over $\hgs \cong \Ext^*_{kG}(k,k)$.  We let $V_G(k)$ denote the
maximal ideal spectrum of $\hgs$.  For any $kG$-module $M$, let $J(M)$ be
the annihilator in $\hgs$ of the cohomology ring $\Ext^*_{kG}(M,M)$.  Let
$V_G(M) = V_G(J(M))$ be the closed subset of $V_G(k)$ consisting of all
maximal ideals that contain $J(M)$.  So $V_G(M)$ is a homogeneous affine
variety. We need some of the properties of
support varieties in essential ways in the course of our proofs.
See the general references \cite{Bbook}, \cite{Ev} for more
explanations and details.

\begin{thm}\label{propvar}  Let $L, M$ and $N$ be $kG$-modules.
\begin{enumerate}
\item $V_G(M) = \{0\}$ if and only if $M$ is projective.
\item If $\, 0 \rightarrow L \rightarrow M \rightarrow N \rightarrow 0$
is exact then the variety of any one of $L, M$ or $N$ is contained in the
union of the varieties of the other two.
Moreover{\rm ,} if $V_G(L) \cap V_G(N) = \{0\}${\rm ,} then
the sequence splits.
\item $V_G(M\otimes N) = V_G(M) \cap V_G(N)$.
\item $V_G(\Omega^n(M)) = V_G(M) = V_G(M^*)$ where $M^* =
\Hom_k(M,k)$ is the $k$-dual of~$M$.
\item If $V_G(M) = V_1 \cup V_2$ where $V_1$ and $V_2$ are nonzero closed
subsets of $V_G(k)$ and $V_1 \cap V_2 = \{0\}${\rm ,} then $M \cong M_1 \oplus
M_2$ where $V_G(M_1) = V_1$ and $V_G(M_2)\break = V_2$.
\item A nonprojective module $M$ is periodic  \/{\rm (}\/i.e.\ for some $n > 0${\rm ,}
$\Omega^n(M)
\cong \Omega^0(M)${\rm )} if and only if its variety $V_G(M)$ is a union of lines
through the origin in $\VG$.
\item Let $\zeta \in \Ext^n_{kG}(k,k) = \HH{n}{G}{k}$ be represented by the
\/{\rm (}\/unique\/{\rm )}\/ 
cocycle $\zeta: \Omega^n(k) \Rarr{} k$ and let $L=\Ker(\zeta)${\rm ,} so that
there is an exact sequence
$$
0 \Rarr{} L \Rarr{} \Omega^n(k) \Rarr{\zeta} k \Rarr{} 0 \,.
$$
Then $V_G(L) = V_G(\zeta)${\rm ,} the variety of the ideal generated by $\zeta${\rm ,}
consisting of all maximal ideals containing $\zeta$.
\end{enumerate}
\end{thm}

We are particularly interested in the case in which the group $G$ is an
elementary abelian group.  First assume that $p = 2$ and
$G = \langle x_1, \dots,
x_n\rangle \cong (C_2)^n$.  Then $\hgs \cong k[\zeta_1,
\dots, \zeta_n]$ is a polynomial ring in $n$ variables.   Here the
elements $\zeta_1, \dots, \zeta_n$ are in degree 1 and by proper choice
of generators we can assume that ${\rm res}_{G, \langle
x_i\rangle}(\zeta_j) = \delta_{ij} \cdot \gamma_i$ where $\gamma_i \in
\HH{1}{\langle x_i\rangle}{k}$ is a generator for the cohomology ring of
$\langle x_i\rangle$.  Indeed if we assume that the generators are chosen
correctly, then for any $\alpha = (\alpha_1, \dots, \alpha_n) \in k^n$,
$u_{\alpha} = 1 + \sum_{i=1}^n \alpha_i(x_i - 1)\in kG$, $U = \langle
u_{\alpha}\rangle$, we have that
$$
\res{G}{U}(f(\zeta_1, \dots, \zeta_n)) = f(\alpha_1, \dots,
\alpha_n) \gamma^t_{\alpha}
$$
where $f$ is a homogeneous polynomial of degree $t$ and $\gamma_{\alpha}
\in \HH{1}{U}{k}$ is a generator of the cohomology ring of $U$.

Now suppose that $p$ is an odd prime and let $G = \langle x_1, \dots,
x_n\rangle \cong (C_p)^n$. Then
$$
\hgs \cong k[\zeta_1, \dots, \zeta_n] \otimes \Lambda(\eta_1, \dots,
\eta_n)\,,
$$ 
where $\Lambda$ is an exterior algebra generated by the elements
$\eta_1 , \dots, \eta_n$ in degree ~1 and the polynomial generators
$\zeta_1, \dots, \zeta_n$ are in degree ~2. We can assume that each
$\zeta_j$ is the Bockstein of the element $\eta_j$ and that the
elements can be chosen so that ${\rm res}_{G, \langle
x_i\rangle}(\zeta_j) = \delta_{ij} \cdot \gamma_i$ where $\gamma_i \in
\HH{2}{\langle x_i\rangle}{k}$ is a generator for the cohomology ring of
$\langle x_i\rangle$.  Similarly,  assuming that the generators are chosen
correctly, for any $\alpha = (\alpha_1, \dots, \alpha_n) \in k^n$,
$u_{\alpha} = 1 + \sum_{i=1}^n \alpha_i(x_i - 1)\in kG$, $U = \langle
u_{\alpha}\rangle$, we have that
$$
\res{G}{U}(f(\zeta_1, \dots, \zeta_n)) = f(\alpha^p_1, \dots,
\alpha^p_n) \gamma^t_{\alpha}
$$
where $f$ is a homogeneous polynomial of degree $t$ and $\gamma_{\alpha}
\in \HH{1}{U}{k}$ is a generator of the cohomology ring of $U$.

Associated to a $kG$-module $M$ we can define a rank variety
$$
V^r_G(M) = \left\{\alpha \in k^n \mid M{\downarrow}_{\langle u_{\alpha}
\rangle} {\rm \ is \ not \ a \ free \ }\langle u_{\alpha}\rangle\hbox{-module  }\right\} \cup \{0\}
$$
where $u_{\alpha}$ is given as above and where $M{\downarrow}_{\langle
u_{\alpha}\rangle}$ denotes the restriction of $M$ to the subalgebra
$k {\langle u_{\alpha} \rangle}$ of $kG$.
Then we have the following result for any $p$.

\begin{thm}
Let $M$ be any $kG$-module. If $p=2$ then{\rm ,} $V^r_G(M) = V_G(M)$ as
subsets of $k^n$.  If $p  > 2$ then the map $V_G(M) \Rarr{} V^r_G(M)$
given by  $\alpha \mapsto \alpha ^p = (\alpha_1^p, \dots, \alpha_n^p)$
is an inseparable isogeny \/{\rm (}\/both injective and surjective\/{\rm ).}\/ In
particular{\rm ,} for $\alpha \ne 0${\rm ,} $\alpha^p  \in V_G(M)$
{\rm (}$\alpha \in V_G(M)$ if $p = 2${\rm )} if and only if
$M{\downarrow}_{\langle u_{\alpha}\rangle}$ is not a free $k\langle
u_{\alpha}\rangle$-module.
\end{thm}

We should emphasize that if $v$ is a unit in $kG$ such that $$v \equiv
u_{\alpha} \ {\rm mod}({\rm Rad}(kG)^2)$$ then $M{\downarrow}_{\langle v \rangle}$
is a free $k \langle v \rangle$-module if and only if $\alpha^p \not \in
V_G(M)$ ($\alpha \not \in V_G(M)$ if $p = 2$).
So for example the element $x_1x_2x_3$ fails to act freely on
$M$ if and only if $(1, 1, 1, 0, \dots, 0) \in V_G(M)$.

\section{Extraspecial groups in characteristic 2}

In this section and the next, we are interested in the structure  and
cohomology of extraspecial and almost extraspecial $p$-groups.  These are
precisely the $p$-groups $G$ with the property that $G$ has a {\it unique}
normal subgroup $Z$ of order $p$ such that $G/Z$ is elementary abelian.
Note that the dihedral group  $D_8$ of order 8 and, more generally, the
Sylow $p$-subgroup of ${\rm GL}(3,p)$ are extraspecial $p$-groups.
The quaternion group $Q_8$ of order 8 and the cyclic group $C_{p^2}$ of
order $p^2$ also have the required property. Indeed, for $p = 2$ any
extraspecial or almost extraspecial group is constructed from copies of
$D_8, Q_8$ and $C_4$ by taking central products. In this section we
concentrate on the case $p=2$ and look more deeply into the structure
of the extraspecial and almost extraspecial group and their cohomology.

Suppose that $G_1$ and $G_2$ are 2-groups with
the property that each has a unique
normal subgroup of order 2.  Let
$\langle z_i\rangle \in G_i$ be the subgroups.
Then the central product $G_1 * G_2$ is defined by
$$
G_1 * G_2 = (G_1 \times G_2)/\langle (z_1, z_2)\rangle.
$$
It is not difficult to check that
$D_8 * D_8 \cong Q_8 * Q_8$ and that $D_8 * C_4
\cong Q_8 * C_4$.  Moreover, $C_4 * C_4$ has
a central elementary abelian subgroup
of order 4 and hence is not of interest to us
(it is neither extraspecial nor almost extraspecial).
We are left with three types. They are:

\underbar{Type 1.}  $G = D_8 * D_8 * \dots * D_8$
of order $2^{2n+1}$ where $n$ is
the number of factors in the central product.

\underbar{Type 2.}  $G = D_8 * \dots * D_8 * Q_8$
of order $2^{2n+1}$ where $n$ is
the number of factors in the central product.

\underbar{Type 3.}  $G = D_8 * \dots * D_8* C_4$
of order $2^{2n+2}$ where
$n$ is the number of factors isomorphic to $D_8$.

The groups of type 1 and 2 are the extraspecial
groups (see \cite{Gor1}) while
the groups of type 3 are what we call the
almost extraspecial groups.

The groups are also characterized by an
associated quadratic form in the following way. Each
group is a central extension
$$
0 \longrightarrow Z \longrightarrow G \stackrel{\mu}{\longrightarrow} E
\longrightarrow 0
$$
where $Z = \langle z\rangle$ is the unique central
normal subgroup of order 2 and $E \cong {\Bbb F}_2^m$ is elementary
abelian. Recall that a quadratic form on $E$ (as a vector space
over~${\Bbb F}_2)$ is a map $q: E \longrightarrow {\Bbb F}_2$ with the
property that
$$
q(x+y) = q(x) + q(y) + b(x,y)
$$
where $b: E \times E \longrightarrow {\Bbb F}_2$
is a symmetric bilinear form.  Here the
quadratic form~$q$ expresses the class of the extension
as given in the above sequence.
That is, if $\tilde{x}$, $\tilde{y}$ are elements of $G$ and if
$\mu(\tilde{x}) = x$ and $\mu(\tilde{y}) = y$, then
$$
\tilde{x}^2 = z^{q(x)} {\rm \ and \ } [\tilde{x}, \tilde{y}] = z^{b(x,y)}.
$$
Notice here that we are writing the operation
in $G$ as multiplication.
Given the structure of the groups, it is
not difficult to write down the associated quadratic forms.
With respect to a choice of basis, $E$ can be identified with
${\Bbb F}_2^m$ and in the sequel we make this identification. Thus
we write $x=(x_1, \dots, x_m)$ for the elements of~$E$.

\begin{lemma}\label{quadforms}  Let $G$ be an extraspecial or almost
extraspecial group of\break order~$2^{m+1}$.  Then the quadratic form $q$
associated  to $G$ is given on $x =\break (x_1, \dots, x_m)\in {\Bbb F}_2^m = E$
as
follows.
\begin{itemize}\item[]
For type {\rm 1,} $\;q(x) = x_1 x_2 + \dots + x_{2n-1} x_{2n} \quad (m = 2n)$.
\item[]
For type {\rm 2,} $\;q(x) = x_1 x_2 + \dots + x_{2n-3} x_{2n-2} + x^2_{2n-1} +
x_{2n-1} x_{2n}$\hfill\break\phantom{For type 2, $q(x)=$} $+ x^2_{2n} \quad  (m = 2n)$.
\item[]
For type {\rm 3,} $\;q(x) = x_1 x_2 + \dots + x_{2n-1}x_{2n} + x^2_{2n+1}
 \quad  (m = 2n+1)$. \end{itemize}
\end{lemma}

Now on the $k$-vector space $V=k^m$ of  dimension $m$, let $q, b$ denote
the same forms but with the field of coefficients  expanded from ${\Bbb
F}_2$
to~$k$. Let $F: k \rightarrow k$ be the Frobenius  homomorphism, $F(a) =
a^2$.  If $\nu = (x_1, \dots, x_m) \in V$, let $F$ act on
$\nu$ by $F(\nu) = (x^2_1, x^2_2, \dots, x^2_m)$.  Recall that a subspace
$W \subseteq V$  is isotropic if $q(w) = 0$ for all $w \in W$. The following
is
not difficult:

\begin{lemma}\label{codim} Let $h$ be the codimension  in $V$ of a
maximal isotropic subspace of~$V$. The values of $h$ for the
quadratic forms associated to the above groups are\/{\rm :}\/
\begin{itemize}
\item[]
$h=n$ for $G$ of type {\rm 1 (}$m = 2n${\rm ),}
\item[]
$h=n+1$ for $G$ of type {\rm 2 (}$m = 2n${\rm )} or type {\rm 3 ($m = 2n+1$).}
\end{itemize}
Moreover $2^h$ is the index in $G$ of a maximal elementary
abelian subgroup.
\end{lemma}

We are now prepared to state the theorem of Quillen
on the cohomology.  See \cite{BCexp} for one treatment.

\begin{thm}[\cite{Qu2}]\label{coho-exp}  Let $G$ be an extraspecial or
almost extraspecial group of order~$2^{m+1}$. If  
$\nu = (x_1, \dots, x_m)${\rm ,} then  
$$
\hgs = k[x_1, \dots, x_m]/(q(\nu), b(\nu, F(\nu)),
\dots, b(\nu, F^{h-1}(\nu))) \;\otimes k[\delta]
$$
where $\delta$ is an element of degree $2^h$
that restricts to a nonzero element of~$Z$.
Moreover the elements $q(\nu), b(\nu, F(\nu)),
\dots, b(\nu, F^{h-1}(\nu))$
form a regular sequence in $k[x_1, \dots, x_m]$
and $\hgs$ is a Cohen-Macaulay ring.
\end{thm}

The following will be vital for the proof of our main results.

\begin{thm}\label{noserrelmt}  Let $G$ be an extraspecial or
almost extraspecial $2$-group. Define $t =t_G$ to be the natural
number given as follows. If $G$ is of type $1$ of order $2^{2n+1}${\rm ,} let
$$
t_G = \begin{cases} 2^{n-1} + 1 &{\it for} \ n \le 4\,,\\
2^{n-1} + 2^{n-4} &{\it for } \ n \ge 4\,.
\end{cases}
$$
If $G$ is of type $2$ of order $2^{2n+1}$ or of type $3$ of
order $2^{2n+2}${\rm ,} then let
$$
t_G = \begin{cases} 3 & {\it for } \ n = 1\,,\\
2^n + 2^{n-2} &{\it for }\ n \ge 2\,.
\end{cases}
$$
Then there exist  nonzero elements
$\zeta_1, \dots, \zeta_t \in \HH{1}{G}{{\Bbb F}_2}$ such that $\zeta_1
\dots
\zeta_t\break = 0$. Moreover{\rm ,} in the isomorphism $\HH{1}{G}{{\Bbb F}_2} \cong
\Hom(G,{\Bbb F}_2)${\rm ,} each $\zeta_i$ corresponds to a homomorphism whose
kernel is a maximal subgroup of $G$ and is the centralizer of a
noncentral involution in $G$.
\end{thm}

\Proof   The proof is contained in the paper \cite{Yal}.
For the groups of type 1, $t_G$ is actually equal to the
cohomological length, that is, the least
number of nonzero elements in $\HH{1}{G}{{\Bbb F}_2}$
such that the product of those elements is zero (see \cite[Th.~1.3]{Yal}).

Now, suppose that $G$ has type 2 or 3.  Then $t_G$ in our
theorem is equal to the cardinality $s(G)$ of a representing
set in $G$ (see \cite[Props.~6.2 and 6.3]{Yal}).
A representing set for $G$ is a collection of elements of
$G$ that contains at least one noncentral element
from each elementary abelian subgroup of $G$.  But
now Proposition 1.1 of \cite{Yal} shows that the
cohomological length is at most~$s(G)$.

The point of the last statement is that the centralizer
of any maximal elementary abelian subgroup of $G$ is
contained in the centralizers of some elements in a
representing set. Because the cohomology ring $\HH{*}{G}{{\Bbb F}_2}$
is Cohen-Macaulay (see Theorem~\ref{coho-exp}), any element whose
restriction to  the centralizer of every maximal elementary abelian subgroup
of $G$ vanishes, is the zero element (see Theorem~3.4 in~\cite{Yal}).
Hence if we choose the elements $\zeta_i$ to correspond to the
centralizers of the elements in a representing set as in the last
statement, then their product is zero as desired.
\Endproof\vskip4pt  

The next theorem will be very important to the proof of
the general case. It is part of the effort to get an
explicit upper bound on the dimensions of critical modules.

\begin{thm} \label{dimcoho}
Let $G$ be an extraspecial or almost extraspecial group of order~$2^{m+3}$
and let $H$ be the centralizer of a noncentral involution in $G$.
Then $H \cong C_2 \times U$ where $U$ is an extraspecial
or almost extraspecial group of order~$2^{m+1}$ of the same type as~$G$.
Assume that $m \geq 2$ and{\rm ,} if $m=2${\rm ,} that $U\not\cong D_8$.
Then for $2 \leq r \leq t_G${\rm ,}
$$
\Dim \HH{r}{H}{k} \leq \binom{m+r}{r} - \binom{m+r-2}{r-2}.
$$
\end{thm}

\Proof 
The structure of the centralizer $H$ can be verified directly
from what we know of $G$. For one thing it can be checked that
all noncentral involutions in $G$ are conjugate by an element
in the automorphism group of $G$ and hence their centralizers are all
isomorphic. 

Throughout the proof we use the notation in Theorem~\ref{coho-exp},
for the cohomology of the group $U$, so that
$\HH{*}{U}{k}$ is generated by $x_1, \dots, x_m$ and $\delta$, with
$\text{deg}(\delta)=2^h$  (where $h$ is the value associated to the
group~$U$ as in Lemma~\ref{codim}). We know that
$$
\HH{*}{H}{k} \cong \HH{*}{U}{k} \otimes \HH{*}{C_2}{k}
$$
and moreover we know that $\HH{*}{C_2}{k} \cong k[y]$
is a polynomial ring in one variable $y$ in degree ~1.
We want to focus on the polynomial ring $S$ generated by
$x_1, \dots, x_m, y$. We have a homomorphism from $S$
to $\HH{*}{H}{k}$ whose kernel contains the elements
$q(\nu)$ and $\beta(\nu,F(\nu))$ where
$\nu = (x_1, \dots, x_m)$. Let $Q$ denote the image
of $S$ in $\HH{*}{H}{k}$.  For this argument, let
$S^{\#} = S/(q(\nu))$ and let $S^{\#\#} = S/(q(\nu),\beta(\nu,F(\nu)))$.
If $R$ denotes any of these graded rings, we let $R_r$ denote the
homogeneous part of $R$ in degree exactly~$r$.  Note that $R_r=0$
if $r<0$.

First notice that $\Dim S_r = \binom{m+r}{m}= \binom{m+r}{r}$.
Because $q(\nu)$ and $\beta(\nu,F(\nu))$ are two terms of
a  regular sequence of elements in $S$ we must have that
$$
\Dim S^{\#}_r = \Dim S_r - \Dim S_{r-2}
$$
and 
$$
\Dim S^{\#\#}_r = \Dim S^{\#}_r - \Dim S^{\#}_{r-3}
$$
for all $r \geq 2$. Moreover $\Dim S_r \geq \Dim S^{\#\#}_r
\geq \Dim Q_r$ for all values of $r$.

By Theorem~\ref{noserrelmt}, $t_G \leq 2t_U$ (with equality in most
cases) and moreover, by Lemma~\ref{codim}, we see that $t_U<2^h$ in all
cases.
 The choice that $r \leq  t_G$ now means  that
 $$
r \leq  t_G \leq 2 t_U< 2 \cdot 2^h = 2 \cdot \text{deg}(\delta)
$$
and this implies that we must have either $\Dim \HH{r}{H}{k} = \Dim Q_r$ ,
if $r < \text{deg}(\delta)$, or $\Dim \HH{r}{H}{k} =
\Dim Q_r + \Dim (\delta \cdot Q_{r-\text{deg}(\delta)})$,
if $\text{deg}(\delta) \leq r < 2 \text{deg}(\delta)$.
Notice also that $\text{deg}(\delta) = 2^h \geq 4$
in all cases because we assumed that $m \geq 2$ and
$U \not\cong D_8$ (if $U \cong D_8$, then
$h=1$ and $\text{deg}(\delta) = 2$). Hence we have that
\begin{align*}
\Dim \HH{r}{H}{k} &\leq \Dim Q_r +\Dim Q_{r-\text{deg}(\zeta)}\\
&\leq  \Dim S^{\#}_r - \Dim S^{\#}_{r-3}
+ \Dim S^{\#}_{r-\text{deg}(\zeta)}
- \Dim S^{\#}_{r-\text{deg}(\zeta)-3}\\
&\leq  \Dim S^{\#}_r - \Dim S^{\#}_{r-3}
+ \Dim S^{\#}_{r-\text{deg}(\zeta)}\\
&\leq  \Dim S^{\#}_r = \binom{m+r}{r} - \binom{m+r-2}{r-2}\,.
\end{align*}
The last inequality follows from the facts that $r- \text{deg}(\delta)
\leq r-3$ and that $\Dim S^{\#}_s$ is an increasing function of $s$.
\hfill\qed

\begin{cor} \label{dimsumomeg}
Suppose that $G$ and $H$ are as in the theorem. If $2 \leq r \leq t_G${\rm ,}
then 
$$
\sum_{i = 0}^r \Dim \Omega^i(k_H){\uparrow}_H^G \leq \binom{m+r-1}{m} \vert G \vert + 2.
$$
\end{cor}

\Proof 
For any $i$ we have an exact sequence
$$
0 \Rarr{} \Omega^{i+1}(k_H) \Rarr{} P_i \Rarr{} \Omega^i(k_H) \Rarr{} 0
$$ 
where $P_i$ is the degree $i$ term in a minimal $kH$-projective
resolution of the trivial $kH$-module $k_H$. Recall that $\Dim P_i
= \Dim \HH{i}{H}{k} \cdot \vert H \vert$. Then by the theorem,
for $r = 2s+1$,  
\begin{align*}
\sum_{i=0}^{r} \Dim \Omega^i(k_H) \ &= \
\sum_{j=0}^{s} \big{(}\Dim \Omega^{2j+1}(k_H) +
\Dim \Omega^{2j}(k_H)\big{)}
= \sum_{j=0}^{s} \Dim P_{2j} \\
&\leq \ \Dim P_0 + \sum_{j=1}^s \Big{[} \binom{m+2j}{2j} -
\binom{m+2j-2}{2j-2} \Big{]} \vert H \vert \\
&= \ \vert H\vert + \Big{[}\binom{m+2s}{2s} - \binom{m}{0}\Big{]}
\vert H\vert \\ 
&= \ \binom{m+r-1}{r-1} \vert H\vert = \binom{m+r-1}{m} \vert H\vert
\,.
\end{align*}
On the other hand if $r = 2s$ is even, then we use the fact that
$\Dim P_1 = \binom{m+1}{1}\vert G \vert$ and we obtain similarly
\begin{eqnarray*}
\sum_{i=0}^{r} \Dim \Omega^i(k_H)   &=& \Dim k + \Dim P_1 +
 \sum_{j=2}^{s} \Dim P_{2j-1} \\
&\leq& 1 + \binom{m+2s-1}{2s-1} \vert H\vert
= 1 + \binom{m+r-1}{m} \vert H\vert.
\end{eqnarray*}
In both cases, inducing from $H$ to $G$, the dimension of
$\Omega^i(k_H){\uparrow}_H^G$ is doubled and the result follows.
\Endproof\vskip4pt

\section{Extraspecial groups in odd characteristic}

Our aim in this section is to get results similar to those of the
last section for extraspecial $p$-groups in the case that the prime
$p$ is not~$2$. As in the characteristic~$2$ case, for any positive
integer $n$ there are two isomorphism types of extraspecial groups
of order $p^{2n+1}$ and one isomorphism type of almost extraspecial
group of order $p^{2n+2}$. For each $n$, one of the two nonisomorphic
groups of order $p^{2n+1}$ has exponent $p^2$ and the other one has
exponent~$p$. In the earlier paper~\cite{CT1} we showed that
Theorem~\ref{mainthm} holds for extraspecial groups
of exponent~$p^2$ and almost extraspecial groups (i.e.\ for these groups
there are no nontrivial critical modules). As a consequence, the only
groups of interest to us are the extraspecial groups of order $p^{2n+1}$ and
exponent $p$.

Up to isomorphism, there is exactly one extraspecial group $G_1$ of
order $p^3$ and exponent $p$. It is generated by elements $x$, $y$ and
$z$, which satisfy the relations that $z$ is in the center of $G_1$,
$z^p = x^p = y^p =1$ and $[x,y] = z$. It is isomorphic to the Sylow
$p$-subgroup of the general linear group ${\rm GL}(3,p)$. For $n > 1$,
the extraspecial group of order $p^{2n+1}$ is a central product
$$
G_n = G_1 \ * \ G_1 \ * \ \dots \ * \ G_1
$$ 
of $n$ copies of $G_1$ as in the last section. That is, $G_n$ is the
quotient group obtained by taking the direct product of $n$ copies of
$G_1$ and then identifying the centers (see \cite{Gor1}).  The center
of $G_n$ is a cyclic subgroup $Z = \langle z \rangle$ of order $p$
and $G_n/Z$ is an elementary abelian $p$-group of order $p^{2n}$.

We need an analogue to Theorem \ref{noserrelmt} for our case.

\begin{thm}\label{noserrelmt-p}
For $G = G_1${\rm ,} let $t_G = 2(p+1)${\rm ,}
while for $G = G_n${\rm ,} $n > 1${\rm ,} let $t_G = (p^2+p-1)p^{n-2}$.
Then there exist  nonzero elements
$\eta_1, \dots, \eta_t \in \HH{1}{G}{{\Bbb F}_p}$ such that
$\beta(\eta_1) \dots \beta(\eta_t) = 0$ where $t = t_G$.
Moreover{\rm ,} in the isomorphism $\HH{1}{G}{{\Bbb F}_p} \cong
\Hom(G,{\Bbb F}_p)${\rm ,} each $\eta_i$ corresponds to a homomorphism whose
kernel is a maximal subgroup of $G$ and is the centralizer of a
noncentral element of order $p$ in $G$.
\end{thm}

\Proof 
The proof of the theorem is contained in the paper by Yal\c{c}in
as Theorem~1.2 of \cite{Yal}. In this case the dimension of
$\HH{1}{G,\bfp}$ is the same as that of $\Hom(G,\bfp)$ which is~$2n$.
\Endproof\vskip4pt  

As in the last section we are going to need estimates on the dimensions
of the cohomology groups $\HH{r}{G_n}{k}$ where $k$ is a field of
characteristic~$p$. We begin with the case of the extraspecial
group $G = G_1$ of order $p^3$. Ian Leary \cite{Lear1} has given a complete
description of the cohomology ring $\hgs$ except that he did not fully
compute the Poincar\'e series, which is something that we need. The
calculation
is, of course, implicit in his work, and he did calculate
it in the special case that $p = 3$ \cite{Lear2}. Note that our results
agree with his in that situation.

\begin{thm} \label{poin-p3}
The Poincar{\rm \'{\hskip-6pt\it e}} series for the cohomology ring of the group
$G= G_1$ is given by the rational function
$$
\sum_{n = 0}^{\infty} \Dim \HH{n}{G}{k} \ t^n =
\dfrac{1 + t + 2t^2 + 2t^3 +t^4 +t^5 + \dots +
t^{2p-1}}{(1-t)(1-t^{2p})}\,.
$$
\end{thm}
 
\Proof 
We will not repeat the long list of relations given by Leary (Theorem~6
of \cite{Lear1}). However we will use exactly the notation of that paper
and the interested reader can follow the computation. The strategy is
first to ignore the contribution of the regular element $z$ in degree
$2p$. This element is a nondivisor of zero as it restricts nontrivially
to the center of $G$. We also know that it is regular from the given
relation
and from the fact that it is represented on the $E_2$ of the spectral
sequence, by an element in $E_2^{0,2p}$ which survives to the $E_{\infty}$
page of the spectral sequence. Consequently, the Poincar\'e series
$f(t)$ of $\hgs$ is obtained by multiplying
$1/(1-t^{2p})$ times the Poincar\'e series
of the subalgebra $A$ generated by all of the given generators other than
$z$.

Next we consider the subalgebra $A$ as a module over the subring $R$
generated by $x$ and $x'$. Note that $x$ and $x'$ are in degree ~2 and
satisfy the relation $x^px' - x{x'}^p=0$ in degree $2p+2$. So the Poincar\'e
series for $R$ is $f_1 = (1-t^{2p+2})/(1-t^2)^2$. This is also the series
for 
the $R$-submodule $M_1$ generated by the element $1$ in degree~$0$. The
first thing that needs to be established from the relations is that the
$R$-generators are the elements of the sequence
$$
S= [1,y,y',Y,Y',X,X',yY', XY',XX',d_4,c_4,d_5, \dots, c_{p-1},d_p]
$$
of length $2p+3$.
Let $M_i$ be the $R$-submodule generated by the first $i$ elements of
the sequence, and let $f_i$ be the Poincar\'e series for $M_i/M_{i-1}$.
Then the desired Poincar\'e series for $A$ is $f_1 + f_2 + \dots +f_{2p+3}$.
Note that $f_1$ has been calculated.
\begin{itemize}
\item For $f_2$, we note that $xy' = x'y$
and $x^py' = {x'}^py$. So $x'(x^{p-1}-{x'}^{p-1})y =0$. Therefore $f_2 =
t(1-t^{2p})/(1-t^2)^2$.

\item Since $xy' = x'y \in M_2$, we have that $f_3 = t/(1-t^2)$.

\item Similarly to the calculation for $f_2$, we have that
$f_4 = t^2(1-t^{2p})/(1-t^2)^2$ and $f_6 = t^3(1-t^{2p})/(1-t^2)^2$.

\item For $f_5$, note that $x^2Y' = xx'Y \in M_4$ and $xx'Y' \in M_4$.
Therefore  $f_5 = t^3 + t^2/(1-t^2)$.

\item The calculation for $f_7$ is similar to that for $f_3$ and we get that
$f_7 = {t^3/(1-t^2)}$.

\item For $i := 8, \dots, 2p+3$, it should be checked that $xS_i, x'S_i \in
M_{i-1}$ where $S_i$ is the $i^{th}$ element of the sequence $S$.
Consequently, $f_i = t^{j_i}$, where $j_i$ is the degree of $S_i$. Note
that $j_8 = 3$ while $j_i = i-4$ for $i \geq 9$.
\end{itemize} 

Finally it is necessary to verify that
$$
f_1 + f_2 + \dots + f_{2p+3} = (1+t+2t^2+2t^3+t^4 + \dots +
t^{2p-1})/(1-t)
$$
by routine but tedious calculation.
\Endproof\vskip4pt  

We need to derive two facts from the above theorem. The first is an upper
bound which is not optimal but will be sufficient for our purposes.

\begin{cor} \label{dimcoho-p3}
For $G = G_1${\rm ,} 
$$
\Dim \HH{r}{G}{k} \leq 2(r+1) = 2 \binom{r+1}{1}.
$$
Moreover{\rm ,} $\Dim \HH{r}{G}{k}=2r$ if $1\leq r \leq 3$ and $\Dim
\HH{r}{G}{k}=r+3$ if $4\leq r \leq 2p-1$.
\end{cor}

\Proof  Consider the series expansion
$$
g(t) = \dfrac{1+t+2t^2+2t^3+t^4 + \dots +t^{2p-1}}{1-t} =
\sum_{r=0}^{\infty} a_r t^r \,.
$$
A routine computation yields the value of the coefficients
$a_0 = 1$, $a_r= 2r$ if $1\leq r \leq 3$, $a_r = r+3$ if $4\leq r \leq
2p-1$, and $a_r = 2p+2$ if $r \geq 2p-1$. The Poincar\'e series for the
cohomology ring of~$G_1$ is obtained by multiplying $g(t)$ with
${1\over{1-t^{2p}}}= \sum_{i=0}^{\infty}  t^{2ip}$. Therefore $\Dim
\HH{r}{G}{k} = a_r$ for $r\leq 2p-1$ and this proves the second statement of
the lemma. Moreover, for arbitrary~$r$, writing $r = j + q(2p)$ with $0 \leq
j < 2p$, we have that
\vskip12pt
\hfill $
\displaystyle{\Dim \HH{r}{G}{k} = a_j + qa_{2p} \leq (j+3) + q(2p+2) \leq 2(r+1)\,.}
$ 
\Endproof\vskip4pt  

\begin{cor} \label{dim-omega-2p}
For $G = G_1${\rm ,} $\Dim \Omega^{2p}(k) = p^3(p+1)+1$.
\end{cor}

\Proof 
If $P_j$ is the $j$-th term of a minimal projective resolution of $k$, we
have
$\Dim(P_j) = \Dim \HH{j}{G}{k} \, \vert G\vert$ and so
$\Dim \Omega^{j+1}(k) = \Dim \HH{j}{G}{k} \vert G\vert - \Dim
\Omega^{j}(k)$.
Using this relation and the dimensions given in the previous corollary, we
obtain $\Dim \Omega^{2}(k)= p^3+1$ and then by induction $\Dim
\Omega^{2j-1}(k) = (j+1)p^3-1$ and $\Dim \Omega^{2j}(k) = (j+1)p^3+1$ for
$2\leq j \leq p$.
\Endproof\vskip4pt  

In the rest of the section, we require the following
well known combinatorial identity.

\begin{lemma} \label{binomident}
For all integers  $c,i,j \geq 0${\rm ,}
$$
\sum_{a+b = c} \binom{a+i}{i}\binom{b+j}{j} = \binom{c+i+j+1}{i+j+1}.
$$
\end{lemma}

\Proof 
Recall that if $P$ is a polynomial ring in $n$ variables,
then the number of monomials of degree $r$ is $\binom{r+n-1}{n-1}$. Now the
tensor product of a polynomial ring in $i+1$ variables with a polynomial
ring in $j+1$ variables yields a polynomial ring in $i+j+2$ variables. The
identity follows by counting the number of monomials of degree~$c$.
\Endproof\vskip4pt  

We also need to know the dimension of the cohomology groups of elementary
abelian groups.

\begin{lemma} \label{dimcoho-elem}
Let $p$ be an odd prime and let $E$ be an elementary abelian $p$-group of
rank~$m$. Then $ \Dim \HH{r}{E}{k} =  \binom{r+m-1}{m-1}$.
\end{lemma}

\Proof 
Recall that $\HH{*}{E}{k} \cong k[\zeta_1, \dots, \zeta_m] \otimes
\Lambda(\eta_1, \dots, \eta_m)$ where $\zeta_1, \dots, \zeta_m$ are in
degree~2 and $\eta_1, \dots, \eta_m$ are in degree~1. A basis of
$\HH{r}{E}{k}$ consists of the elements $\zeta_1^{a_1}
\dots,\zeta_m^{a_m}\eta_1^{e_1}, \dots, \eta_m^{e_m}$ where $0\leq a_i\leq
r/2$, $0\leq e_i\leq 1$ and ${\sum_{i=1}^m (2a_i+e_i)} = r$.
This basis is in bijection with the set of monomials of degree~$r$ in
$k[x_1,\dots,x_m]$ by mapping the above basis element to
$x_1^{2a_1+e_1}\dots x_m^{2a_m+e_m}$. Now the number of monomials of degree
$r$ is $\binom{r+m-1}{m-1}$. \phantom{odd}
\Endproof\vskip4pt  

Our main result in this section gives estimates for the dimensions of the
cohomology of the centralizers of $p$-elements.

\begin{thm} \label{dimcoho-p}
Let $G = G_n$ be an extraspecial group of order~$p^{2n+1}$ and exponent $p$.
Let $H$ be the centralizer of a noncentral element of order $p$ in $G$.
Then $H \cong C_p \times G_{n-1}$.
Moreover{\rm ,} 
$$
\Dim \HH{m}{H}{k} \leq 2\binom{m+2n-2}{2n-2}.
$$
\end{thm}
 
\Proof 
As with the characteristic ~2 case, the structure
of the centralizer $H$ can be verified directly
from what we know of $G$. All noncentral elements of order $p$ in $G$
are conjugate by an element in the automorphism group of $G$ and hence
their centralizers are isomorphic.

Next we need to approximate the dimensions of the cohomology groups
of the group $G_{n-1}$ for $n \geq 1$. The estimate in Corollary
\ref{dimcoho-p3} will serve in the case that $n=2$. Let $N$ be a  normal
subgroup of $G_{n-1}$ such that $N \cong G_1$. We can take $N$ to be the
first factor in the central product that expresses
$G_{n-1}$. Then $G_{n-1}/N \cong C_p^{2(n-2)}$, an elementary abelian
group of order $p^{2(n-2)}$. The Lyndon-Hochschild-Serre
spectral sequence of the extension of $G_{n-1}/N$ by $N$ has $E_2$ term
$$
E_2^{r,s} = \HH{r}{G_{n-1}/N}{\HH{s}{N}{k}}  \Rightarrow
\HH{r+s}{G_{n-1}}{k}.
$$
As $k$-vector spaces, it is true that $E_2^{r,s} \cong \HH{r}{G_{n-1}/N}{k}
\otimes
\HH{s}{N}{k}$ because $N$ commutes with the other factors of the central
product. So we have that
\begin{align*}
\Dim \HH{m}{G_{n-1}}{k} & \leq \sum_{r+s=m}\Dim(E_2^{r,s}) \\
& =  \sum_{r+s=m} \Dim \HH{r}{G_{n-1}/N}{k}  \ \Dim \HH{s}{N}{k} \\
& \leq \sum_{r+s=m} \binom{r+2(n-2)-1}{2(n-2)-1} 2 \binom{s+1}{1}
= 2\binom{m+2n-3}{2n-3},
\end{align*}
using Lemma~\ref{dimcoho-elem}, Corollary~\ref{dimcoho-p3} and the
combinatorial identity of Lemma~\ref{binomident}.

Now $\HH{m}{H}{k} \cong \bigoplus_{r+s=m} \HH{r}{G_{n-1}}{k} \otimes
\HH{s}{C_p}{k}$. Therefore,
\begin{align*} 
\Dim \HH{m}{H}{k} &=  \sum_{r+s=m} \Dim \HH{r}{G_{n-1}}{k} \cdot
\Dim \HH{s}{C_p}{k} \\
&\leq \sum_{r+s=m} 2\binom{r+2n-3}{2n-3}\binom{s}{0} =
2\binom{m+2n-2}{2n-2},
\end{align*}
again by Corollary~\ref{dimcoho-p3} and Lemma~\ref{binomident}.
\hfill\qed

\begin{cor}\label{dimsumomeg-p}
Suppose that $G$ and $H$ are as in the theorem. If $r \geq 1${\rm ,}
then
$$
\sum_{i = 0}^r \Dim \Omega^i(k_H){\uparrow}_H^G \leq 2p^{2n+1} \binom{r+2n-2}{2n-1}.
$$
\end{cor}

\Proof 
Suppose that $\quad \dots \rightarrow P_1 \rightarrow
P_0 \rightarrow k \rightarrow 0$ is a minimal $kH$-projective
resolution of the trivial module $k$. Then we know that
$\Dim \Omega^0(k) + \Dim \Omega^1(k) = \Dim P_0$. For $j \geq 2$,
$\Omega^j(k_H)$ is a submodule of $P_{j-1}$. The dimension of $P_j$ is
precisely $\vert H\vert \Dim \HH{j}{H}{k}$ and the dimension
of $\Omega^j(k_H){\uparrow}_H^G$ is $p$ times the dimension of
$\Omega^j(k_H)$. So from the theorem we have that
\begin{align*}
\sum_{i = 0}^r \Dim \Omega^i(k_H){\uparrow}_H^G & \leq
p \vert H\vert \sum_{i=0}^{r-1} \Dim \HH{i}{H}{k} \\
& \leq p^{2n+1} \sum_{i=0}^{r-1} 2\binom{i+2n-2}{2n-2}\binom{r-1-i}{0} \\
& = 2p^{2n+1} \binom{r+2n-2}{2n-1},
\end{align*}
by the identity \ref{binomident}.
\hfill\qed

\section{New endo-trivial modules from old endo-trivial modules}

Here we start the proof of Theorem~\ref{mainthm}. Suppose that $G$ is an
extraspecial or almost extraspecial $p$-group and that $G\not\cong Q_8$.
Let $Z = \langle z \rangle$ be the Frattini subgroup of $G$, of order $p$,
with elementary abelian quotient $\overline{G} = G/Z$ of rank~$m$. Let
$x_1, \dots, x_m \in G$ such that $\overline{G} =
\langle \overline{x}_1, \dots, \overline{x}_m\rangle$.
Recall that $Z$ is the unique normal subgroup of order~$p$. Moreover
every maximal subgroup of $G$
contains~$Z$ and $G$ is not elementary abelian.
Some of the results in this section hold more generally if $G$ has a
Frattini subgroup $Z$ of order $p$, but we leave this generalization to the
reader.

Let $M$ be an endo-trivial $kG$-module whose class in $T(G)$ lies in the
kernel of the restriction to proper subgroups. This means that
$M{\downarrow}_H^G \cong k \oplus ({\rm free})$ for every maximal subgroup
$H$ of~$G$. For the purpose of the proof of Theorem~\ref{mainthm}
(Sections 5--11), we make the following definition:

\begin{defi}\label{critical}
We say that a $kG$-module $M$ is {\it critical\/}
if it is an indecomposable endo-trivial module such that
$M{\downarrow}_H^G \cong k \oplus ({\rm free})$ for every maximal subgroup
$H$ of~$G$. 
\end{defi}

Actually, the last condition implies that the module $M$ is
endo-trivial because its restriction to every elementary abelian subgroup
is isomorphic to $k \oplus ({\rm free})$, hence endo-trivial (see
Lemma~2.9 of~\cite{CT1}). In fact $M$ is a torsion endo-trivial module by
a theorem of Puig~\cite{Pu1}, but we do not need this fact in our
arguments. By factoring out all free summands of an endo-trivial module
$M$, one can always assume that $M$ is indecomposable and this is why we
do so. We shall often omit to mention this indecomposability condition, to
the effect that we shall usually only prove that a module satisfies the
condition on restriction to maximal subgroups in order to deduce that it
is critical. Since our aim is to prove that the kernel above is trivial,
we have to show that any critical $kG$-module $M$ is isomorphic to $k$ as
a $kG$-module. We will often assume, by contradiction, the existence of a
nontrivial critical $kG$-module.

In this section, we prove several results concerning the structure of
a critical module~$M$ and the construction of new modules with the same
property. For some of the results, we only need to assume that
$M{\downarrow}_H^G \cong k \oplus ({\rm free})$ for a single
subgroup $H$ of~$G$.

For any critical $kG$-module $M$, and more generally for any $kG$-module
$M$ such that $M{\downarrow}_Z^G \cong k \oplus ({\rm free})$, we let
$M' =\{m \in M \;|\;(z - 1)^{p-1}m = 0\}$ and we set
$$\overline{M} = M/M'\,.$$
We let $^-: M \longrightarrow \overline{M}\,$ be the
quotient map. Since $(z-1)\overline{M}=0$, the module $\overline{M}$ can be
viewed as a $k\overline{G}$-module. A large part of this paper is devoted
to an analysis of the properties of the module~$\overline{M}$.

\begin{lemma} \label{struc-critical} Let $M$ be a $kG$-module. Suppose
that $M{\downarrow}_Z^G \cong k \oplus ({\rm free})$.
\begin{itemize}
\ritem{(a)} The module $M$ has two filtrations
$$
\begin{array}{ccccccccccc}
 &  & K_1 & \subset & K_2 & \subset & \dots & \subset & K_{p-1} &
 \subset & K_p=M \\
 &  & \cup &  & \cup & & & & \cup & & \\
 \{0\} & \subset & I_{p-1} & \subset & I_{p-2} & \subset & \dots & \subset
& I_1 & &
\end{array}
$$
where $K_i = \{ m \in M \;|\;(z - 1)^i m = 0\}$ is the kernel of
multiplication by $(z-1)^i$ \/{\rm (}\/in particular $K_{p-1} = M'$\/{\rm )}\/ and $I_i =
(z-1)^iM$ is the image of multiplication by $(z-1)^i$.

\ritem{(b)} $K_i/I_{p-i} \cong k$ for any $i = 1, \dots, p{-}1$. Moreover
$K_{p-1}/I_{p-1} \cong k \oplus (I_1/I_{p-1})$.

\ritem{(c)} The module $I_1 = (z-1)M$ is free as a module over the ring
$kZ/(z-1)^{p-1}$. Moreover\/{\rm ,}\/ $I_i/I_{i+1} \cong \overline{M}$
for any $i = 1, \dots, p{-}1$.

\ritem{(d)} The module $M/K_1$ is isomorphic to~$I_1$. In particular it is free as
a module over the ring $kZ/(z-1)^{p-1}$ and $K_{i+1}/K_{i} \cong
\overline{M}$ for any $i = 1, \dots, p{-}1$.

\ritem{(e)} $\Dim(M) = p \Dim(\overline{M}) + 1$.
\end{itemize}
\end{lemma}

\Proof  (a) Note that $K_i$ and $I_i$ are submodules
because $z$ is central in~$kG$. We have $I_{p-i} \subset K_i$ because
$(z-1)^p =0$. The filtrations are clear.
\smallbreak
(b) In order to prove (b), it suffices to restrict to the subgroup~$Z$.
But we have $M{\downarrow}_Z^G = k \oplus F$ for some free $kZ$-module~$F$,
and therefore
$$K_i = k \oplus (z-1)^{p-i}F\,, \quad I_{p-i} = (z-1)^{p-i}F\,.$$
Moreover it is clear that $K_{p-1}/I_{p-1} = K_1/I_{p-1} \oplus
(I_1/I_{p-1}) \cong k \oplus (I_1/I_{p-1})$.

\smallbreak (c) Multiplication by $(z-1)^i$ induces a map
$$ M \Rarr{} (z-1)^iM/(z-1)^{i+1}M = I_i/I_{i+1} $$
and we claim that its kernel is $M'$. Again, in order to prove this, it
suffices to restrict to the subgroup~$Z$ and consider the decomposition
$M{\downarrow}_Z^G = k \oplus F$ as above. Then the kernel is
$k \oplus (z-1)F = M'$. It is also clear that
$$ (z-1)M = (z-1)F \cong F/(z-1)^{p-1}F$$
and this is free over the ring $kZ/(z-1)^{p-1}$.

\smallbreak (d) Multiplication by $(z-1)$ induces an isomorphism $M/K_1 \cong I_1$.

\smallbreak (e) Since $M{\downarrow}_Z^G = k \oplus F$, we have that 
$\Dim(\overline{M}) = \Dim(F/(z{-}1)F) = \Dim(F)/p$ and
$\Dim(M) = p\Dim(\overline{M}) + 1$.
\hfill\qed

\begin{lemma} \label{barperiod} Let $M$ be a $kG$-module. Suppose that
there is a maximal subgroup $H$ of $G$ such that
$M{\downarrow}_H^G \cong k \oplus ({\rm free})$.
\begin{itemize}
\ritem{(a)} $M\cong k\oplus ({\rm free})$ as a $kG$-module if and only if
$\overline{M}$ is a free $k\overline{G}$-module. More precisely{\rm ,} $M$ has a
free summand with $r$ generators as a $kG$-module if and only if
$\overline{M}$ has a free summand with $r$ generators as a
$k\overline{G}$-module. In particular{\rm ,} if $M$ is indecomposable{\rm ,} then
$\overline{M}$ has no projective summands.

\ritem{(b)} $M\not\cong k\oplus ({\rm free})$ as a $kG$-module if and only if
$\overline{M}$ is a periodic $k\overline{G}$-module.
\end{itemize}
\end{lemma}

\Proof  (a) It is easy to see that if $M$ has a free
summand $L \cong (kG)^r$ as a $kG$-module then $\overline{M}$ has a free
summand $L/(z-1)L\cong (k\overline{G})^r$ as a
$k\overline{G}$-module.

The converse is essentially contained in
Lemma~3.3 of~\cite{CT1} and we recall the argument.
Assume that $M = N \oplus L$ where $L$ is free and $N$ has no free
summands. Then $t_1^G \cdot N = 0$ where
$$t_1^G = \sum_{g\in G} g =
(z-1)^{p-1} \displaystyle\prod^m_{i=1}(x_i-1)^{p-1}\,,$$
 $x_i$ being a lift in $G$ of the generator $\overline{x}_i$
of~$\overline{G}$.
 Let $X = \displaystyle\prod^m_{i=1}(x_i-1)^{p-1}$. If $\overline{N}$
has a free submodule then $X \cdot \overline{N} \ne 0$,
since $\overline{X} = t_1^{\overline{G}}$.  But if $X \cdot
\overline{N} \ne 0$ then, via the isomorphism $\overline{N} \cong
(z-1)^{p-1}N$ of Lemma~\ref{struc-critical}, we would obtain
$(z-1)^{p-1}X  \cdot N = t_1^G \cdot N \ne 0$, which is a contradiction.

(b) The hypothesis on $M{\downarrow}_H^G$ implies that $\overline{M}$ is
free on restriction to $H/Z$.  But $\overline{H}=H/Z$ is a maximal
subgroup of $\overline{G} = G/Z$, so $\overline{G}/\overline{H}$ is a
cyclic group of order~$p$. Tensoring with $\overline{M}$ the exact sequence
$$ 0 \Rarr{} k \Rarr{} k[\overline{G}/\overline{H}] \Rarr{}
k[\overline{G}/\overline{H}] \Rarr{} k \Rarr{} 0
\,,$$
we obtain an exact sequence with $\overline{M}$ at both ends and free
$k\overline{G}$-modules in the middle, because
$k[\overline{G}/\overline{H}] \otimes \overline{M} \cong
\overline{M}{\downarrow}_{\overline{H}}^{\overline{G}}
{\uparrow}_{\overline{H}}^{\overline{G}}$.
If now $M \not\cong k\oplus ({\rm free})$, then $\overline{M}$ is not zero
and is not free as a $k\overline{G}$-module, by part (a), so
$\overline{M}$ is periodic. If conversely $\overline{M}$ is periodic,
then $\overline{M}$ is not free and $M \not\cong k\oplus ({\rm free})$ by
part (a).
\phantom{overthere}\hfill\qed

\begin{lemma} \label{struc-two} Suppose that $p=2$ and that $M$ is a
nontrivial critical\break $kG$-module. Then the number of
generators of $M$ is the same as the number of generators
of $\overline{M}$ and is equal to $4\Dim(\overline{M})/\vert G \vert$.
Moreover $\Dim(\Omega(M)) = \Dim(\Omega^{-1}(M)) = \Dim(M) -2$.
\end{lemma}

\Proof  Let $H$ be a maximal subgroup of~$G$. Since
$M{\downarrow}_H^G\cong k \oplus ({\rm free})$, we know that
$\overline{M}$ is free as a module over~$k\overline{H}$. Thus, the number of
generators of
$\overline{M}$ as a $k\overline{H}$-module is $\Dim(\overline{M})/\vert
\overline{H} \vert$. Our first claim is that $G$ acts trivially on
$\overline{M}/{\rm Rad}(k\overline{H})\overline{M}$. Thus, the number of
generators of $\overline{M}$ as a $k\overline{G}$-module is also
$\Dim(\overline{M}/{\rm Rad}(k\overline{H})\overline{M}) =
\Dim(\overline{M})/\vert\overline{H} \vert$. In order to prove the claim,
we note that the group $G/H$ acts on
$\overline{M}/{\rm Rad}(k\overline{H})\overline{M}$. If there were a free
summand generated by the class of an element $\overline{m}$, then
$\overline{m}$ would generate a free summand of $\overline{M}$ as a module
over $k\overline{G}$, contrary to part~(c) of the previous lemma. Since the
group $G/H$ has order~2, the only possibility is that $G/H$ acts trivially
on $\overline{M}/{\rm Rad}(k\overline{H})\overline{M}$.

Now our second claim is that, given a set of generators of
$\overline{M}$, some lifts of those generators in $M$ will
generate~$M$. If we asume this, it follows that the number of
generators of $M$ is $\Dim(\overline{M})/\vert \overline{H} \vert =
4\Dim(\overline{M})/\vert G \vert$.
If $r = 4\Dim(\overline{M})/\vert G \vert$, then the projective cover of
$M$ is the free module $(kG)^r$. Using Lemma~\ref{struc-critical} we obtain
\begin{eqnarray*}
\Dim(\Omega(M))&  =& \Dim((kG)^r) - \Dim(M) \\
&=& 4\Dim(\overline{M}) -
2\Dim(\overline{M}) - 1  = \Dim(M) -2
\end{eqnarray*}
 as desired. Finally, since the dual module $M^*$ also satisfies the
assumptions of the lemma, we have that
 \begin{eqnarray*}
\Dim(\Omega^{-1}(M)) &= &\Dim(\Omega^{-1}(M)^*) \\
&=& \Dim(\Omega(M^*))
=  \Dim(M^*) - 2 = \Dim(M) - 2 
\end{eqnarray*}
and this completes the proof. 

We are left with the proof of the second claim. Let $L$ be the submodule
of~$M$ generated by some lifts in $M$ of the generators of~$\overline{M}$.
Assume by contradiction that $L\neq M$. Since $M{\downarrow}_H^G = k
\oplus F$ for some free $kH$-module~$F$, we have
$\overline{M}{\downarrow}_{\overline{H}}^{\overline{G}} = F/(z-1)F$ and so
we can choose the lifts of the generators of~$\overline{M}$ so that
$L{\downarrow}_H^G = F$. Now for any other maximal subgroup $H'$ of~$G$,
we have $M{\downarrow}_{H'}^G = k \oplus F'$ for some free
$kH'$-module~$F'$. The subgroup $H\cap H'$ is nontrivial because it
contains~$Z$ and there are two decompositions
 $$M{\downarrow}_{H\cap H'}^G
= T{\downarrow}_{H\cap H'}^H\oplus F{\downarrow}_{H\cap H'}^H
= T'{\downarrow}_{H\cap H'}^{H'}\oplus F'{\downarrow}_{H\cap H'}^{H'}$$
where $T$, respectively~$T'$, denotes a trivial one-dimensional module for
$kH$, respectively~$kH'$. By comparing the fixed points
$M^{H\cap H'}$ and the relative traces $t_1^{H\cap H'}\cdot M$ in both
decompositions, we see that
$T'{\downarrow}_{H\cap H'}^{H'}$ cannot be contained in
$F{\downarrow}_{H\cap H'}^H$ and therefore
$$M{\downarrow}_{H\cap H'}^G
= T'{\downarrow}_{H\cap H'}^{H'}\oplus F{\downarrow}_{H\cap H'}^H$$
(see Lemma~8.2 in~\cite{CT1} for details). Since $F$ is the restriction of
a $kG$-submodule, this is a decomposition of~$M$ as a $kH'$-module, namely
$$M{\downarrow}_{H'}^G = T' \oplus L{\downarrow}_{H'}^G\,.$$
By the Krull-Schmidt theorem, we deduce that $L{\downarrow}_{H'}^G$ is
free. Since this holds for any maximal subgroup~$H'$ and since $G$ is not
elementary abelian, Chouinard's theorem (see \cite{Bbook} or \cite{Ev})
implies that $L$ is free as a $kG$-module and so $M \cong k\break \oplus L$. But
$M$ is indecomposable and nontrivial by assumption. This contradiction
completes the proof of the claim.
\Endproof\vskip4pt  

For our next theorem, we first need a technical lemma.

\begin{lemma} \label{split}
Let $W$ be a $kG$-module satisfying the following two conditions\/{\rm :}\/
\begin{itemize}
\ritem{(a)} $W/(z-1)W = U_1 \oplus U_2$ where $U_1$ and $U_2$ are
$k\overline{G}$-submodules such that the varieties satisfy
$V_{\overline{G}}(U_1) \cap
V_{\overline{G}}(U_2) = \{0\}$.

\ritem{(b)} For some $r\leq p$\/{\rm ,}\/ there is $(z-1)^r W = 0$ and $W$ is free as a module
over the ring $kZ/(z-1)^r$.
\end{itemize}
Then $W = W_1 \oplus W_2$ where $W_1$ and $W_2$ are $kG$-submodules of $W$
such that ${W_i/(z-1)W_i} \cong U_i$ for $i=1,2$.
\end{lemma}

\Proof 
We use induction on~$r$. There is nothing to prove if $r=1$
so we assume $r\geq 2$. By induction, $W/(z-1)^{r-1}W = V_1 \oplus V_2$
where $V_1$ and $V_2$ are $kG$-submodules of $W/(z-1)^{r-1}W$ such that
${V_i/(z-1)V_i} \cong U_i$ for $i=1,2$. Now, since $W$ is free as a module
over $kZ/(z-1)^r$, multiplication by $(z-1)$ induces an
isomorphism $W/(z-1)^{r-1}W \cong (z-1)W$ and we write $L_i$ for the image
of~$V_i$. So $(z-1)W = L_1 \oplus L_2$.

Let $\pi: W \to W/(z-1)W = U_1\oplus U_2$ be the canonical surjection.
Passing to the quotient by~$L_1$, we obtain a short exact sequence
$$ 0 \Rarr{} L_2 \Rarr{} W/L_1 \Rarr{\widetilde\pi} U_1\oplus U_2
\Rarr{} 0 $$
where $\widetilde\pi$ is induced by~$\pi$.
Let $K = \{ x\in W/L_1 \,\mid\, (z-1)x = 0 \}$. We claim that
$\widetilde\pi(K)=U_1$. Let $x\in K$ and let $w\in W$ be a lift of~$x$.
Then $(z-1)w\in L_1$. Since multiplication by $(z-1)$ induces an
isomorphism $W/(z-1)^{r-1}W \cong (z-1)W$, the class of~$w$ in
$W/(z-1)^{r-1}W$ is in~$V_1$. It follows that $\pi(w)\in U_1$, hence
$\widetilde\pi(x)\in U_1$, proving the claim.

Therefore we obtain a short exact sequence
$$ 0 \Rarr{} (z-1)^{r-2}L_2 \Rarr{} K \Rarr{\widetilde\pi} U_1 \Rarr{} 0 $$
because $L_2 \cap \Ker(z-1) = (z-1)^{r-2}L_2$. This is a sequence of
$k\overline G$-modules since $(z-1)K=0$ by construction. Now
multiplication by $(z-1)^{r-1}$ induces an isomorphism
$W/(z-1)W \cong (z-1)^{r-1}W$ mapping $U_2$ onto $(z-1)^{r-2}L_2$.
By applying our assumption on the varieties of $U_1$ and~$U_2$
we deduce that the sequence splits (see Theorem~\ref{propvar}).
Let $\sigma$ be a section of $\widetilde\pi:K \to U_1$ and let $W_1$ be the
inverse image of $\sigma(U_1)$ in~$W$, so that $W_1/L_1 = \sigma(U_1)$.
We have obtained a short exact sequence
$$ 0 \Rarr{} L_1 \Rarr{} W_1 \Rarr{\pi} U_1 \Rarr{} 0 \,.$$
We can construct similarly a submodule $W_2$ and a short exact sequence
$$ 0 \Rarr{} L_2 \Rarr{} W_2 \Rarr{\pi} U_2 \Rarr{} 0 \,.$$
Then $\pi(W_1\cap W_2)=0$, so that  $W_1\cap W_2 \subseteq \Ker(\pi)=L_1 \oplus
L_2$. But since $W_i\cap \Ker(\pi) = L_i$, we obtain $W_1\cap W_2=0$.
For reasons of dimensions (or by a direct argument), the direct sum
$W_1\oplus W_2$ must be the whole of~$W$.
\hfill\qed

\begin{thm} \label{separate-var}
Let $M$ be a critical $kG$-module and suppose that $\overline{M} =
\overline{M}_1 \oplus \overline{M}_2$ where $\overline{M}_1$
and $\overline{M}_2$ are $k\overline{G}$-submodules. Suppose that the
varieties satisfy
$$
V_{\overline{G}}(\overline{M}_1) \cap
V_{\overline{G}}(\overline{M}_2) = \{0\}.
$$
Then there exist critical $kG$-modules $N_1$ and $N_2$ such that
$\overline{N}_i \cong \overline{M}_i$ for $1 \leq i \leq 2$.
\end{thm}

\Proof  As before, let $M' = \{m \in M
\;|\;(z - 1)^{p-1}m = 0\}$.
Let $M_1 \subseteq M$ be the inverse
image of $\overline{M}_1$ under the quotient map $M
\longrightarrow M/M' = \overline{M}$.  Let $M_2$ be the
inverse image of $\overline{M}_2$.  Then $M' = M_1 \cap M_2$
and $M_1/M' \cong \overline{M}_1$, $M_2/M' = \overline{M}_2$.

By Lemma~\ref{struc-critical}, $(z-1)M$ is free over $kZ/(z-1)^{p-1}$ and
$$ (z-1)M/(z-1)^2M \cong M/M' = \overline{M} = \overline{M}_1 \oplus
\overline{M}_2 \,.$$
Therefore Lemma~\ref{split} applies and we have $(z-1)M = W_1 \oplus W_2$
such that ${W_i/(z-1)W_i \cong \overline{M}_i}$ for $i=1,2$. Now define
$N_1 = M_1/W_2$ and $N_2 = M_2/W_1$. If $r_i = \Dim(\overline{M}_i)$, then
$\Dim(\overline{M})= r_1 + r_2$ and by
Lemma~\ref{struc-critical} we obtain $\Dim(M) = pr_1 + pr_2 +1$ and
$\Dim((z-1)M) = (p-1)r_1 + (p-1)r_2$. Therefore we have $\Dim(M_1) = pr_1
+ (p-1)r_2 +1$ and $\Dim(M_2) = (p-1)r_1 + pr_2 +1$. Also
$\Dim(W_i) = (p-1)r_i\,$; hence $\Dim(N_i) = pr_i +1$ for $i=1,2$.

We claim that $N_1{\downarrow}_H^G \cong k \oplus ({\rm free})$ for every
maximal subgroup $H$ of $G$ (and similarly for $N_2$).
Let $H = \langle z, y_1, \dots, y_{m-1}\rangle$ where
$\overline{y}_1, \dots , \overline{y}_{m-1}$ are generators of
$\overline{H}=H/Z$. The assumption on $M{\downarrow}_H^G$ implies that
$\overline{M}$ is free as a $k\overline{H}$-module. Therefore
$\overline{M}_1$ and $\overline{M}_2$ must be free as
$k\overline{H}$-modules.
Let $Y = \displaystyle\prod_{i=1}^{m-1}(y_i-1)^{p-1}$
so that $\overline{Y} = t_1^{\overline{H}}$ and $Y(z-1)^{p-1} = t_1^H$.
Then we get
$$\Dim(\overline{M}_1) = \vert \overline{H} \vert \cdot \Dim(\overline{Y}
\cdot \overline{M}_1) \,.$$
Now $(z-1)^{p-1}N_1 \cong (z-1)^{p-1}M_1$ because $N_1=M_1/W_2$ and
$(z-1)^{p-1}W_2=0$. Therefore
 $$t_1^H\cdot N_1 = Y(z-1)^{p-1}N_1 \cong Y(z-1)^{p-1}M_1
\cong Y\cdot\overline{M}_1 = \overline{Y}\cdot\overline{M}_1\,.$$
It follows that
 $$\vert H \vert\Dim(t_1^H\cdot N_1)
= p\cdot \vert \overline{H} \vert \cdot
\Dim(\overline{Y}\cdot\overline{M}_1)
= p\cdot \Dim(\overline{M}_1) = pr_1 = \Dim(N_1) -1 \,.$$
 Therefore $N_1{\downarrow}_H^G$ has a free submodule of dimension
$\Dim(N_1) - 1$. The only way this can happen is if $N_1{\downarrow}_H^G
\cong k \oplus ({\rm free})$.

Now we prove that $\overline{N}_1 \cong \overline{M}_1$ (and
similarly for $N_2$). We have to compute the submodule
$N'_1 = \{ x\in N_1 \;\mid\, (z-1)^{p-1}x =0 \}$. But $N_1 = M_1/W_2$ and
we have $W_2\subseteq M' \subseteq M_1$ and $(z-1)^{p-1}M'=0$. Therefore
$M'/W_2 \subseteq N'_1$ and $\overline{N}_1 = N_1/N'_1$ is a quotient of
$N_1/(M'/W_2) \cong M_1/M' = \overline{M}_1$. In order to prove that this
is not a proper quotient, it suffices to prove that $\overline{N}_1$ and
$\overline{M}_1$ have the same dimension. But by the previous part of
the proof, we know that $N_1{\downarrow}_H^G
\cong k \oplus ({\rm free})$ for every maximal subgroup~$H$. By
Lemma~\ref{struc-critical} this implies
 $$\Dim(\overline{N}_1) = {\Dim(N_1)-1 \over p} = r_1
= \Dim(\overline{M}_1)\,,$$
as was to be shown.

Finally we conclude that $N_1$ is critical. Indeed, since $\overline{M}$
has no free summand as a $k\overline{G}$-module, $\overline{N}_1$ cannot
have a free summand and therefore $N_1$ has no free summand as a
$kG$-module by Lemma~\ref{barperiod}. This implies that $N_1$ is critical
since we know that $N_1{\downarrow}_H^G
\cong k \oplus ({\rm free})$ for every maximal subgroup~$H$.
\phantom{endofpage}\hfill\qed

\begin{thm} \label{module-products}
Let $M_1$ and $M_2$ be critical $kG$-modules and suppose that the
varieties satisfy
$$
V_{\overline{G}}(\overline{M}_1) \cap
V_{\overline{G}}(\overline{M}_2) = \{0\}.
$$
Then $M_1 \otimes M_2 \cong M \oplus ({\rm free})$ where $M$ is a critical
$kG$-module such that $\overline{M} \cong
\overline{M}_1 \oplus \overline{M}_2$.
\end{thm}

\Proof 
Let $r_j = \Dim(\overline{M}_j)$ for $j=1,2$. Thus $\Dim(M_j)=pr_j+1$.
Consider the filtration of $M_1$ as in Lemma~\ref{struc-critical}
$$ \{0\} \subset (z-1)^{p-1} M_1 \subset K_1 \subset \dots \subset
K_{p-1} \subset K_p = M_1\,,$$
where $K_i = \{ m \in M_1 \;|\;(z - 1)^i m = 0\}$.
This induces a filtration on $M_1 \otimes M_2$
$$
\{0\} \subset (z-1)^{p-1} M_1 \otimes M_2 \subset K_1 \otimes M_2 \subset
\dots \subset K_{p-1} \otimes M_2 \subset M_1 \otimes M_2\,,
$$
with all quotients but one isomorphic to $\overline{M}_1 \otimes M_2$.
We need to prove the following.

\begin{lemma} \label{tenprod1}  $\overline{M}_1 \otimes M_2 =
F \oplus L$ where $L\cong \overline{M}_1$ and $F$ is a free $kG$-module
of dimension $pr_1r_2$ such that $(z-1)^{p-1}F = \overline{M}_1 \otimes
(z-1)^{p-1} M_2$.
\end{lemma}

\Proof  By hypothesis
$V_{\overline{G}}(\overline{M}_1)
\cap V_{\overline{G}}(\overline{M}_2) = \{0\}$ and hence
$\overline{M}_1 \otimes \overline{M}_2$ is projective as a
$k\overline{G}$-module. Choose elements $m_1, \dots , m_r \in
\overline{M}_1 \otimes M_2$ such that $\overline{m}_1, \dots ,
\overline{m}_r$ is a free
$k\overline{G}$-basis for $\overline{M}_1 \otimes \overline{M_2}$.
Here $\overline{m}_i = m_i + (\overline{M}_1 \otimes M'_2)$ denotes the
class of $m_i$ in $\overline{M}_1 \otimes \overline{M_2} =
(\overline{M}_1 \otimes M_2)/(\overline{M}_1 \otimes M'_2)$.

As before, let
$X = \displaystyle\prod_{i=1}^{m}(x_i-1)^{p-1}$ so that $\overline{X} =
t_1^{\overline{G}}$ and $X(z-1)^{p-1}\break = t_1^G$.
Then $\overline{X}
\overline{m}_1, \dots, \overline{X} \overline{m}_r$ are linearly
independent in $\overline{M}_1 \otimes \overline{M}_2$. Since $z$ acts
trivially on $\overline{M}_1$, multiplication by $(z-1)^{p-1}$ induces an
isomorphism
$\overline{M}_1 \otimes \overline{M}_2
\cong \overline{M}_1 \otimes (z-1)^{p-1} M_2$
and it follows that $X(z-1)^{p-1} m_1,
\dots, X(z-1)^{p-1} m_r$ are linearly independent in $\overline{M}_1
\otimes (z-1)^{p-1}M_2$. Therefore $t_1^G m_1,
\dots, t_1^G m_r$ are linearly independent in $\overline{M}_1
\otimes M_2$.  So the $kG$-submodule $F$ of $\overline{M}_1
\otimes M_2$ generated by $m_1, \dots, m_r$ is a free $kG$-module.
Moreover we have $(z-1)^{p-1}F = \overline{M}_1 \otimes (z-1)^{p-1} M_2$.

Consider now the exact sequence of $k\overline{G}$-modules
$$
0  \Rarr{} \overline{M}_1 \otimes (z-1)^{p-1} M_2 \Rarr{} \overline{M}_1
\otimes (M_2)^Z \Rarr{} \overline{M}_1 \otimes k \Rarr{} 0 \,,
$$
where $(M_2)^Z = \{ x\in M_2 \;|\; (z - 1) x = 0\}$.
Since the kernel is free over $k\overline{G}$, the sequence splits and we
have $\overline{M}_1 \otimes (M_2)^Z \cong
(z-1)^{p-1} F \oplus L$ where $L$ is a submodule
isomorphic to $\overline{M}_1$. If we had $F\cap L \neq 0$, then we would
have $\Soc(F)\cap L \neq 0$; hence $(z-1)^{p-1} F\cap L \neq 0$, a
contradiction. Therefore $F\cap L = 0$ and $\overline{M}_1 \otimes M_2$
contains a submodule $F \oplus L$.

We now show that $F \oplus L = \overline{M}_1 \otimes M_2$ by proving that
both modules have the same dimension. We have
$\Dim(F) = p \cdot \Dim(\overline{M}_1 \otimes \overline{M}_2)
= p r_1r_2$ and therefore
$$
\Dim(F \oplus L) = pr_1r_2 + r_1 = r_1(pr_2 +1)= \Dim(\overline{M}_1)
\Dim(M_2)\,,
$$
as was to be shown.
\Endproof\vskip4pt  

Now, continuing with the proof of the theorem, we note that each quotient
$(K_{i+1} \otimes M_2) / (K_{i} \otimes M_2)$ is isomorphic to
$\overline{M}_1 \otimes M_2$, hence contains a free submodule
$\overline{F_i}$ of dimension $pr_1r_2$ by the lemma.
Now remember that projective modules are also injective and, as a result,
if a projective module is a direct summand of a section of a module~$V$,
then it is a direct summand of~$V$.\break Thus we can lift the
free module $\overline{F_i}$ and obtain a free submodule $F_i$ of\break
$M_1\otimes M_2$ mapping isomorphically onto $\overline{F_i}$ under the
quotient map $M_1\otimes M_2
\to (M_1\otimes M_2) / (K_{i} \otimes M_2)$. Similarly, $(z-1)^{p-1} M_1
\otimes M_2$ is isomorphic to $\overline{M}_1 \otimes M_2$, hence contains
a free submodule $F_0$ of dimension $pr_1r_2$ by the lemma. Therefore we
have
$$M_1\otimes M_2 = M \oplus F$$
where $F = F_0 \oplus \dots \oplus F_{p-1}$ is free of dimension
$p^2r_1r_2$ and $M$ is a submodule of dimension
$(pr_1+1)(pr_2+1)-p^2r_1r_2 = p(r_1+r_2) +1$.

Since, for any maximal subgroup $H$ of~$G$, we have $M_j{\downarrow}_H^G
\cong k \oplus ({\rm free})$ for $j=1,2$, the same holds for $M_1\otimes
M_2$ and hence $M{\downarrow}_H^G \cong k \oplus ({\rm free})$.
We are going to prove that $\overline{M} \cong \overline{M}_1
\oplus \overline{M}_2$. This will imply that $M$ is critical. Indeed
$\overline{M}_j$ has no $k\overline{G}$-free summand, because $M_j$ is
critical ($j=1,2$), so $\overline{M}_1 \oplus \overline{M}_2$ has no
$k\overline{G}$-free summand and therefore $M$ has no $kG$-free
summand by Lemma~\ref{struc-critical}. This forces the endo-trivial
module~$M$ to be indecomposable.

Instead of working with $\overline{M}$, we consider the isomorphic module\break
$(z-1)^{p-1}M$ and our goal now is to prove that
$(z-1)^{p-1}M \cong \overline{M}_1 \oplus \overline{M}_2$.
We work with the submodule $K_1 \otimes M_2$ of our filtration and
we first analyze its submodule $K_1 \otimes K_1'$,
where $K_1' = \{m \in M_2\,|\,(z-1)m =0\}$ is the analog of $K_1$ for~$M_2$.  Notice that
$K_1 \otimes K_1'$ is a $k\overline{G}$-module with a filtration
\begin{eqnarray*}
(z-1)^{p-1} M_1 \otimes (z-1)^{p-1} M_2& \subset&
\big((z-1)^{p-1} M_1 \otimes K_1'\big)
\\
&&+ \big(K_1 \otimes (z-1)^{p-1} M_2\big)
\subset K_1 \otimes K_1'\,.
\end{eqnarray*}
In the filtration, the bottom submodule is 
free over $k\overline{G}$ and is equal to\break
$(z-1)^{p-1}F_0$ by the lemma. The middle quotient of this filtration is
the direct sum of
\begin{eqnarray*}
&&((z{-}1)^{p-1} M_1 \otimes K_1')/((z{-}1)^{p-1} M_1 \otimes (z{-}1)^{p-1}
M_2)\\
&&\hskip.5in\cong (z{-}1)^{p-1} M_1 \otimes k \cong \overline{M}_1\end{eqnarray*}
and
\begin{eqnarray*}
&&(K_1 \otimes (z{-}1)^{p-1} M_2)/((z{-}1)^{p-1} M_1 \otimes (z{-}1)^{p-1}
M_2)\\
&&\hskip.5in\cong k \otimes (z{-}1)^{p-1} M_2 \cong \overline{M}_2\,.
\end{eqnarray*}
This direct sum can be lifted in $K_1 \otimes K_1'$, because
the submodule $(z-1)^{p-1}F_0 = (z{-}1)^{p-1} M_1 \otimes
(z{-}1)^{p-1} M_2$ is $k\overline{G}$-free and so the sequence
$$ 0  \Rarr{} (z-1)^{p-1}F_0 \Rarr{} K_1 \otimes K_1' \Rarr{}
(K_1 \otimes K_1')/(z-1)^{p-1}F_0 \Rarr{} 0$$
splits. Therefore $K_1 \otimes K_1'$ contains a submodule $V_1 \oplus
V_2$ with $V_j \cong \overline{M}_j$ and $(V_1 \oplus V_2)\cap
(z-1)^{p-1}F_0 = 0$. It follows that $(V_1 \oplus V_2)\cap \Soc(F_0) = 0$
and so $(V_1 \oplus V_2)\cap F_0 = 0$.

We now have $V_1 \oplus V_2 \oplus F_0 \subset K_1 \otimes M_2$ and
therefore
$V_1 \oplus V_2 \oplus F_0$ intersects trivially
$F_1 \oplus \dots \oplus F_{p-1}$ because this free module has been
lifted from quotients of
${(M_1 \otimes M_2) / (K_1 \otimes M_2)}$. This shows that $M_1 \otimes
M_2$ contains the submodule
$V_1 \oplus V_2 \oplus F_0 \oplus F_1 \oplus \dots \oplus F_{p-1}
= V_1 \oplus V_2 \oplus F$. Therefore $V_1 \oplus V_2$ is isomorphic to a
submodule of~$M$.

Now we show that $V_1 \oplus V_2 \subset (z-1)^{p-1} (M_1 \otimes M_2)$.
Since $z$ acts trivially on~$K_1'$, we have $(z-1)^{p-1} (M_1 \otimes K_1')
= (z-1)^{p-1} M_1 \otimes K_1'$ and this contains $V_1$ by construction
of~$V_1$. Similarly $V_2 \subset K_1 \otimes (z-1)^{p-1} M_2 = (z-1)^{p-1}
(K_1 \otimes M_2)$. Passing to the quotient by~$F$, we deduce that
$V_1 \oplus V_2$ is isomorphic to a submodule of
$$(z-1)^{p-1} \big((M_1 \otimes M_2)/F\big)
= (z-1)^{p-1} \big((M \oplus F)/F\big)
\cong (z-1)^{p-1} M\,.$$
In order to prove that this submodule is the whole of $(z-1)^{p-1} M$, it
suffices to prove that they have the same dimension. But $V_j \cong
\overline{M}_j$ has dimension~$r_j$ (for $j=1,2$) and we know that $\Dim(M)
= p(r_1+r_2)+1$. Therefore $\Dim(\overline{M}) = r_1+r_2$ and we are done.
This shows that $(z-1)^{p-1} M \cong V_1 \oplus V_2$ and completes the
proof of the theorem.
\Endproof\vskip4pt  

Theorem \ref{separate-var} and Theorem \ref{module-products} provide
the basic tools for constructing a large critical module from any
given finite set of such modules, as follows.

\begin{thm}  \label{dim-endo}
For every $i=1, \dots, t$, let $M_i$ be a nontrivial critical
$kG$-module. Let $\ell_i$ be a line in the
variety of the periodic $k\overline{G}$-module~$\overline{M}_i$ and
assume that $\ell_i \neq \ell_j$ for $i \neq j$. Then there exists a
nontrivial critical $kG$-module $M$ such that
$V_{\overline{G}}(\overline{M}) = \bigcup^t_{i=1}\ell_i$. Moreover{\rm ,}
$\Dim(M) \ge t\vert G \vert/2 + 1$ if $p=2$ and $\Dim(M) \ge t\vert G
\vert + 1$ if $p$ is odd.
\end{thm}

\Proof   Recall that $\overline{M}_i$ is periodic by Lemma
\ref{barperiod}, and hence $V_{\overline{G}}(\overline{M}_i)$ is a union
of lines. By Theorem~\ref{propvar}, $\overline{M}_i = L_i \oplus N_i$ such
that $V_{\overline{G}}(L_i) = \ell_i$ and the variety of $N_i$ is the
union of the other lines (if any; otherwise  simply set  $L_i =
\overline{M}_i$). By Theorem~\ref{separate-var}, there exists a
critical $kG$-module~$U_i$ such that $\overline{U}_i = L_i$.

Now by Theorem \ref{module-products} and the assumption that the lines
$\ell_i$ are distinct, we obtain a critical $kG$-module $M$ such
that
 $$U_1 \otimes U_2 \otimes \dots \otimes U_t = M \oplus ({\rm free})$$
and $\overline{M} = \overline{U}_1 \oplus \overline{U}_2 \oplus \dots
\oplus \overline{U}_t$ so that $V_{\overline{G}}(\overline{M}) =
\bigcup^t_{i=1}\ell_i$.

Since $U_i{\downarrow}_H^G \cong k \oplus ({\rm free})$ where $H$ is a
maximal subgroup of~$G$, $\Dim(U_i)-1$ is a multiple of~$\vert H \vert$
and therefore $\Dim(\overline{U}_i)$ is a multiple of
$\vert H \vert/p = \vert G \vert/p^2$. It follows that $\Dim(\overline{M})
\geq t \vert G \vert/p^2$ and $\Dim(M) \geq t \vert G \vert/p +1$. We can
do better if $p$ is odd because $U_i$ is an endo-trivial module and so
$\Dim(U_i) \equiv \pm1\pmod{\vert G \vert}$ by Lemma~2.10 in~\cite{CT1}.
A plus sign is forced here and therefore $\Dim(U_i)-1$ is a multiple
of~$\vert G \vert$. The same argument then yields $\Dim(M) \geq t \vert G
\vert +1$.
\Endproof\vskip4pt  

{\it Remark.}
By a theorem of Puig~\cite{Pu1}, the torsion subgroup
$T_t(G)$ is finite. Therefore, there are actually finitely many possible
choices for the modules $M_i$ in the last theorem. It then follows from
the theorem that one can construct an indecomposable torsion endo-trivial
module $M$ such that $V_{\overline{G}}(\overline{M})$ contains
$V_{\overline{G}}(\overline{N})$ for any torsion endo-trivial module $N$.
Moreover, $\Dim(M) \ge t\vert G \vert/2 + 1$, respectively $t\vert G
\vert + 1$, where $t$ is the number of components of
$V_{\overline{G}}(\overline{M})$. However, in view of the main theorem of
this paper, it will turn out that $T_t(G)=0$ and so
$M \cong k$.

\section{Lower bounds on dimensions of critical modules}

In this section we prove a theorem that is
essential to the general cases of our main result. Basically it says
that, if an extraspecial group or an almost extraspecial group
has a nontrivial critical module, then it has one of large
dimension. For the proof, we need a few lemmas.

\begin{lemma} \label{fprationality}
Suppose that $M$ is a nontrivial critical $kG$-module and let $\ell$ be a
line in $V_{\overline{G}}(\overline{M})$. Then $\ell$ is not contained in
any ${\Bbb F}_p$-rational subspace of $V_{\overline{G}}(k)$.
\end{lemma}

\Proof 
Note that $V_{\overline{G}}(k) = k^m$ where $\vert G \vert = p^{m+1}$.
An ${\Bbb F}_p$-rational subspace (i.e.\ a subspace defined
by a linear equation with ${\Bbb F}_p$-coefficients)  corresponds to
a maximal subgroup $\overline{H} \subseteq \overline{G}$. That is, the
${\Bbb F}_p$-rational subspaces of $V_{\overline{G}}(k)$ are precisely the
subspaces of the form
$\res{\overline{G}}{\overline{H}}^*(V_{\overline{H}}(k))$.  If $\ell$
were in $\res{\overline{G}}{\overline{H}}^*(V_{\overline{H}}(k))$
then it would have to be the case that
$V_{\overline{H}}(\overline{M}{\downarrow}_{\overline{H}}^{\overline{G}})
\neq \{0\}$ and hence
$\overline{M}{\downarrow}_{\overline{H}}^{\overline{G}}$ would not be
free as  a $k\overline{H}$-module. By Lemma~\ref{barperiod}, this would
contradict the hypothesis that $M{\downarrow}_H^G \cong k \oplus ({\rm
free})$ where $H$ is the inverse image of $\overline{H}$ in $G$.
\Endproof\vskip4pt  

Recall that a $p'$-group is a group of order prime to~$p$.

\begin{lemma} \label{stabilizer}
Suppose that $\ell$ is a line through the origin in
$V_{\overline{G}}(k) = k^m$ and suppose that $\ell$ is not contained
in any ${\Bbb F}_p$-rational subspace of~$k^m$. Then the stabilizer $S$ of
$\ell$ for the action of ${\rm GL}_m({\Bbb F}_p)$ on $k^m$ is a cyclic
$p'$-subgroup.
\end{lemma}

\Proof 
Suppose that $y \in {\rm GL}_m({\Bbb F}_p)$ stabilizes $\ell$ and that $v$
is a point on $\ell$. Then $v$ is an eigenvector of $y$ with eigenvalue
$\lambda$. That is, simply, $y \cdot v = \lambda v$. So the line $\ell$
is a $kS$-submodule for the  action of $S$ on $k^m$, corresponding to a
homomorphism $\rho:{S \Rarr{} {\rm GL}(\ell) \cong k^*}$ mapping
$y\in S$ to the eigenvalue~$\lambda$.

We claim that $\rho$ is injective on the stabilizer $S$.
For suppose that $\rho(y) = \lambda = 1$. Then, viewing $y$ as a matrix,
we have that $(y - I_m)v = 0$. If $y$ is not the 
identity then some row $(a_1, a_2,
\dots, a_m)$ of $y - I_m$ is not zero. But then $v$ is in the subspace
defined by the equation $a_1x_1 + a_2x_2 + \dots + a_mx_m = 0$. Because
the coefficients of $y$ are in ${\Bbb F}_p$ we have a contradiction.

Now $S$ is isomorphic to a finite subgroup of $k^*$ and therefore it must
consist of roots of unity. Thus it is a cyclic $p'$-group and we are done.
\Endproof\vskip4pt  

Recall that for any automorphism $\alpha$ of $G$, the conjugate
module $N^\alpha$ is defined to be the $k$-vector space $N$ with the
action of $G$ given by $g\cdot n = \alpha(g) n$  for $g \in G$ and $n \in
N$. If $\alpha$ is an inner automorphism of $G$, then $N^\alpha \cong N$
and it follows that the group $\Out(G)$ of outer automorphisms of
$G$ acts on the set of isomorphism classes of $kG$-modules. We shall
also write $N^y$ for a conjugate module defined by an outer automorphism
$y \in \Out(G)$.

Since $G$ is extraspecial or almost extraspecial, we control $\Out(G)$
in the following sense. Recall that if $p=2$, there is an associated
quadratic form on the ${\Bbb F}_2$-vector space $G/Z(G)$ (see
Lemma~\ref{quadforms}). If $p$ is odd, there is a symplectic form~$b$ on the
${\Bbb F}_p$-vector space $G/Z(G)$ defined by $[\tilde{x}, \tilde{y}] =
z^{b(x,y)}$, where $x,y \in G/Z(G)$, $\tilde{x},\tilde{y} \in G$ are
elements of $G$ that lift $x$ and $y$, and $z$ is a generator of~$Z(G)$.

\begin{lemma} \label{out}
Let $G$ be an extraspecial or almost extraspecial $p$-group.
Let $\Out_0(G)$ be the subgroup of $\Out(G)$ consisting of outer
automorphisms fixing the center $Z(G)$ pointwise.
\begin{itemize}
\ritem{(a)} If $p$ is odd and $G$ is extraspecial of exponent $p$\/{\rm ,}\/ then $\Out_0(G)$
is isomorphic to the symplectic group $O_G$ associated to the symplectic
form corresponding to~$G$.

\ritem{(b)} If $p=2${\rm ,} $\Out_0(G)$ is isomorphic to the orthogonal group $O_G$
associated to the quadratic form corresponding to~$G$. \end{itemize}
\end{lemma}

\Proof 
When $G$ is extraspecial, this is one of the main results in Winter's
paper~\cite{Win}. If $G$ is almost extraspecial, the arguments given in
Sections 3F and~4 of~\cite{Win} extend and yield the same result.
Alternatively, this appears explicitly in Exercise~5 of Chapter~8
of~\cite{Asch}.
\hfill\qed

\begin{thm}\label{lowercriterion} \hskip-4pt
Suppose that there exists a nontrivial critical $kG$-module~$N$.
\begin{itemize}
\ritem{(a)} If $p$ is odd{\rm ,} there exists a critical $kG$-module $M$ such that
$$
\Dim(M) > \vert G\vert \cdot
\dfrac{\vert O_G\vert}{\vert C\vert}
$$
where $O_G$ is the symplectic group associated to $G$ and $C$ is a cyclic\break
$p'$-subgroup of $O_G$ of maximal order.

\ritem{(b)} If $p=2$\/{\rm ,}\/ there exists a critical $kG$-module $M$ such that
$$
\Dim(M) > \dfrac{\vert G\vert}{2} \cdot
\dfrac{\vert O_G\vert}{\vert C\vert}
$$
where $O_G$ is the orthogonal group associated to $G$ and $C$ is an odd
order cyclic subgroup of $O_G$ of maximal order. \end{itemize}
\end{thm}

\Proof 
Let $\ell$ be a line in $V_{\overline{G}}(\overline{N})$. Notice that if
$y \in O_G$ then $N^y$ is also a nontrivial $kG$-module such that
$N^y{\downarrow}_H^G \cong k \oplus ({\rm free})$. But then
$y(\ell)$ is in the variety $V_{\overline{G}}(\overline{N}^y)$. If $B$
denotes the stabilizer of $\ell$ in $O_G$, we obtain a family of
modules $N^y$ indexed by the set of cosets $O_G/B$. So
by Theorem \ref{dim-endo}, there exists a critical $kG$-module
$M$ such that 
$V_{\overline{G}}(\overline{M}) = \bigcup_{y\in O_G/B} y(\ell)$. Moreover,
$\Dim M > \dfrac{\vert O_G\vert}{\vert B\vert} \cdot
\dfrac{\vert G\vert}{2}$ if $p=2$ and $\Dim M > \dfrac{\vert O_G\vert}{\vert
B\vert} \cdot \vert G\vert$ if $p$ is odd.

By Lemma \ref{fprationality}, the line $\ell$ is not contained in any
${\Bbb F}_p$-rational subspace of $V_{\overline{G}}(k) = k^m$. Thus by
Lemma~\ref{stabilizer}, the group $B = S \cap O_G$ is cyclic of order prime
to~$p$. If $C$ is of maximal order among cyclic $p'$-subgroups of $O_G$, we
deduce the lower bound of the statement.
\Endproof\vskip4pt  

For use in the following sections, we need to have some estimates
of the orders of the orthogonal and symplectic groups and their cyclic
$p'$-subgroups.

\begin{prop} \label{ogrouporders}
Let $G$ be an extraspecial or almost extraspecial $p$-group. Let $O_G$
be the orthogonal or symplectic group associated to $G$.
\begin{itemize}
\item[{\rm 1.}] 
If $p$ is odd and $G$ is extraspecial of exponent $p$ and
order $p^{2n+1}${\rm ,} then $O_G = {\rm Sp}(2n,{\Bbb F}_p)$ and
$$ 
\vert O_G \vert = p^{n^2}  \prod_{i=1}^{n} (p^{2i} -1)
\,.
$$
\item[{\rm 2.}] If $p=2$ and $G \cong D_8* \dots *D_8$ is extraspecial of
order $2^{2n+1}$ {\rm (}type~$1$\/{\rm ),}\/ then $O_G = O^+(2n,{\Bbb F}_2)$ and
$$ 
\vert O_G \vert = 2 \cdot 2^{n(n-1)} (2^n-1) \prod_{i=1}^{n-1} (2^{2i}-1)
\,.
$$
\item[{\rm 3.}] 
If $p=2$ and $G \cong D_8* \dots *D_8*Q_8$ is extraspecial of
order $2^{2n+1}$ \/{\rm (}\/type~{\rm 2),} then $O_G = O^-(2n,{\Bbb F}_2)$ and
$$ 
\vert O_G \vert = 2 \cdot 2^{n(n-1)}(2^n+1) \prod_{i=1}^{n-1} (2^{2i}-1)
\,.
$$
\item[{\rm 4.}] 
If $p=2$ and $G \cong D_8* \dots *D_8*C_4$ is almost extraspecial of
order $2^{2n+2}$ {\rm (}type~{\rm 3),} then $O_G = {\rm Sp}(2n,{\Bbb F}_2)$ and
$$ 
\vert O_G \vert = 2^{n^2}  \prod_{i=1}^{n} (2^{2i} -1)
\,.
$$
\end{itemize}
Moreover{\rm ,} if $C$ is any cyclic $p'$-subgroup of $O_G${\rm ,} then
$\vert C \vert \leq (p+1)^n$.
\end{prop}

\Proof 
In the first three cases, we have $O_G ={\rm Sp}(2n,{\Bbb F}_p)$, respectively
$O_G = O^{\pm}(2n,{\Bbb F}_2)$, essentially
by definition (see also~\cite{Win}). In the third case, we obtain
$O_G = O(2n+1,{\Bbb F}_2) \cong {\rm Sp}(2n,{\Bbb F}_2)$ (see Theorem~11.9 of
Taylor's book~\cite{Tay}) where  orders of the four groups appear on
pages~70 and~141. The list can also be found in
any of a number of text books on Chevalley groups or finite simple groups
(e.g.\ Gorenstein's book~\cite{Gor2}). The types of the groups of Lie type
in the four cases listed are $C_n(p)$, $D_n(2)$, $^2\!D_n(2)$ and $C_n(2)$,
respectively. In the first case the corresponding simple group is
$O_G/\{\pm 1\}$. In the next two cases the corresponding simple group has
index 2 in~$O_G$, while the group is simple in the fourth case.

For the statement about the cyclic $p'$-subgroups, note first that
elements of order prime to $p$ are semi-simple, hence contained in a maximal
torus. Now, for a Chevalley group of rank~$n$ over the field~${\Bbb F}_q$,
the order of a maximal torus is equal to
$$\vert \det(w^{-1}F-1)\vert \ = \ \vert \det(F-w)\vert$$
where $F(x)=x^q$ is the Frobenius morphism, $w$ is an element of the
corresponding Weyl group, and where $F$ and $w$ act on the cocharacter
group of a fixed maximal torus of the corresponding algebraic group (see
Proposition~3.3.5 in Carter's book~\cite{Cart}). Since
$w$ has finite order, we obtain a product
$\prod_{i=1}^n (q - \zeta_i)$ for suitable roots of unity~$\zeta_i$ (the
eigenvalues of~$w$). In our case, $q=p$ and the rank~$n$ is the same as
the integer~$n$ of the statement. Therefore if $C$ is any cyclic
$p'$-subgroup of $O_G$, we get
$$\vert C \vert \leq \vert \det(F-w)\vert =
 \big\vert \prod_{i=1}^n (p - \zeta_i) \big\vert
\leq \prod_{i=1}^n (p + 1) = (p+1)^n \,,$$
as was to be shown.
\hfill\qed

\section{Upper bounds on dimensions of critical modules}

Throughout the section we assume that $G$ is a $p$-group and that $k$ is
an algebraically closed field of characteristic $p$.

We will need the following results. Recall that a nonzero element $\zeta$
of $\HH{1}{G}{{\Bbb F}_p}$ corresponds to a maximal subgroup of $G$ in
the sense that there is a unique maximal subgroup $H$ of
$G$ such that $\resgh(\zeta) = 0$. When $p$ is odd, we also need the
Bockstein map $\beta:\HH{1}{G}{{\Bbb F}_p} \longrightarrow
\HH{2}{G}{{\Bbb F}_p}$ (see \cite{Bbook} or \cite{Ev}).

\begin{thm}\label{serrefilter}
Suppose that $G$ is a $p$-group which is
not elementary abelian.  Suppose that $\eta_1, \dots, \eta_t$
are nonzero elements in $\HH{1}{G}{{\Bbb F}_p}$ and have the property
that 
\begin{eqnarray*}
\eta_1  \dots \eta_t = 0 &&\quad \hbox{\it if $p=2$},
\\
\beta(\eta_1)  \dots \beta(\eta_t) = 0&& \quad \hbox{\it if $p$ is odd.}
\end{eqnarray*}
\begin{itemize}
\ritem{(a)} Assume that $p=2$. For each $i$\/{\rm ,}\/ let $H_i$ be the maximal
subgroup of $G$ corresponding to $\eta_i$.
Then there is a projective module $P$ such that $k \oplus
\Omega^{1-t}(k) \oplus P$ has a filtration
$$
\{0\} = L_0 \subseteq L_1 \subseteq \dots \subseteq L_t \cong k
\oplus \Omega^{1-t}(k) \oplus P
$$
where $L_i/L_{i-1} \cong (\Omega^{1-i} (k)){\uparrow}_{H_i}^G$ for
each $i=1,\dots, t$.

\ritem{(b)} Assume that $p$ is odd. For each $i${\rm ,} let $K_i$ be the maximal
subgroup of $G$ corresponding to $\eta_i$ and set $H_{2i}=H_{2i-1}=K_i$.
Then there is a projective module $P$ such that $k \oplus
\Omega^{1-2t}(k) \oplus P$ has a filtration
$$
\{0\} = L_0 \subseteq L_1 \subseteq \dots \subseteq L_{2t} \cong k
\oplus \Omega^{1-2t}(k) \oplus P
$$
where $L_i/L_{i-1} \cong (\Omega^{1-i} (k)){\uparrow}_{H_i}^G$ for
each $i=1,\dots, 2t$. \end{itemize}
\end{thm}

\Proof  This is the essence of Lemma ~3.10 of \cite{Celeab}.
That lemma is stated for $\bZ G$-modules but this does not really
matter since we can tensor the whole thing with~$k$.
Because the emphasis of our theorem is different from that of
the results of \cite{Celeab} we give a brief sketch of the proof
here. However, all of the ideas as well as the details are given
in the paper \cite{Celeab}.

(a) We first give the proof when $p=2$ and then indicate how to modify the
arguments for odd~$p$. Each of the cohomology elements $\eta_i$ corresponds
to an exact sequence
$$
0 \Rarr{} {\Bbb F}_2 \Rarr{} {\Bbb F}_2{\uparrow}_{H_i}^G
\Rarr{} {\Bbb F}_2 \Rarr{} 0 \,.
$$ 
Now we splice all of these together
and tensor with $k$ to get a sequence of the form
$$
0 \Rarr{} k \Rarr{} k{\uparrow}_{H_t}^G \Rarr{} \dots
\Rarr{} k{\uparrow}_{H_2}^G \Rarr{} k{\uparrow}_{H_1}^G \Rarr{}
k \Rarr{} 0 \,,
$$ 
which represents the element $\eta_1  \dots \eta_t = 0$ in
$\HH{t}{G}{k}$. Note that we are using the same notation $\eta_i$
for the element of $\HH{1}{G}{{\Bbb F}_2}$ and its image under the change
of rings in $\HH{1}{G}{k}$. Now we consider the complex $\CC$
obtained by truncating the ends off of the sequence. That is,
$\CC_i = k{\uparrow}_{H_{i+1}}^G$ for $i = 0, \dots, t-1$ and
$\CC_i = 0$ otherwise. We see that the homology of $\CC$ is a
result of the truncations. That is, $\HHH_i(\CC) = k$ if either
$i=0$ or $i=t-1$ and $\HHH_i(\CC) = 0$ otherwise.

The next step is to collapse the complex $\CC$ into a single
module. This is accomplished exactly as in the paragraphs preceding
Proposition ~3.7 of \cite{Celeab}. That is, we tensor, over $k$,
the complex $\CC$ with a
projective resolution of the trivial module $k$. This gives us a
projective resolution of the complex $\CC$ and it has the same
homology as $\CC$. Thus, in degrees above $t$, it is exact and is the
projective resolution of a module $U$, which we can take to be the image
of the $t^{th}$ boundary map of the total complex. The only problem
with $U$ is that it is in the wrong degree. So we take
$W = \Omega^{-t}(U)$. This is the module that we want.

There are now two things to note about $W$. First because the terms
of the complex $\CC$ are induced from the maximal subgroup $H_1,
\dots, H_t$, the module $W$ has a filtration by the modules
$k{\uparrow}_{H_i}^G$ suitably translated by~$\Omega$, exactly as
described in the statement of the theorem. That is, the projective
resolution of the complex $\CC$ as constructed above is filtered by the
projective resolutions of the terms of the complex, suitably translated.
See the proof of Proposition~3.8 of \cite{Celeab} for this part.

Next we note that the module $W$ is isomorphic to $k \oplus
\Omega^{1-t}(k) \oplus P$ for some projective module $P$. This is because
the original sequence that represented $\eta_1  \dots \eta_t$ splits and
hence the projective resolution of the complex is, in high degrees,
a projective resolution of the homology groups of the complex,
suitably translated. See Proposition ~3.8 of \cite{Celeab} for this
part. This proves the theorem if $p=2$.

(b) If $p$ is odd, the cohomology element $\eta_i$ has to be replaced by
its Bockstein $\beta(\eta_i)$ which corresponds to an exact sequence
$$
0 \Rarr{} {\Bbb F}_p \Rarr{} {\Bbb F}_p{\uparrow}_{K_i}^G \Rarr{}
{\Bbb F}_p{\uparrow}_{K_i}^G \Rarr{} {\Bbb F}_p \Rarr{} 0 \,.
$$ 
Again we splice all of these together and tensor with $k$. Using our
numbering of the subgroups $H_i$, we obtain a sequence of the form
$$
0 \Rarr{} k \Rarr{} k{\uparrow}_{H_{2t}}^G \Rarr{} \dots
\Rarr{} k{\uparrow}_{H_2}^G \Rarr{} k{\uparrow}_{H_1}^G \Rarr{}
k \Rarr{} 0 \,,
$$ 
which represents the element $\beta(\eta_1) \dots \beta(\eta_t) = 0$ in
$\HH{2t}{G}{k}$. The complex~$\CC$ is obtained by truncating the ends off
of the sequence and the rest of the argument is the same, except that the
integer $t$ has to be replaced by~$2t$ throughout.
\Endproof\vskip4pt  

The upper bounds wanted for the dimensions of our
critical modules is contained in the following.

\begin{thm}\label{uppercriterion}
Suppose that $G$ is a $p$-group which is not elementary abelian.
Suppose that $\eta_1, \dots, \eta_t \in \HH{1}{G}{{\Bbb F}_p}$ are
nonzero and have the property that
\begin{eqnarray*}
\eta_1  \dots \eta_t = 0&& \quad \hbox{\it if $p=2$},
\\
\beta(\eta_1)  \dots \beta(\eta_t) = 0&& \quad \hbox{\it if $p$ is odd.}
\end{eqnarray*}
Let $r=t$ if $p=2$ and $r=2t$ if $p$ is odd.
Let $H_1, \dots, H_r$ be the maximal
subgroups of $G$ as in the previous theorem. Suppose that $M$ is an
indecomposable $kG$-module with the property that
$M{\downarrow}_{H_i}^G \cong k \oplus ({\rm free})$ for every~$i$. Then
for any~$s${\rm ,}
$$
\Dim \ \Omega^s(M) + \Dim \ \Omega^{s-r+1}(M) \leq
\sum_{i = 1}^{r} \Dim \ (\Omega^{s+1-i}(k){\uparrow}_{H_i}^G)\,.
$$
\end{thm}

\Proof 
Let $P$ be a projective module such that $W = k
\oplus \Omega^{1-r}(k) \oplus P$ has a filtration as in the last
theorem. Then tensoring $W$ and all of the
factors in the filtration with $\Omega^s(M)$ we get that
\begin{align*}
\{0\} = L_0 \otimes \Omega^s(M)& \subseteq L_1 \otimes \Omega^s(M)
\subseteq \dots \subseteq L_r \otimes \Omega^s(M) \\
&\cong W \otimes \Omega^s(M) \cong \Omega^s(M)
\oplus \Omega^{s+1-r}(M) \oplus ({\rm free}).
\end{align*}
Then we have 
\begin{align*}
(L_i \otimes \Omega^s(M))/(L_{i-1} \otimes \Omega^s(M))
& \cong (L_i/(L_{i-1}) \otimes \Omega^s(M) \\
& \cong \Omega^{1-i} (k){\uparrow}_{H_i}^G \otimes \Omega^s(M) \\
& \cong \Omega^{1-i} (\Omega^s(M{\downarrow}_{H_i}^G)){\uparrow}_{H_i}^G
\oplus Q \\ & \cong \Omega^{s+1-i} (k){\uparrow}_{H_i}^G \oplus Q'
\end{align*}
for some projective modules $Q$ and $Q'$. Now the important thing
to remember is that $kG$ is a self injective algebra and hence
projective modules are also injective. As a result, if a projective
module is a direct summand of a section of a module~$V$, then
it is a direct summand of~$V$.  The consequence of this is that
(after stripping away the unnecessary projective modules~$Q'$) we can get
that, for some projective module~$R$, the module
$\Omega^s(M) \oplus \Omega^{s+1-r}(M) \oplus R$ has a filtration
$$
\{0\} = X_0 \subseteq X_1 \subseteq \dots \subseteq X_r \cong
\Omega^s(M) \oplus \Omega^{s+1-r}(M) \oplus R
$$
where $X_i/X_{i-1} \cong \Omega^{s+1-i}(k){\uparrow}_{H_i}^G$.
The statement about dimensions follows immediately.
\hfill\qed

\section{Special cases of 2-groups of small order}

In this section we consider some special cases of 2-groups that we need to
treat separately, for they are not covered by the general argument of
Section~10. For each of the groups we show 
Theorem~\ref{mainthm} directly. We discuss the groups of order~8, the almost
extraspecial group
$D_8*C_4$ of order~16 and the extraspecial group
$D_8*D_8$ of order~32 (type~1).

Let us start with the groups of order 8. First $Q_8$ is excluded
by assumption (and there is actually a nontrivial critical $kQ_8$-module
of dimension~5; see~\cite{CT1}). For $G=D_8$ the structure of $T(D_8)$ is
known (see~\cite{CT1}) and every nontrivial endo-trivial $kD_8$-module
is nontrivial on restriction to one of the two elementary abelian
2-subgroups of~$D_8$. Thus the only critical module is the trivial one.
Alternatively, we can also prove the result in the following way.

\begin{prop} \label{d8}
Let $G = D_8$. Then there exists no nontrivial critical $kG$-module.
\end{prop}

\Proof 
Let $M$ be a critical $kG$-module. By Theorem~\ref{noserrelmt}, the number
of cohomology classes whose product vanishes is equal to $t_G=2$.
Applying Theorem~\ref{uppercriterion} with $s=1$, we get
$$
\Dim \ \Omega^1(M) + \Dim \ M \leq
 \Dim \ \Omega^{1}(k_{H_1}){\uparrow}_{H_1}^G
+ \Dim k{\uparrow}_{H_2}^G
$$
for some maximal subgroups $H_1$ and $H_2$. Since $H_i$ has
order~4, $\Omega^{1}(k_{H_i})$ has dimension~3 and we obtain
$$
\Dim \ \Omega^1(M) + \Dim \ M \leq 6 + 2 = 8\,.
$$
By Lemma~\ref{struc-two}, $\Dim \ \Omega^1(M) =
\Dim \ M - 2$. So $\Dim \ M \leq 5$. This part of the argument is
essentially the same as the one appearing in Theorem~5.3 of \pagebreak \cite{CT1}.

If we assume now that there exists a nontrivial critical $kG$-module,
then by Theorem~\ref{lowercriterion}, there exists a nontrivial critical
$kG$-module $M$ of dimension
 $$
\Dim M > \dfrac{\vert G\vert}{2} \cdot
\dfrac{\vert O_G\vert}{\vert C\vert} = 4\cdot 2 = 8 \,,
$$
since $\vert O_G\vert = 2$ by Proposition \ref{ogrouporders}. This
contradicts the previous upper bound. \phantom{endofpage}
\Endproof 

We turn now to the group $G = D_8*D_8$ of order 32.

\begin{prop} \label{e32}
Let $G = D_8*D_8$. Then there exists no nontrivial critical $kG$-module.
\end{prop}

\Proof 
Let $M$ be a critical $kG$-module. By Theorem~\ref{noserrelmt}, the number
of cohomology classes whose product vanishes is equal to $t_G=3$.
Applying now Theorem~\ref{uppercriterion} with $s=1$, we get
$$
\Dim \ \Omega^1(M) + \Dim \ \Omega^{-1}(M) \leq
 \Dim \ \Omega^{1}(k){\uparrow}_{H_1}^G
+ \Dim k{\uparrow}_{H_2}^G + \Dim \ \Omega^{-1}(k){\uparrow}_{H_3}^G
$$
for some maximal subgroups $H_1$, $H_2$, and $H_3$. Since $H_i$ has
order~16, $\Omega^{\pm 1}(k_{H_i})$ has dimension~15 and we obtain
$$
\Dim \ \Omega^1(M) + \Dim \ \Omega^{-1}(M) \leq 30 + 2 + 30 = 62\,,
$$
so that $\Dim \ \Omega^1(M) \leq 62$.

If we assume now that there exists a nontrivial critical $kG$-module,
then by Theorem~\ref{lowercriterion}, there exists a nontrivial critical
$kG$-module $M$ of dimension
 $$
\Dim M > \dfrac{\vert G\vert}{2} \cdot
\dfrac{\vert O_G\vert}{\vert C\vert} \geq 16\cdot \dfrac{72}{9} = 128 \,,
$$
by Proposition \ref{ogrouporders}.
So $\Dim \ \Omega^1(M) > 126$ by Lemma~\ref{struc-two}, a contradiction.
\Endproof 

In the last case, $G$ is the almost extraspecial group $D_8*C_4$ of
order~16. The method of the previous cases does not work because the
orthogonal group~$O_G$ is too small. Instead of using the action of~$O_G$,
we shall give an argument using the action of a Galois group.

\begin{lemma} \label{gps16dims}
Let $G=D_8*C_4$. If $M$ is a critical $kG$-module{\rm ,} then
$\Dim\ M\break \leq 17$. 
\end{lemma}

\Proof  
By Theorem~\ref{noserrelmt}, the number of cohomology classes whose
product vanishes is equal to $t_G=3$. Applying now
Theorem~\ref{uppercriterion} with $s=1$, we get
$$
\Dim \ \Omega^1(M) + \Dim \ \Omega^{-1}(M) \leq
 \Dim \ \Omega^{1}(k){\uparrow}_{H_1}^G
+ \Dim k{\uparrow}_{H_2}^G + \Dim \ \Omega^{-1}(k){\uparrow}_{H_3}^G
$$
for some maximal subgroups $H_1$, $H_2$, and $H_3$. Since $H_i$ has
order~8, $\Omega^{\pm 1}(k)$ has dimension~7 and we obtain
$$
\Dim \ \Omega^1(M) + \Dim \ \Omega^{-1}(M) \leq 14 + 2 + 14 = 30\,.
$$
By Lemma~\ref{struc-two}, $\Dim \ \Omega^1(M) = \Dim \ \Omega^{-1}(M) =
\Dim \ M - 2$. So $\Dim \ \Omega^1(M) \leq 15$ and $\Dim \ M \leq 17$.
\hfill\qed

\begin{prop} \label{ae16}
Let $G = D_8*C_4$. Then there exists no nontrivial critical $kG$-module.
\end{prop}

\Proof 
Suppose that there is such a module $N$.
We need to look at $V_{\overline{G}}(\overline{N}) \subseteq
V_{\overline{G}}(k) \cong k^3$. Suppose that $p = (\alpha,\beta,\gamma)$
is a  point in $V_{\overline{G}}(\overline{N})$. By dividing by
$\alpha$ we may assume that $\alpha =1$, so that
$p = (1,\beta,\gamma) \in V_{\overline{G}}(\overline{N})$.
Notice that $p \notin \ress{G}{H}(V_{\overline{H}}(\overline{N}_H))$
for any maximal subgroup $H$ since $\overline{N}_H$ is a
free $k\overline{H}$-module. Therefore $p$ is not in any
${\Bbb F}_2$-rational subspace of $k^3$, and hence $\beta$ and $\gamma$
cannot both be in the field with four elements (otherwise
$1,\beta,\gamma$ would be linearly dependent over~${\Bbb F}_2$). It
follows
that if $F: k^3 \Rarr{} k^3$ is the Frobenius map,
$F(a,b,c) = (a^2,b^2,c^2)$, then $p$, $F(p)$ and $F^2(p)$
lie on different lines in $\VG$.

Next we need to notice that using the Frobenius homomorphism
we can create a new module from $N$, by letting it act
on the coefficients of the action of the elements of $G$ on $N$.
That is, if the module $N$ is defined by a representation $G \Rarr{}
{\rm GL}(N)$, and if we consider the homomorphism $F:{\rm GL}(N) \Rarr{} {\rm GL}(N)$ that
takes a matrix $(a_{ij})$ to $(a_{ij}^2)$, we let $N^F$ be the module
defined by the composition. It is not difficult to see that
$N^F$ is also critical. Moreover, $F(p)$ is a point in
$V_{\overline{G}}(\overline{N^F})$. It follows that the lines through
$p$, $F(p)$, and $F^2(p)$ are all lines in the variety of the quotient
module $\overline{L}$ for some nontrivial critical module $L$. Thus by
Theorem~\ref{dim-endo}, $kG$ has a nontrivial critical module of
dimension at least~25. This contradicts Lemma~\ref{gps16dims}.
\hfill\qed

\section{The groups of order $p^3$ for odd $p$}

When the prime $p$ is odd, there is one special case in the
proof of Theorem~\ref{mainthm} that must be handled
with extra care. This involves the groups of order $p^3$. The problem
is that the general estimates of the dimensions of critical
modules  used later are not sufficient to handle this case.
The result that we want is the following.

\begin{prop} \label{pcubed}
Let $G = G_1${\rm ,} an extraspecial group of order $p^3$ and exponent~$p${\rm ,} for
$p$ an odd prime. Then there exists no nontrivial critical $kG$-module.
\end{prop}

The proof proceeds in several steps. Throughout assume that a nontrivial
critical $kG$-module exists and use Theorem~\ref{lowercriterion} to obtain
one of large dimension, as follows.

\begin{lemma}\label{lowb-p3}
If a nontrivial critical $kG$-module exists{\rm ,} then there exists a critical
$kG$-module $M$ whose dimension is at least equal to $(p-1)p^4+1$. Moreover
$\Dim\Omega(M) \geq (p-1)p^4-1$ and
$\Dim\Omega^{-1}(M) \geq (p-1)p^4-1$.
\end{lemma}

\Proof 
By Theorem \ref{lowercriterion} there exists a critical module~$M$ whose
dimension is at least $\vert G \vert \ \vert {\rm Sp}(2,\bfp)\vert/ \vert C \vert$
where $C$ is a cyclic $p'$-subgroup of the symplectic group ${\rm Sp}(2,\bfp)$
of maximal order. Now ${\rm Sp}(2,\bfp) = {\rm SL}(2,\bfp)$ has order\break $p(p^2-1)$
and its cyclic $p'$-subgroup of maximal order has order~$p+1$. So the
dimension of $M$ must be greater than $(p-1)p^4$
and must be congruent to~1 modulo~$p$.

Now to compute the dimension of $\Omega(M)$, we notice from the proof
of Theorem~\ref{lowercriterion} that the variety of the
module $\overline{M}$ is the union of at least $p(p-1)$ distinct lines
in $V_{\overline{G}}(k) = k^2$. So $\overline{M} = \overline{U}_1 \oplus
\dots \oplus \overline{U}_t$ where, for each $i$,
$V_{\overline{G}}(\overline{U}_i)$ is a single
line and $t > p(p-1)$. Now, as in the proof of Theorem~\ref{dim-endo},
$\Dim \overline{U}_i = r_ip^2$ for some~$r_i$ (we use here the fact that
$U_i$ is endo-trivial and $p$ is odd). Because $\overline{U}_i$ is not a
free $k\overline{G}$-module
(and, in fact, has no free submodules) and because a projective cover
of~$\overline{U}_i$ has dimension $p^2  \Dim
\overline{U}_i/{\rm Rad}(\overline{U}_i)$, we must have
$\Dim \overline{U}_i/{\rm Rad}(\overline{U}_i) > r_i$.
Therefore $\overline{U}_i$ is minimally generated by at least $r_i+1$
generators and the number of generators of $\overline{M}$ is at least
$$m=\sum_{i = 1}^t (r_i+1)= \left(\sum_{i=1}^t r_i\right) +t \,.$$
Now $\overline{M}/{\rm Rad}(\overline{M})$ is a quotient of $M/{\rm Rad}(M)$, so the
minimal number of generators of $M$ is at least $m$.
As a result, the number of copies of $kG$ appearing in the projective cover
of $M$ must be at least $m$. Now the
dimension of $M$ is $p^3(\sum_{i=1}^t r_i) + 1$ and so the dimension of
$\Omega(M)$ is at least
$$p^3m- \Dim(M) = p^3\left(\left(\sum_{i=1}^t r_i\right) +t\right) - p^3\left(\sum_{i=1}^t r_i\right) - 1 =
tp^3-1 \geq (p-1)p^4-1 \,.$$
By applying the same argument to the dual module $M^*$ (which also satisfies
the properties we need), we obtain
$$\Dim\Omega^{-1}(M) = \Dim\Omega^{-1}(M)^* = \Dim\Omega(M^*) \geq
(p-1)p^4-1 \,.$$
This proves the lemma.
\hfill\qed

\begin{lemma}\label{upper-p3}
$$
\Dim \Omega^{2p}(M) + \Dim \Omega^{-1}(M) \leq p^3(p^2+p+1).
$$
\end{lemma}

\Proof 
From any one of the papers \cite{Lear1}, \cite{Yal}, \cite{BCexp} we have that
there exist $\eta_1, \dots \eta_{p+1} \in \HH{1}{G}{k}$ such that
$\beta(\eta_1) \dots \beta(\eta_{p+1}) = 0$. In Leary \cite{Lear1}
the relation is given as $x^px' -x{x'}^p =0$. Now applying
Theorem~\ref{uppercriterion} with $t=p+1$ (hence $r=2t=2(p+1)$) and choosing
$s=2p$ in that theorem, we get
$$
\Dim \Omega^{2p}(M) + \Dim \Omega^{-1}(M) \leq \sum_{i=1}^{2p+2} \Dim (\Omega^{2p+1-i}(k){\uparrow}_{H_i}^G) \,,
$$
where $H_i$ is a maximal subgroup of $G$ corresponding to the
appropriate~$\eta_j$. In our case, every $H_i$ is an elementary
abelian group of order $p^2$, and hence the dimensions on the right-hand
side of the inequality are independent of the particular~$\eta_j$.
Because $\Dim \HH{j}{H_i}{k} = j+1$ (see Lemma~\ref{dimcoho-elem}), we have
that (for $H_i = H$)
$$
\Dim \Omega^{2j-1}(k_H) + \Dim \Omega^{2j}(k_H) =
p^2\Dim \HH{2j-1}{H}{k} = p^2 (2j) \,.
$$
Induction to $G$ multiplies the dimensions by $p$.
Consequently the right-hand side of the above inequality has the form
\begin{align*}
\sum_{i=1}^{2p+2} \Dim &(\Omega^{2p+1-i}(k){\uparrow}_{H_i}^G) \\
& = p\Dim \Omega^{-1}(k_H) + p\Dim k
+ p\sum_{j=1}^p \big(\Dim \Omega^{2j-1}(k_H)
+ \Dim \Omega^{2j}(k_H)\big) \\
& = p\big(p^2-1+1+ \sum_{j=1}^p 2p^2j\,\big) \\
& = p^3 + 2p^3(p)(p+1)/2 = p^3(1 + p^2+ p)
\end{align*}
as desired.
\Endproof\vskip4pt  

At this point we should notice that the two lemmas above are not
sufficient to give us the contradiction wanted. We need some
further analysis of the dimension of $\Omega^{2p}(M)$. For this
purpose we  recall that there exists an element $\zeta
\in \HH{2p}{G}{k}$ which has the property that its restriction
$\res{G}{Z}(\zeta)$ is not zero where $Z = \langle z \rangle$ is the center
of~$G$. In Leary's paper \cite{Lear1}, the element that he calls $z$ will
do. The element $\zeta$ can also be obtained by applying the Evens norm map
to an element in the degree~2 cohomology of a maximal elementary abelian
subgroup whose restriction to $Z$ is not trivial.

The element $\zeta$ can be represented by a unique
cocycle $\zeta: \Omega^{2p}(k) \Rarr{} k$. Hence we have an
exact sequence 
$$
0 \Rarr{} L \Rarr{} \Omega^{2p}(k) \Rarr{\zeta} k \Rarr{} 0
$$
where $L$ is the kernel of $\zeta$. Now by Theorem \ref{propvar},
$V_G(L) = V_G(\zeta)$, the variety of the ideal generated by $\zeta$.
In particular, the restriction  $L{\downarrow}_Z^G$
is free as a $kZ$-module. This fact can also be derived from
the observation that the above sequence is split as a sequence
of $kZ$-modules because the restriction of $\zeta$ to $Z$ is not
zero and $\Omega^{2p}(k){\downarrow}_Z^G \cong k \oplus ({\rm free})$.

Let $\overline{L} = L/(z-1)L \cong (z-1)^{p-1}L$. Then $\overline{L}$
is a $k\overline{G}$-module where $\overline{G} = G/Z$.

\begin{lemma} \label{varLbar}
The $k\overline{G}$-module $\overline{L}$ has no projective submodules{\rm ,} and
moreover{\rm , }
$$
V_{\overline{G}}(\overline{L}) \subseteq \;
\bigcup \, \ress{\overline{G}}{\overline{E}} V_{\overline{E}}(k)
$$ 
where the union is over the set of all subgroups
$\overline{E} = E/Z$ where $E$ is a maximal subgroup of $G$.
\end{lemma}

Notice that every maximal subgroup of $G$ is elementary abelian and the
union in the lemma is over all subgroups of order $p$ in
$\overline{G}$. Thus the
right-hand side of the containment is the union of all of the $\bfp$-rational lines in $V_{\overline{G}}(k)
\cong k^2$. It can be proved that the two sides are actually equal, but we do not need to know this.

\Proof  
If $\overline{L}$ had a $k\overline{G}$-projective submodule then
$L$ and hence also $\Omega^{2p}(k)$ would have projective $kG$-submodules.
That is, if $t_1^{\overline{G}}\overline{L} \neq 0$ then also
$t_1^{G}L \neq 0$. But clearly this is impossible.

Now suppose that $\ell \subseteq V_{\overline{G}}(k)$ is a line that
is not $\bfp$-rational. Let $N$ be a $k\overline{G}$-module such that
$V_{\overline{G}}(N) = \ell$ (e.g.\ take $N=k\overline{G}/ (\sigma -1)$ where
$\langle \sigma \rangle$ is a cyclic shifted subgroup corresponding to the
line~$\ell$). Then the restriction
$N{\downarrow}_{\overline{E}}^{\overline{G}}$
is a free $k\overline{E}$-module for any maximal subgroup $E$ of~$G$.
So, viewing $N$ as a $kG$-module by inflation, we have that
$V_E(N{\downarrow}_E^G)$ is the line determined by the center~$Z$, because
$Z$ acts trivially on~$N{\downarrow}_E^G$.
Therefore $N$ is periodic as a $kG$-module and we must have that $V_G(N) =
\ress{G}{Z}(V_Z(k))$,
the line determined by the center~$Z$. Because $L$ is free
on restriction to $Z$ we know that $V_G(L) \cap V_G(N) = \{0\}$ and
hence $L \otimes N$ is a free $kG$-module. Now $Z$ acts trivially on
$N$ and hence $(z-1)(L \otimes N) = ((z-1)L) \otimes N$. Thus,
$\overline{L \otimes N} \cong \overline{L} \otimes N$ is a free
$k\overline{G}$-module. It follows from Theorem~\ref{propvar} that
$V_{\overline{G}}(\overline{L}) \cap V_{\overline{G}}(N) = \{0\}$.
Hence the line $\ell$ is not in $V_{\overline{G}}(\overline{L})$ and
this holds for all lines in $V_{\overline{G}}(k)$ which are not
$\bfp$-rational. Thus the variety $V_{\overline{G}}(\overline{L})$
must be contained in the union of the $\bfp$-rational lines.
\hfill\qed

\begin{lemma} \label{LtensonM}
If $M$ is a critical $kG$-module{\rm ,} $V_{\overline{G}}(\overline{M}) \cap
V_{\overline{G}}(\overline{L}) = \{0\}$ and $M \otimes L \ \cong L \ \oplus
({\rm free})$. 
\end{lemma}

\Proof 
We first show that $\overline{M} \otimes L$ is a free $kG$-module.
That is, $L$ is free as a $kZ$-module and $\overline{\overline{M}
\otimes L} \cong \overline{M} \otimes \overline{L}$. But from Lemmas
\ref{varLbar} and \ref{fprationality} we have that
$V_{\overline{G}}(\overline{M}) \cap V_{\overline{G}}(\overline{L})
= \{0\}$. Hence $\overline{M} \otimes \overline{L}$ is free
as a $k\overline{G}$-module. Thus $\overline{M} \otimes L$ is free
as a $kG$-module.

It follows that $M \otimes L$ has a filtration
$$
0 \subseteq ((z-1)^{p-1}M)\otimes L \subseteq  \dots \subseteq
((z-1)M)\otimes L \subseteq M' \otimes L \subseteq M \otimes L
$$ 
where $M' = \{m \in M \vert (z-1)^{p-1}m = 0\}$. All of the factors
are isomorphic to $\overline{M} \otimes L$ and hence are projective,
except for the factor
$$(M' \otimes L)/((z-1)M \otimes L) \cong (M'/(z-1)M)
\otimes L \cong k \otimes L \cong L \,.$$
The lemma follows from the fact that free modules are also injective and
hence any free composition factor is a direct summand.
\Endproof\vskip4pt  

Now tensoring the sequence given above with $M$ we get an exact
sequence
$$
0 \Rarr{} M \otimes L \Rarr{} M \otimes \Omega^{2p}(k)
\Rarr{1 \otimes \zeta} M \Rarr{} 0 \,.
$$
Any projective submodule of $M\otimes L$ is also a direct summand of
the middle term and can be factored out. So we have an exact sequence
of the form
$$
0 \Rarr{} L \Rarr{} \Omega^{2p}(M) \oplus P  \Rarr{} M \Rarr{} 0 \,,
$$
for some projective module $P$. It remains to prove the following.

\begin{lemma} \label{Peq0}
In the preceding exact sequence{\rm ,} the projective module $P$ is zero.
\end{lemma}

\Proof 
Because the module $L$ is free as a $kZ$-module the sequence is split
as a sequence of $kZ$-modules. So multiplication by $z-1$ is an exact
functor on this sequence. Hence we have a sequence
$$
0 \Rarr{} (z-1)^{p-1}L \Rarr{} (z-1)^{p-1}\Omega^{2p}(M) \oplus (z-1)^{p-1}P
\Rarr{} (z-1)^{p-1}M \Rarr{} 0 \,;
$$
that is,
$$
0 \Rarr{} \overline{L} \Rarr{} \overline{\Omega^{2p}(M)} \oplus \overline{P}
\Rarr{} \overline{M} \Rarr{} 0 \,,
$$
which is a sequence of $k\overline{G}$-modules.
Because $V_{\overline{G}}(\overline{L}) \cap V_{\overline{G}}(\overline{M})
= \{0\}$ by the previous lemma, we must have that the sequence splits. Thus,
$$ 
\overline{L} \oplus \overline{M} \cong \overline{\Omega^{2p}(M)} \oplus
\overline{P} \,.
$$
But $\overline{L} \oplus \overline{M}$ has no projective 
$k\overline{G}$-submodules by
Lemma~\ref{varLbar}. Hence $\overline{P} = \{0\}$ and therefore also $P =
\{0\}$.
\Endproof\vskip4pt  

  {\it Proof of Proposition} \ref{pcubed}.
By Lemma \ref{Peq0},  $\Dim \Omega^{2p}(M) = \Dim L + \Dim M$.
By Lemma~\ref{dim-omega-2p}, $\Dim \Omega^{2p}(k) = p^3(p+1)+1$, and so
$\Dim L = p^3(p+1)$ by definition of~$L$. Now by Lemma~\ref{lowb-p3},
$\Dim M \geq p^4(p-1)+1$ and $\Dim\Omega^{-1}(M) \geq p^4(p-1)-1$.
Hence we have that 
\begin{eqnarray*}
\Dim \Omega^{2p}(M) + \Dim \Omega^{-1}(M)& \geq& p^3(p+1) + p^4(p-1)+1 +
p^4(p-1)-1 \\
&=& p^3(2p^2-p+1) \,.
\end{eqnarray*}
This inequality, however, is a contradiction
to Lemma \ref{upper-p3} since we are assuming that $p \geq 3$.

\section{The general case in characteristic 2}

We are now prepared to prove the general case by induction and
complete the proof of the detection Theorem~\ref{mainthm} when $p=2$.
Throughout, $k$ has characteristic~2. Let $G$ be an
extraspecial or almost extraspecial group of order $2^{m+1}$.
The theorem that we are trying to prove is the following. It is equivalent
to Theorem~\ref{mainthm}.

\begin{thm} \label{mainthm2}
If $G$ is an extraspecial or almost extraspecial $2$-group and if
$G$ is not isomorphic to $Q_8${\rm ,} then there are no nontrivial critical
$kG$-modules. 
\end{thm}

Three cases have to be treated separately, namely the
groups of order at most~16 as well as~$D_8*D_8$. But these cases have
been dealt with in Section~8. Therefore we can now assume
that $m\geq 4$ and that $m>4$ for the groups of type~1. This allows us
to use Corollary~\ref{dimsumomeg}.

The strategy of the proof is expressed in the following.

\begin{prop} \label{strategy}
Let $G$ be an extraspecial or almost extraspecial group of
order~$2^{m+1}${\rm ,} with $m=2n$. Assume that $m\geq 4$ and $m>4$ if $G$ is of
type~$1$. Let $t_G$ be the number of cohomology classes whose product
vanishes{\rm ,} as described in Theorem~{\rm \ref{noserrelmt},} and let
$$
\sigma_G = \binom{t_G+m-4}{m-2}\vert G \vert+2 \qquad \text{and} \qquad
\tau_G = \dfrac{\vert G \vert}{2}  \cdot
\dfrac{\vert O_G \vert}{3^n} \,.
$$ 
If $\tau_G > \sigma_G$ then there exists no nontrivial critical
$kG$-module.
\end{prop}

\Proof  Let $t=t_G$. In view of Theorem~\ref{noserrelmt}, there exist 
nonzero elements $\eta_1, \dots, \eta_t\in \HH{1}{G}{{\Bbb F}_2}$ such
that $\eta_1\dots\eta_t=0$ and each $\eta_i$ corresponds to a maximal
subgroup~$H_i$. Moreover each subgroup $H_i$ is the
centralizer of a noncentral involution in~$G$ and by
Theorem~\ref{dimcoho}, $H_i \cong C_2 \times U$ where $U$ has the same
type as~$G$. So $H_i\cong H_1$ for each~$i$.

Suppose that $M$ is a critical $kG$-module. Then by Theorem
\ref{uppercriterion} with $t=t_G$ and $s=t-1$, we have
$$
\Dim \ M \leq \Dim \ \Omega^{t-1}(M) + \Dim \ M \leq \sum_{i = 1}^{t} \Dim \ (\Omega^{t-i}(k){\uparrow}_{H_i}^G)\,.
$$
Since all the subgroups $H_i$ are isomorphic to $H_1$, we obtain
$$
\Dim \ M \leq \sum_{j = 0}^{t-1} \Dim \
(\Omega^{j}(k){\uparrow}_{H_1}^G)\,.
$$
 Now by Corollary \ref{dimsumomeg}, which
applies in view of our assumption on~$m$ (with $m$, in the corollary,
replaced by $m-2$ and $r = t_G-1$), we obtain
$$
\sum_{j = 0}^{t-1} \Dim \ (\Omega^{j}(k){\uparrow}_{H_1}^G) \leq \vert G \vert \cdot \binom{m-2+t_G-1-1}{m-2} + 2 = \sigma_G \,.
$$
It follows that $\Dim M\leq \sigma_G$.

If  there exists a nontrivial critical $kG$-module,
then by Theorem~\ref{lowercriterion}, there exists a nontrivial critical
$kG$-module $M$ of dimension
 $$
\Dim \ M> \dfrac{\vert G\vert}{2} \cdot
\dfrac{\vert O_G\vert}{\vert C\vert}
\ \geq \ \dfrac{\vert G\vert}{2} \cdot
\dfrac{\vert O_G\vert}{3^n}
= \tau_G> \sigma_G\,.
$$
This contradicts the upper-bound obtained above.
\Endproof\vskip4pt  

We have now reduced the problem to the proof that $\tau_G > \sigma_G$ for
all the groups~$G$ as above. This is a purely numerical problem which
only requires estimating the numbers $\tau_G$ and~$\sigma_G$. We
start with a lemma which will be useful for estimating~$\sigma_G$.

\begin{lemma} \label{estimation}
Let $t$ and $m$ be integers with $t\geq 4$ and $m\geq 6$. Then
$$
\dfrac{\displaystyle\binom{2t+m-2}{m}} {\displaystyle\binom{t+m-4}{m-2}}
< 2^{m-3} \ t^2 \,.
$$
\end{lemma}

\Proof 
Expanding the left-hand side and eliminating the common factor~$(m-2)!$,
we get the following expression. Notice that we can bound each of the first
$m{-}5$ fractions by~2 (using $m\geq 6$), the next three by~3 (using
$t\geq 4$), and bound $1/m(m-1)$ by~$1/30$. Thus we get the following.

\begin{gather*}
  \dfrac{2t+m-2}{t+m-4} \cdot \dfrac{2t+m-3}{t+m-5} \cdot \ldots
 \cdot \dfrac{2t+4}{t+2} \cdot \dfrac{2t+3}{t+1} \cdot \dfrac{2t+2}{t}
\cdot
\dfrac{2t+1}{t-1} \cdot \dfrac{2t}{m} \cdot \dfrac{2t-1}{m-1} \\
 < 2^{m-5} \, 3^3 \, \dfrac{4t^2}{30} < 2^{m-3} t^2\,.
\end{gather*}
\vglue-20pt
\Endproof\vskip12pt  

For the proof that $\tau_G > \sigma_G$, we proceed with cases.

\Subsec{Groups of type $1$}
Let $G_n = D_8* \dots *D_8$ be the central product of $n$ copies
of~$D_8$, with $n\geq 3$. Remember that the cases $n = 1$ and $n=2$ were
treated in Propositions~\ref{d8} and~\ref{e32}. For convenience, we write
$G = G_n$ and let
$\sigma_n = \sigma_{G_n}$ and $\tau_n = \tau_{G_n}$.
We prove that $\sigma_n < \tau_n$ by induction, starting with two cases.

If $n=3$, then $t=t_G=5$ by Theorem~\ref{noserrelmt} and we have that
$$ \tau_3 = 2^6\cdot \dfrac{2^7\cdot7\cdot3\cdot15}{3^3}> 2^{13}\cdot11 \,,
$$
$$
\sigma_3 = \binom{5+6-4}{6-2}\cdot 2^7+2 = 35\cdot2^7+2<
2^{13}\cdot11< \tau_3\,.$$
 If $n=4$, then $t=t_G=9$ by Theorem~\ref{noserrelmt} and
\begin{eqnarray*}
 \tau_4&=&   2^8\cdot \dfrac{2^{13}\cdot15\cdot3\cdot15\cdot63}{3^4}
= 2^{21}\cdot525 \,,
\\
\sigma_4&=& \binom{9+8-4}{8-2}\cdot 2^9+2 = 1716\cdot2^9+2<
\tau_4\,.
\end{eqnarray*}
 For $n\geq 4$, we have $t_{G_{n+1}}=2t_{G_n}$ by
Theorem~\ref{noserrelmt}, and this allows for an inductive argument. We
assume that $\sigma_n <
\tau_n$ and we prove that $\sigma_{n+1} < \tau_{n+1}$.
The course of our proof is to show that
$$
\dfrac{\sigma_{n+1}-2}{\sigma_n-2}<
2^{4n}< \dfrac{\tau_{n+1}}{\tau_n}
$$
from which we get $\sigma_{n+1} - 2 < 2^{4n} \sigma_n - 2^{4n+1} <
2^{4n} \tau_n - 2 < \tau_{n+1} - 2$ and we are done. So we are left with
the proof of the two inequalities above.

From the value of $\tau_n$ given by Proposition \ref{ogrouporders}, we
obtain
\begin{eqnarray*}
\dfrac{\tau_{n+1}}{\tau_n}   &=&
  \dfrac{\vert G_{n+1}\vert}{\vert G_n\vert}\cdot
\dfrac{2^{(n+1)n+1}}{2^{n(n-1)+1}}\cdot\dfrac{2^{n+1}-1}{2^n-1}\cdot
\dfrac{2^{2n}-1}{3} \\
&>&    2^2\cdot 2^{2n} \cdot 2 \cdot \dfrac{2^{2n}-1}{4}
= 2^{4n+1}  - 2^{2n+1}> 2^{4n} \,.
\end{eqnarray*}
On the other hand, setting $m=2n$ and $t_{G_n} = t_n = 2^{n-1}+2^{n-4}$
(Theorem~\ref{noserrelmt}), we obtain by Lemma~\ref{estimation}
\begin{eqnarray*}
\dfrac{\sigma_{n+1}-2}{\sigma_n-2}    & = &
\dfrac{\vert G_{n+1}\vert}{\vert G_n\vert} \cdot
\dfrac{\displaystyle\binom{2t_n+m-2}{m}}
{\displaystyle\binom{t_n+m-4}{m-2}} \\
 &  <&
4\cdot 2^{m-3}\cdot t_n^2
= 2^{2n-1} (2^{n-1}+2^{n-4})^2< 2^{2n-1}\, 2^{2n}< 2^{4n} \,.
\end{eqnarray*}

\Subsec{Groups of type $2$}
Let $G_n = D_8* \dots *D_8*Q_8$ be the central product
of $n-1$ copies of $D_8$ and one of~$Q_8$, with $n\geq 2$. Let $G = G_n$,
$\sigma_n = \sigma_{G_n}$, $\tau_n = \tau_{G_n}$, and $t_n = t_{G_n} =
2^n + 2^{n-2}$ (see Theorem~\ref{noserrelmt}).

We start with the case $n=2$, for which we need to replace $\tau_2$ by
the slightly larger value
 $$ \tau_2' = \dfrac{\vert G_2 \vert}{2}  \cdot
\dfrac{\vert O_G \vert}{\vert C \vert}$$
where $C$ is a cyclic subgroup of~$O_G$ of maximal odd order. By
Corollary~12.43 of Taylor's book~\cite{Tay}, $O_G$ has a simple subgroup
of index~2 isomorphic\break to ${\rm PSL}(2,{\Bbb F}_4)$ (that is, $A_5$, and in fact
$O_G$ is isomorphic to the symmetric group~$S_5$). Therefore $\vert C
\vert = 5$ and we get $\tau_2' = 384$. On the other hand $t_2 = 5$ and we
have
$$\sigma_2 = \binom{5+4-4}{4-2}\cdot 2^5+2 = 322<
\tau_2'\,.$$
 The argument of Proposition~\ref{strategy} goes through with $\tau_2'$
instead of~$\tau_2$.

Now we prove that $\sigma_n < \tau_n$ by induction, starting with $n=3$:
\begin{eqnarray*}
 \tau_3 & =& 2^6\cdot \dfrac{2^7\cdot9\cdot3\cdot15}{3^3}
 = 2^{13}\cdot15 \,,
\\
\sigma_3& =& \binom{10+6-4}{6-2}\cdot 2^7+2 = 495\cdot2^7+2<
\tau_3\,.
\end{eqnarray*}
 If now $n\geq 3$ the course of our proof is to show that
$$
\dfrac{\sigma_{n+1}-2}{\sigma_n-2}  < 
2^{4n} - 2^{2n}< \dfrac{\tau_{n+1}}{\tau_n}
$$
from which we conclude the proof as in the previous case.
Here is the computation:
\begin{eqnarray*}
\dfrac{\tau_{n+1}}{\tau_n}  &=& \dfrac{\vert G_{n+1}\vert}{\vert G_n\vert}\cdot
\dfrac{2^{(n+1)n+1}}{2^{n(n-1)+1}}\cdot\dfrac{2^{n+1}+1}{2^n+1}\cdot
\dfrac{2^{2n}-1}{3} \\
 & >& 2^2\cdot 2^{2n} \cdot 1 \cdot \dfrac{2^{2n}-1}{4}
 = 2^{4n}  - 2^{2n} \,.
\end{eqnarray*}
On the other hand, we obtain by Lemma~\ref{estimation}
\begin{eqnarray*}
\dfrac{\sigma_{n+1} - 2}{\sigma_n - 2}   & =&
\dfrac{\vert G_{n+1}\vert}{\vert G_n\vert} \cdot
\dfrac{\displaystyle\binom{2t_n+m-2}{m}}
{\displaystyle\binom{t_n+m-4}{m-2}} \\
 & <&
4\cdot 2^{m-3}\cdot t_n^2
 = 2^{2n-1} (2^{n}+2^{n-2})^2 \\
 & =& 2^{4n-1} + 2^{4n-2} + 2^{4n-5} + 2^{2n}- 2^{2n}< 2^{4n}- 2^{2n} \,.
\end{eqnarray*}

\Subsec{Groups of type $3$}
Let $G_n = D_8* \dots *D_8*C_4$ be the central product
of $n$ copies of $D_8$ and one of~$C_4$. Let $G = G_n$,
$\sigma_n = \sigma_{G_n}$, $\tau_n = \tau_{G_n}$, and $t_n = t_{G_n} =
2^n + 2^{n-2}$ (see Theorem~\ref{noserrelmt}). Note that $m = 2n + 1$ for
type~3.

We prove that $\sigma_n < \tau_n$ by induction, starting with $n=2$.
Remember that the case in which $n = 1$ was treated in
Proposition~\ref{ae16}. First we have that
\begin{eqnarray*}
 \tau_2 & = & 2^5\cdot \dfrac{2^4\cdot3\cdot15}{3^2}
= 2560 \,,
\\
\sigma_2 &=&  \binom{5+5-4}{5-2}\cdot 2^6+2 = 1282< \tau_2\,.
\end{eqnarray*}
If now $n\geq 2$ we show that
$$
\dfrac{\sigma_{n+1}-2}{\sigma_n-2}   < 
2^{4n+2}  <  \dfrac{\tau_{n+1}}{\tau_n}
$$
from which we make our conclusions as in the previous cases.
Here is the computation. First note that
\begin{eqnarray*}
\dfrac{\tau_{n+1}}{\tau_n}   &=
& \dfrac{\vert G_{n+1}\vert}{\vert G_n\vert}\cdot
\dfrac{2^{(n+1)^2}}{2^{n^2}}\cdot
\dfrac{2^{2(n+1)}-1}{3} \\
 & >& 2^2\cdot 2^{2n+1} \cdot \dfrac{2^{2n+1}}{4}
 = 2^{4n+2} \,.
\end{eqnarray*}
On the other hand, we obtain by Lemma~\ref{estimation}
\begin{eqnarray*}
\dfrac{\sigma_{n+1} - 2}{\sigma_n - 2}   & =&
\dfrac{\vert G_{n+1}\vert}{\vert G_n\vert} \cdot
\dfrac{\displaystyle\binom{2t_n+m-2}{m}}
{\displaystyle\binom{t_n+m-4}{m-2}} \\
 & < &
4\cdot 2^{m-3}\cdot t_n^2
   =   2^{2n} (2^{n}+2^{n-2})^2
  <  2^{2n}\, 2^{2(n+1)}   =  2^{4n+2} \,.
\end{eqnarray*}
This completes the proof of Theorem~\ref{mainthm2} and hence also the
proof of Theorem~\ref{mainthm} when $p=2$.

\section{The general case in odd characteristic}

In this section we complete the proof of Theorem \ref{mainthm} for
odd~$p$.  We assume throughout that the field $k$ has characteristic
$p$ and that $G = G_n$ is an extraspecial group of order $p^{2n+1}$ and
exponent~$p$. Our aim is to prove the following.

\begin{thm} \label{mainthm-p}
If $G = G_n${\rm ,} then there are no nontrivial critical $kG$-modules.
\end{thm}

If $n=1$, the result follows from Section~9. Thus we can assume $n\geq 2$.
The proof follows the same basic pattern as in the last section.
We define $\sigma_n$ and $\tau_n$ such that $\sigma_n$ is an
upper bound for the dimension of any critical module and $\tau_n$
is a lower bound for the dimension of some nontrivial critical
module if nontrivial critical modules exist. Then we prove that
$\sigma_n < \tau_n$. First we give the definitions.

For $n \geq 2$ let 
$$
\sigma_n = 2\vert G_n \vert \binom{t_n +2n -3}{2n-1}
$$
where $t_n = 2(p^2+p-1)p^{n-2}$. Let $\tau_n$ be given by the rule
$$
\tau_n = \dfrac{\vert G_n \vert \ \vert {\rm Sp}(2n,\bfp) \vert}{c_n}
$$ 
where $c_n = (p+1)^n$ except in the case in which $p=3$ and
$n=2$. In that case let $c_n = p^2+1 = 10$. Then we have the
following.

\begin{prop} \label{strategy-p}
If $n \geq 2$ and $\tau_n > \sigma_n${\rm ,} then there exists no
nontrivial critical $kG_n$-module.
\end{prop}

\Proof 
Let $t=t_n/2= (p^2+p-1)p^{n-2}$. By Theorem~\ref{noserrelmt-p}, we know that
there exist  nonzero elements $\eta_1, \dots, \eta_t\in \HH{1}{G}{{\Bbb
F}_2}$ such that $\beta(\eta_1)\dots\beta(\eta_t)\break
=0$ and each $\eta_i$
corresponds to a maximal subgroup~$H_i$. Moreover each subgroup
$H_i$ is the centralizer of a noncentral element of order~$p$ in~$G$ and by
Theorem~\ref{dimcoho-p}, $H_i \cong C_p \times G_{n-1}$. So $H_i\cong H_1$
for each~$i$.

If $M$ is a critical $kG_n$-module, then by Theorem~\ref{uppercriterion}
with $r=2t=t_n$ and $s=t_n-1$,  
$$
\Dim \ M \leq \Dim \ \Omega^{t_n-1}(M) + \Dim \ M \leq \sum_{i = 1}^{t_n} \Dim \ (\Omega^{t_n-i}(k){\uparrow}_{H_i}^G)\,.
$$
Since all the subgroups $H_i$ are isomorphic to $H_1$, we obtain
$$
\Dim \ M \leq \sum_{j = 0}^{t_n-1} \Dim \
(\Omega^{j}(k){\uparrow}_{H_1}^G)\,.
$$
By Corollary~\ref{dimsumomeg-p},  
$$
\sum_{j = 0}^{t_n-1} \Dim \ (\Omega^{j}(k){\uparrow}_{H_1}^G) \leq 2\vert G \vert \cdot \binom{t_n-1+2n-2}{2n-1}= \ \sigma_n \,.
$$
It follows that $\Dim M\leq \sigma_n$.

On the other hand, if we assume that there exists a nontrivial
critical $kG_n$-module, then by Theorem~\ref{lowercriterion}, there
exists a nontrivial critical $kG$-module $M$ of dimension
 $$
\Dim \ M> \vert G\vert \cdot
\dfrac{\vert O_G\vert}{\vert C\vert}
$$
where $C$ is a cyclic $p'$-subgroup in ${\rm Sp}(2n,\bfp)$ of maximal
order. In the case that $p = 3$ and $n =2$, we know from character
tables or from direct analysis on ${\rm Sp}(4,{\Bbb F}_3)$ that $C$ has
order at most~10. In all other cases we know by
Proposition~\ref{ogrouporders} that the order of $C$ is at most $(p+1)^n$.
So in either case, $\Dim M > \tau_n$. Hence if $\sigma_n < \tau_n$
then we have a contradiction.
\Endproof\vskip4pt  

So it remains to prove that $\tau_n > \sigma_n$. We will proceed
by induction beginning with the following.

\begin{lemma}
$\tau_2 > \sigma_2$.
\end{lemma}

\Proof  
If $p=3$, then $\sigma_2 = 860,706$ while $\tau_2 = 1,259,712$, and so
the lemma holds in that case (note that it is here that we need the special
choice made for $\tau_2$). So suppose that $p \geq 5$. Then
\begin{align*}
\dfrac{\sigma_2}{p^{11}} \ &= \
\dfrac{1}{p^{11}}2p^5  \binom{2(p^2+p-1)+1}{3} \\
&= \ \dfrac{2}{3!}(2+\dfrac{2}{p}-\dfrac{1}{p^2})
(2+\dfrac{2}{p}-\dfrac{2}{p^2})(2+\dfrac{2}{p}-\dfrac{3}{p^2})
< \dfrac{1}{3}3 \cdot 3 \cdot 3 = 9 \,.
\end{align*}
On the other hand,
$$
\dfrac{\tau_2}{p^{11}} =
\dfrac{p^5}{p^{11}} \dfrac{p^4(p^2-1)(p^4-1)}{(p+1)^2}
 = p^2 \dfrac{(1-\dfrac{1}{p^2})(1 - \dfrac{1}{p^4})}
{(1+\dfrac{1}{p})^2}> p^2
\dfrac{(\dfrac{4}{5})^2}{(\dfrac{4}{3})^2} = 9 \dfrac{p^2}{25}> 9
$$
by the fact that $p \geq 5$ and hence $1+ 1/p <4/3$ and $1-1/p \geq 4/5$.
So again $\tau_2 > \sigma_2$.
\hfill\qed

\begin{lemma} 
For $n \geq 2${\rm ,} $\dfrac{\sigma_{n+1}}{\sigma_n} <
\dfrac{\tau_{n+1}}{\tau_n}$.
\end{lemma}

\Proof 
Notice first that a special computation is needed if $p=3$ and $n=2$. In
that case, by direct calculation, we
have that $\sigma_3 = 49,157,255,862$ while $\tau_3 = 313,380,128,880$.
It is then easy to check the lemma in this particular case.

More generally, we calculate that
\begin{align*} 
\dfrac{\tau_{n+1}}{\tau_n} & = \dfrac{\vert G_{n+1} \vert}{ \vert G_n \vert}
\dfrac{p^{(n+1)^2}}{p^{n^2}}
\dfrac{(p^2-1) \dots (p^{2n+2}-1)/(p+1)^{n+1}}
{(p^2-1) \dots (p^{2n}-1)/(p+1)^n} \\
&= p^2\cdot p^{2n+1} \cdot (p^{2n+2}-1)/(p+1)> \dfrac{1}{2} p^{4n+4}.
\end{align*}
The above estimate is that, since $p \geq 3$, we have
$1/(p+1) > 1/(\sqrt{2}p)$ and $p^{2n+2} -1 > p^{2n+2}/\sqrt{2}$.

At the same time, setting $t = t_n$ and noting that $t_{n+1}=pt_n$, we have
\begin{align*}
\dfrac{\sigma_{n+1}}{\sigma_n} &= \dfrac{2p^{2n+3}}{2p^{2n+1}}
\dfrac{\binom{tp+2n-1}{2n+1}}{\binom{t+2n-3}{2n-1}} \\
&= \dfrac{p^2}{(2n+1)(2n)} (tp+2n-1)(tp+2n-2)\dfrac{(tp+2n-3)}{(t+2n-3)}
\dots \dfrac{tp}{t}\dfrac{(tp-1)}{(t-1)} \,.
\end{align*}
Now we note that $(tp+b)/(t+b) \leq tp/t = p$ for all~$b\geq 0$.
Also $(tp-1)/(t-1)< \dfrac{3}{2}p$ because $t\geq 3$.
Moreover,
\begin{align*}
tp + 2n -1 & = 2(p^2 + p -1)p^{n-2}p + 2n-1 \\
& = 2p^{n+1}\left( 1 + \dfrac{1}{p} - \dfrac{1}{p^2} +
\dfrac{2n-1}{2p^{n+1}}
\right) \\
& < 2p^{n+1}(2) = 4p^{n+1}.
\end{align*}
So we have that 
\begin{align*}
\dfrac{\sigma_{n+1}}{\sigma_n}
&< \dfrac{p^2}{(2n+1)(2n)}(4p^{n+1})(4p^{n+1})(p^{2n-2})( \dfrac{3}{2}p) \\
&= \dfrac{16 \cdot 3/2}{(2n+1)(2n)}p^{4n+3} \leq \dfrac{24}{20}p^{4n+3}
< \dfrac{1}{2}p^{4n+4}  \,.
\end{align*}
Finally $\dfrac{\sigma_{n+1}}{\sigma_n} < \dfrac{1}{2}p^{4n+4} <
\dfrac{\tau_{n+1}}{\tau_n}$, as required.
\Endproof\vskip4pt  

  {\it Proof of Theorem} \ref{mainthm-p}. 
Remember   the case in which $n=1$ was treated in
Proposition~\ref{pcubed}.
We have shown that $\tau_2 > \sigma_2$ and
that $\tau_{n+1}/\tau_n > \sigma_{n+1}/\sigma_n$ for all $n\geq 2$.
So, by induction, assume that $\tau_n > \sigma_n$.
We get that $\tau_{n+1} =
(\tau_{n+1}/\tau_n)\tau_n > (\sigma_{n+1}/\sigma_n)\sigma_n = \sigma_{n+1}$.
Therefore, $\tau_n > \sigma_n$ for all $n$. The theorem follows from
Proposition~\ref{strategy-p}.

The proof of Theorem~\ref{mainthm} is now complete in all cases.

\section{The detection theorem and the vanishing theorem}

Having now settled Theorem~\ref{mainthm}, we can move to
the main detection theorem (Theorem~\ref{finaldetection}) and the vanishing
theorem (Theorem~\ref{no-torsion}). Recall that they assert that if $G$ is
not cyclic, quaternion or semi-dihedral, then $T(G)$ is detected on
restriction to all elementary abelian subgroups $E$ of rank~2, and that the
torsion subgroup of~$T(G)$ is trivial.

Let us first prove a general version of the detection theorem.

\begin{thm} \label{newdetection}
For any $p$-group $G${\rm ,} the restriction homomorphism
$$
\prod_H \Res_H^G \; : T(G) \Rarr{} \prod_{H} T(H)
$$
is injective{\rm ,} where $H$ runs through the set of all subgroups of $G$ which
are elementary abelian of rank~{\rm 2,} cyclic of order~$p$ with $p$ odd{\rm ,}  cyclic
of order~{\rm 4,} and quaternion of order~$8$.
\end{thm}

\Proof 
First note that that there is nothing to prove if $G$ is cyclic of order 1
or~2, because $T(G)=\{0\}$. There is also nothing to prove if $G$ is in the
detecting family of the statement. So we can assume that $G$ is not
elementary abelian of rank~2, $C_p$, $C_4$, or~$Q_8$. By an obvious
induction argument, it suffices to prove that
$$
\prod_H \Res_H^G \; : T(G) \Rarr{} \prod_{H} T(H)
$$
is injective, where $H$ runs through the set of all maximal subgroups of
$G$.

If $G=C_{p^n}$ is cyclic (with $n\geq 2$ for $p$ odd and $n\geq 3$ if
$p=2$), then
$\Res_{C_{p^{n-1}}}^{C_{p^n}}:T(C_{p^n}) \Rarr{} T(C_{p^{n-1}})$ is an
isomorphism (both groups are isomorphic to~$\bZ/2\bZ$ generated by the class
of~$\Omega^1(k)$).
If $G$ is extraspecial or almost extraspecial, the result follows from the
main theorem of this paper (Theorem~\ref{mainthm}) if either $p=2$ or $p$ is
odd and $G$ has exponent~$p$. If $p$ is odd and $G$ has exponent~$p^2$
(extraspecial or almost extraspecial), then the result was proved in
Section~4 of~\cite{CT1}.

So we can assume that $G$ is neither cyclic, nor elementary abelian of
rank~2, nor extraspecial, nor almost extraspecial. In that case, the result
was proved as Theorem~3.2 of~\cite{CT1}.
\Endproof\vskip4pt  

This theorem provides a direct proof of the following result, which was
first proved by Puig~\cite{Pu1} using an argument of commutative algebra.

\begin{cor}
The abelian group $T(G)$ is finitely generated.
\end{cor}

As observed by Puig, this easily implies the finite generation of the Dade
group of all endo-permutation modules (see Corollary~2.4 in
Puig~\cite{Pu1}).

\Proof 
$T(H)$ is finitely generated whenever $H$ is in the detecting family. Now a
subgroup of a finitely generated group is finitely generated.
\Endproof\vskip4pt  

Theorem~\ref{newdetection} is the intermediate statement which we need
for our inductive proof of Theorem~\ref{finaldetection}. We first need to
prove the result in two special cases.

\begin{prop} \label{gps16}
Suppose that $G \cong Q_8 \times C_2$ or $G \cong D_8*C_4$.
Then $T(G)$ is detected on restriction to all elementary abelian subgroups
$E$ of rank~$2$.
\end{prop}

\Proof 
Suppose that $M$ is a nontrivial endo-trivial module such that
$M{\downarrow}_{E}^G \cong k\oplus \text{(free)}$ for every elementary
abelian subgroup $E$ of rank~2. Assume that $M$
has minimal dimension among such modules. On restriction to a maximal
subgroup of the form $C_4{\times}C_2$, we must have that
$M{\downarrow}_{C_4{\times}C_2}^G \cong k\oplus \text{(free)}$, because
$T(C_4{\times}C_2) \Rarr{} T(E)$ is an isomorphism for $E=C_2{\times}C_2
\subset C_4{\times}C_2$. It follows that $\Dim(M)\equiv {1\pmod{8}}$. It
also follows that $M{\downarrow}_{C_4}^G \cong k\oplus \text{(free)}$ for
any cyclic subgroup $C_4$, because $C_4$ is contained in a maximal
subgroup of the form $C_4{\times}C_2$.

Since $M$ is nontrivial, it must be detected on some restriction
(Theorem~\ref{newdetection}). So there exists a quaternion
subgroup $H \cong Q_8$ in $G$ such that
$M{\downarrow}_H^G$ is nontrivial. Then
$M{\downarrow}_H^G \cong \Omega^2(k_H) \oplus \text{(free)}$, because
$\Omega^2(k_H)$ is the only indecomposable endo-trivial
$kH$-module other than $k_H$ itself whose dimension is congruent to~1
modulo~8 (see \cite[\S 6]{CT1}).

Let $z$ be the generator of the center of~$H$ (which is also central
in~$G$). We consider the variety $V_{\overline{G}}(\overline{M})
\subseteq V_{\overline{G}}(k) \cong k^3$
where, as in Section~5, $\overline{G} = G/\langle z \rangle$ and
$\overline{M} \cong (z-1)M$. On restriction to $H$, we have
$V_{\overline{H}}(\overline{M}) =
V_{\overline{H}}(\,\overline{\Omega^2(k_H)}\,)$ and
$\overline{\Omega^2(k_H)}$ is a periodic $k\overline{H}$-module by
Lemma~\ref{barperiod}. Since $\Omega^2(k_H)$ is invariant
under Galois automorphisms, so is $\overline{\Omega^2(k_H)}$, and
therefore $V_{\overline{H}}(\overline{M})$ is a union of lines permuted
by Galois automorphisms. But these lines are not ${\Bbb F}_2$-rational (by
Lemma~\ref{fprationality} applied to the $kH$-module $\Omega^2(k_H)$,
which is critical), hence not fixed by Galois automorphisms. It follows
that there are at least two lines in
$V_{\overline{H}}(\overline{M})$ (and in fact exactly two, which are
${\Bbb F}_4$-rational, because this is the only possibility for the
4-dimensional module $\overline{\Omega^2(k_H)}\;$). Now
$V_{\overline{G}}(\overline{M})$ also contains at least two lines since
it contains 
$\text{res}_{\overline{G},\overline{H}}^*(V_{\overline{H}}(\overline{M}))$.
So $\overline{M} \cong \overline{M_1} \oplus \overline{M_2}$  where
$V_{\overline{G}}(\overline{M_1})$ is one of the two lines. Now following
the procedure of Theorem~\ref{separate-var} we can construct a
nontrivial endo-trivial $kG$-module $N_1$ such that
$\overline{N_1} = \overline{M_1}$. Moreover $N_1$ is trivial on
restriction to every elementary abelian subgroup. But $\Dim(N_1)<\Dim(M)$,
contrary to the choice of~$M$.
\Endproof\vskip4pt  

We also need a group-theoretical lemma.

\begin{lemma} \label{subgroups}
Let $G$ be a semi-direct product $G = Q_{2^n}\rtimes C_2$ for some
$n\geq 3$ and some action of $C_2$ on $Q_{2^n}$. Then one of the following
properties holds\/{\rm :}\/
\begin{itemize}
\ritem{(a)} $G$ contains a semi-dihedral subgroup $S$ such that
$S \supseteq Q_8 \subseteq Q_{2^n}$.

\ritem{(b)} $G$ contains a subgroup $Q_8*C_4$ with $Q_8\subseteq Q_{2^n}$.

\ritem{(c)} $G$ contains a subgroup $Q_8\times C_2$ with $Q_8\subseteq Q_{2^n}$.
\end{itemize}
\end{lemma}

\Proof 
Let $u$ be a generator of~$C_2$. We use induction on $n$ and first consider
the case $n=3$. If the action of $u$ on $Q_8/Z(Q_8)$ is nontrivial, then
we can choose two generators $x$ and $y$ of~$Q_8$ such that $uxu^{-1} =
y$. In that case $G$ is semi-dihedral and we are in case~(a). If now $u$
acts trivially on $Q_8/Z(Q_8)$, then $u$ fixes each of the three cyclic
subgroups of order~4 of~$Q_8$. If $u$ acts trivially on~$Q_8$, then $G =
Q_8\times C_2$ and we are in case~(c). Otherwise it easy to see that $u$
must invert two of the cyclic subgroups of order~4 and fix pointwise the
third one, say~$\langle x \rangle$. But then the actions of $u$ and~$x$
coincide, so that  $ux^{-1}$ acts trivially and $G = Q_8*C_4$, which is case~(b).

Assume now that $n\geq 4$. Let $x$ and $y$ be generators of $Q_{2^n}$ with
$x^{2^{n-1}}=1$, $y^2=x^{2^{n-2}}$ and $yxy^{-1}=x^{-1}$. All elements of
the form $x^by$ have order~4 (where $0 \leq b \leq 2^{n-1}$). Conjugation
by~$u$ must satisfy  $uxu^{-1}=x^a$ for some odd integer~$a$ and $uyu^{-1}
= x^by$ for some~$b$. Since $u^2=1$, we must have the following congruences
modulo $2^{n-1}\,$:
$$ a \equiv \pm 1 \,, 2^{n-2} \pm 1
\qquad \text{and} \qquad
(a+1)b \equiv 0 \,.$$
If $a \equiv -1$ and $b$ is odd, we can replace $x$ by $x^b$ and we
get a standard presentation of the semi-dihedral group ${\rm SD}_{2^{n+1}}$, so
we are in case~(a). Otherwise $b$ must be even, because this is forced by
the condition $(a+1)b \equiv 0$ if $a \not\equiv -1$. Therefore
conjugation by
$u$ stabilizes the subgroup $Q_{2^{n-1}}$ generated by $x^2$ and~$y$.
The result now follows by induction applied to the group
$Q_{2^{n-1}}\rtimes \langle u \rangle$.
\Endproof\vskip4pt

Now we come to the detection theorem (Theorem~\ref{finaldetection} of the
introduction).

\begin{thm} 
Suppose that $G$ is a $p$-group which is not cyclic{\rm ,} quaternion{\rm ,} or
semi-dihedral. Then $T(G)$ is detected on restriction to all elementary
abelian subgroups $E$ of rank~$2$.
\end{thm}

\Proof 
We use induction on the order of~$G$. First recall that the result is
known if $G$ is abelian or dihedral (see~\cite{CT1}); so we assume
that $G$ is neither abelian nor dihedral.

Let $M$ be an endo-trivial module such that
$\Res_E^G[M]=0$ for every elementary abelian subgroup $E$ of rank~2, where
$[M]$ denotes the class of~$M$ in~$T(G)$. It suffices to prove that
$\Res_H^G[M]=0$ for every maximal subgroup $H$ of~$G$, because then $[M]=0$
by Theorem~\ref{newdetection}. For every maximal subgroup~$H$ which is not
cyclic, quaternion, or semi-dihedral, $M{\downarrow}_H^G$ satisfies the same
assumption as~$M$, so that $\Res_H^G[M]=0$ by induction. Now, we are left with the
cases where the maximal subgroup~$H$ is cyclic, quaternion, or
semi-dihedral.

Assume first that $H \cong C_{p^n}$ is cyclic. By a well-known result of
group theory (see Theorem~4.4 in Chapter~5 of~\cite{Gor1}), $G$ is either
abelian, or isomorphic to a group~$P$ to be described below, or in addition
when $p=2$, isomorphic to $D_{2^{n+1}}$, $Q_{2^{n+1}}$, or
${\rm SD}_{2^{n+1}}$. The cases of the cyclic group $C_{p^{n+1}}$, the quaternion
group $Q_{2^{n+1}}$, or the
semi-dihedral group ${\rm SD}_{2^{n+1}}$, are excluded by our hypothesis. The
cases of an abelian group or a dihedral group $D_{2^{n+1}}$ have  
already been dealt with. So we are left with the case
$G = P = H \rtimes C_p$, with respect to the action $uxu^{-1} =
x^{1+p^{n-1}}$, where $x$ is a generator of~$H$ and $u$ is a generator
of~$C_p$. This case occurs if $n\geq 2$ when $p$ is odd and $n \geq 3$ when
$p=2$. Now $G$ also contains a maximal
subgroup $K = \langle x^p \rangle \times \langle u \rangle \cong
C_{p^{n-1}}\times C_p$ and we already know that $\Res_K^G[M]=0$.
Therefore
$$\Res_{C_{p^{n-1}}}^G[M] =
\Res_{C_{p^{n-1}}}^K \; \Res_K^G[M] =0 \,.$$
But we also have $\Res_{C_{p^{n-1}}}^G =
\Res_{C_{p^{n-1}}}^H \; \Res_H^G$ and
$$ \Res_{C_{p^{n-1}}}^H : T(H) \Rarr{} T(C_{p^{n-1}})$$
is an isomorphism since both $T(H)$ and $T(C_{p^{n-1}})$ are cyclic
of order~2 generated by the class of~$\Omega^1(k)$ (because $n\geq 2$ and
$n\geq 3$ if $p=2$). It follows that
$\Res_H^G[M] = 0$.

Assume now that $H \cong {\rm SD}_{2^n}$ is semi-dihedral. We know that the
torsion subgroup $T_t(H)$ is cyclic of order~2 generated by the class of an
endo-trivial module whose dimension is congruent to 1 modulo~$2^{n-1}$
(see~\cite[\S 7]{CT1}). This class cannot be in the image of
$\Res_H^G$, because all endo-trivial modules for~$G$ have dimension
congruent to $\pm 1$ modulo~$2^n$, by Lemma~2.10 in~\cite{CT1}. It
follows that the image of $\Res_H^G$ is
contained in $\langle \,[\Omega_H^1(k)]\, \rangle \cong \bZ$, because $T(H)
= T_t(H) \oplus \langle \,[\Omega_H^1(k)]\, \rangle$. But now the
restriction map
 $$\Res_{E}^H \,: \langle \,[\Omega_H^1(k)]\, \rangle \Rarr{}
T(E) = \langle \,[\Omega_{E}^1(k)]\, \rangle $$
is an isomorphism where $E$ is an elementary abelian subgroup of rank~2.
Since $\Res_{E}^H \, \Res_H^G[M]=0$, we must have
$\Res_H^G[M] = 0$ as required.
Note that the same argument shows that $\Res_S^G[M] = 0$ for any
semi-dihedral subgroup~$S$ of~$G$.

Assume finally that $H \cong Q_{2^n}$ is quaternion. Since $G$ is neither
cyclic nor quaternion, its 2-rank cannot be~1 (see Chapter~5
of~\cite{Gor1}) and so there exists an element of order~2 outside~$H$.
Therefore $G \cong Q_{2^n}\rtimes C_2$ for some action
of $C_2$ on~$Q_{2^n}$. By Lemma~\ref{subgroups}, $G$ contains a
subgroup~$R$ which is isomorphic to
$Q_8*C_4$, $Q_8\times C_2$, or semi-dihedral, and such that $R\supseteq
Q_8 \subseteq H$. In the first two cases we have $\Res_R^G[M] = 0$ by
Proposition~\ref{gps16} and in the third we have $\Res_R^G[M] = 0$ by the
argument above. It follows that
$$ \Res_{Q_8}^H \; \Res_H^G[M]  =\Res_{Q_8}^G[M]  = \Res_{Q_8}^R
\;\Res_R^G[M] =0 \,. $$
We know that $T(H) \cong \bZ/4\bZ \oplus \bZ/2\bZ$, where
$\bZ/4\bZ$ is generated by the class of~$\Omega_H^1(k)$ and $\bZ/2\bZ$ is
generated by the class of an endo-trivial module of
dimension~$2^{n-1}+1$ (see~\cite[\S 6]{CT1}). Again this class cannot
be in the image of~$\Res_H^G$, because all endo-trivial modules for~$G$ have
dimension congruent to $\pm 1$ modulo~$2^n$. Thus the image of $\Res_H^G$ is
contained in $\bZ/4\bZ = \langle \,[\Omega_H^1(k)]\, \rangle$. But now the
restriction map
 $$\Res_{Q_8}^H \,: \langle \,[\Omega_H^1(k)]\, \rangle \Rarr{}
\langle \,[\Omega_{Q_8}^1(k)]\, \rangle $$
is an isomorphism. Since $\Res_{Q_8}^H \, \Res_H^G[M] =0$, we must have
$\Res_H^G[M] = 0$ as required.
\Endproof\vskip4pt  

We immediately deduce the vanishing theorem (Theorem~\ref{no-torsion} of the
introduction).

\begin{cor}
If $G$ is not cyclic{\rm ,} quaternion or semi-dihedral{\rm ,} then the torsion subgroup
of~$T(G)$ is trivial.
\end{cor}

\Proof 
By the theorem, we know that $T(G)$ is embedded in a product of copies of
${T(E)\cong\bZ}$, where $E$ is elementary abelian of rank~2.
\Endproof\vskip4pt  

We can now prove Corollary \ref{rank-one} of the introduction.

\begin{cor}
Suppose that $G$ is a finite $p$-group for which every maximal elementary
subgroup has rank at least~$3$. Then $T(G) \cong \bZ$, generated by the class
of the module $\Omega^1(k)$.
\end{cor}

\Proof 
The assumption implies that $G$ cannot be cyclic, quaternion or
semi-dihedral. Therefore, by the theorem, $T(G)$ is detected on restriction
to elementary abelian subgroups of rank~2. The rest of the proof follows
Alperin~\cite{Alpsyz} and we recall the argument (also used in~\cite{Both}).
The partially ordered set of all elementary abelian subgroups of rank at
least~2 is connected, in view of the assumption and by a well-known
result of the theory of\break $p$-groups. For any such subgroup~$H$, the
restriction map $T(H) \Rarr{} T(E)\cong \bZ$ to an elementary abelian
subgroup of rank~2 is an isomorphism. It follows that all restrictions to
such rank~2 subgroups~$E$ are equal.
\hfill\qed

\section{The Dade group}

In this section, we prove detection theorems for the Dade group $D(G)$
of all endo-permutation modules and we determine its torsion subgroup when
$p$ is odd. We refer to~\cite{Both} for details
about~$D(G)$. Let us only mention that the torsion-free rank of~$D(G)$ has
been determined in~\cite{Both} so that we are particularly interested in
the torsion subgroup~$D_t(G)$. We first state an easy consequence of
Theorem~\ref{newdetection}.

\begin{thm} \label{easydetection}
Let $G$ be a finite $p$-group.
\begin{itemize}
\ritem{(a)} The product of all restriction-deflation maps
 $$\prod_{K/H} \Def_{K/H}^K \; \Res_K^G \; : \; D(G) \Rarr{}
\prod_{K/H} D(K/H)$$
is injective{\rm ,} where $K/H$ runs through the set of all sections of~$G$
which are elementary abelian of rank~{\rm 2,} cyclic of order~$p$ with $p$ odd{\rm ,}
cyclic of order~{\rm 4,} or quaternion of order~{\rm 8}.

\ritem{(b)} For the torsion subgroup{\rm ,} the product of all restriction-deflation maps
 $$\prod_{K/H} \Def_{K/H}^K \; \Res_K^G \; : \; D_t(G) \Rarr{}
\prod_{K/H} D_t(K/H)$$
is injective{\rm ,} where $K/H$ runs through the set of all sections of~$G$
which are cyclic of order~$p$ if $p$ is odd{\rm ,} quaternion of order~{\rm 8} or cyclic
of order~{\rm 4} if $p=2$.\end{itemize}
\end{thm}

\Proof 
The argument is exactly the same as the one given in Theorem~1.6
of~\cite{Both} or in Theorem~10.1 of~\cite{CT1}.
\Endproof\vskip4pt

We now deduce Corollary~\ref{self-dual} of the introduction.

\begin{cor}
Let $G$ be a finite $p$-group.
\begin{itemize}
\ritem{(a)} If $p$ is odd\/{\rm ,}\/ any nontrivial torsion element in $D(G)$ has order {\rm 2}. In
other words{\rm ,} for any indecomposable endo-permutation $kG$-module~$M$ with
vertex~$G${\rm ,} the class of~$M$ is a torsion element if and only if $M$ is
self-dual.

\ritem{(b)} If $p=2$\/{\rm ,}\/ any nontrivial torsion element in $D(G)$ has order \/{\rm 2}\/ or~\/{\rm 4}\/.
\end{itemize}
\end{cor}

\Proof 
The nontrivial elements of $D(C_p)$ have order~2, while those of $D(Q_8)$
and $D(C_4)$ have order 2 or~4. Moreover, an element of order~2 corresponds
to a self-dual module by definition of the group law.
\Endproof\vskip4pt  

If $p$ is odd, the detection theorem above allows for a complete description
of the torsion subgroup of~$D(G)$ (Theorem~\ref{odd-Dade-group} of the
introduction), by  the partial results already obtained in~\cite{Both}.

\begin{thm}
If $p$ is odd and $G$ is a finite $p$-group{\rm ,} the torsion subgroup of $D(G)$
is isomorphic to $(\bZ/2\bZ)^s${\rm ,} where $s$ is the number of conjugacy
classes of nontrivial cyclic subgroups of~$G$.
\end{thm}

Note that explicit generators are described in~\cite{Both}.

\Proof 
Theorem~6.2 in~\cite{Both} asserts that a certain quotient
$\overline{D_t(G)}$ of the torsion subgroup~$D_t(G)$ is isomorphic to
$(\bZ/2\bZ)^s$, where $s$ is as above. So we only have to prove that
$\overline{D_t(G)} = D_t(G)$. But by definition, $\overline{D_t(G)} =
D_t(G)/\Ker(\psi)$, where
$\psi$ is the product of all restriction-deflation maps
 $$\psi = \prod_{K/H} \Def_{K/H}^K \; \Res_K^G \; : \; D_t(G)\Rarr{} \prod_{K/H} D_t(K/H)$$
where $K/H$ runs through the set of all sections of~$G$
which are cyclic of order~$p$. Now $\psi$ is injective by
Theorem~\ref{easydetection} and the result follows.
\Endproof\vskip4pt  

Our purpose now is to improve Theorem~\ref{easydetection} by restricting the
kind of section needed on the right-hand side. However, we will also
change the target by including all groups having torsion endo-trivial
modules, namely cyclic, quaternion, and semi-dihedral groups.

If $S=\langle x \rangle$ is cyclic of order~$p^n$, then
$$ 
D(S) = D_t(S) \cong \prod_{i=1}^n T_t(S/\langle x^{p^i} \rangle)\,,
$$
and we let $\pi_S : D_t(S) \to T_t(S)$ denote the projection onto the
factor indexed by $i=n$. The situation is easier if $S$ is a quaternion or
semi-dihedral group, since $D_t(S) = T_t(S)$ by~\cite[\S 10]{CT1}. In
this case, we write $\pi_S : D_t(S) \to T_t(S)$ for the identity map.

\begin{thm} \label{Dadedetection}
Let $G$ be a finite $p$-group.
If $p$ is odd{\rm ,} let $\CX$ be the class of all subgroups
$H$ of~$G$ such that $N_G(H)/H$ is cyclic.
If $p=2${\rm ,} let $\CX$ be the class of all subgroups
$H$ of~$G$ such that $N_G(H)/H$ is cyclic of order~$\geq 4${\rm ,} quaternion of
order~$\geq 8${\rm ,} or semi-dihedral of order~$\geq 16$.
Let $[\CX/G]$ be \pagebreak a system of representatives of conjugacy classes of
subgroups in~$\CX$. Then the map
\begin{eqnarray*}
&&\hskip-35pt\prod_{H\in[\CX/G]} \pi_{N_G(H)/H} \; \Def_{N_G(H)/H}^{N_G(H)} \;
\Res_{N_G(H)}^G : \\
&&\qquad D_t(G)  \Rarr{}  
\prod_{H\in[\CX/G]} T_t(N_G(H)/H)
\end{eqnarray*}
 is injective.
\end{thm}

\Proof 
Let $\varphi$ denote the map in the statement and let
$a \in \Ker(\varphi)$, so that
$\pi_{N_G(H)/H}\,\Defres_{N_G(H)/H}^G(a)=0$ for every $H\in\CX$, where we
write for simplicity
$\Defres_{K/H}^G=\Def_{K/H}^K\,\Res_K^G$ for every section~$K/H$. By
Theorem~\ref{easydetection} above, it suffices to prove that
$\Defres_{K/H}^G(a)=0$ for every section~$K/H$ isomorphic to $C_p$, $C_4$
or~$Q_8$. We are going to show that
$\Defres_{N_G(H)/H}^G(a)=0$ and the result will follow from this since
$\Defres_{K/H}^G=\Res_{K/H}^{N_G(H)/H}\, \Defres_{N_G(H)/H}^G$. For
simplicity of notation, we write now $L=N_G(H)$.

We use induction on the index~$|G:H|$. If $H$ has index~$p$, there is
nothing to prove because $L=G$, $ \pi_{G/H}={\rm id}$, and
 $$\Defres_{G/H}^G(a) = \pi_{G/H} \,\Defres_{G/H}^G(a)=0 \,,$$
by assumption if $p$ is odd and by the fact that $D(G/H)=\{0\}$ if $p=2$.
Let $F$ be a subgroup such that $H<F\leq L$.
By induction, $\Defres_{N_G(F)/F}^G(a)=0$ and consequently
$\Defres_{N_L(F)/F}^G(a)=0$. This holds for every such $F$ and therefore
$$
\Defres_{L/H}^G(a) \in \bigcap_{H<F\leq L}
\Ker(\Defres_{N_L(F)/F}^{L/H})\, = \, T(L/H)\,.
$$
The last equality is a well-known characterization of $T(L/H)$ as a
subgroup of~$D(L/H)$ (see Lemma~2.1 in~\cite{CT1} and note that this
characterization is also at the heart of the proof of
Theorem~\ref{easydetection}). Since $a$ was chosen to be a torsion element
in~$D(G)$, we have proved that $\Defres_{L/H}^G(a) \in T_t(L/H)$.

If $L/H$ is not cyclic, quaternion, or semi-dihedral, then $T_t(L/H) =
\{0\}$ by Theorem~\ref{no-torsion} and so $\Defres_{L/H}^G(a)=0$. The
same holds if $L/H$ is cyclic of order~2. If $L/H$ is quaternion or
semi-dihedral, then $\pi_{L/H}$ is the identity map and
$\pi_{L/H}\,\Defres_{L/H}^G(a)=0$ by assumption, so
$\Defres_{L/H}^G(a)=0$. If $L/H$ is cyclic of order $\geq 3$, then
$\pi_{L/H}:D_t(L/H)\to T_t(L/H)$ restricts to the identity on~$T_t(L/H)$.
Since $\pi_{L/H}\,\Defres_{L/H}^G(a)=0$ by assumption, we obtain again
${\Defres_{L/H}^G(a)=0}$.
\Endproof\vskip4pt  

In order to illustrate the efficiency of Theorem~\ref{Dadedetection}
compared to Theorem~\ref{easydetection}, suppose that $G$ is abelian.
Then there are numerous sections of~$G$ isomorphic to $C_p$ or~$C_4$ and the
map in Theorem~\ref{easydetection} is an injection in a much larger group,
whereas the map in Theorem~\ref{Dadedetection} hits exactly every cyclic
quotient of~$G$ and is an isomorphism (Dade's theorem).

If $G$ is a dihedral 2-group, there are many sections of~$G$ isomorphic
to~$C_4$, but $N_G(H)/H$ is never cyclic of order~$\geq 4$, quaternion, or
semi-dihedral, so that $\CX$ is empty and $D_t(G)=\{0\}$, a result also
obtained in~\cite[\S 10]{CT1}.

Theorem~\ref{Dadedetection} allows us to handle also a case where the
structure of $D_t(G)$ was not previously known.

\begin{prop}
Suppose that $G$ is an extraspecial $2$-group of type~{\rm 1,} that is{\rm ,} a central
product of copies of~$D_8$. Then $D_t(G)=\{0\}$.
\end{prop}

\Proof 
We claim that $\CX$ is empty and so $D_t(G)=\{0\}$. If $H$ is a subgroup
of~$G$ containing~$Z(G)$, then $H$ is a normal subgroup, $G/H$ is
elementary abelian, and $H\notin\CX$. If $H$ does not contain~$Z(G)$,
then for any $g\in N_G(H)$, we have that $[g,h]\in H\cap Z(G) =\{1\}$. Thus
$N_G(H) = C_G(H)$ and in particular $H$ is abelian, actually elementary
abelian, since the square of every element of $H$ belongs to ${H\cap Z(G)}
=\{1\}$. Using the quadratic form on $G/Z(G)$, it is not hard to prove
that if $n$ is the number of copies of~$D_8$ in the central product and if
$H=(C_2)^k$, then $C_G(H) = H\times L$ where $L$ is a central product of
$n{-}k$ copies of~$D_8$ (possibly $n{-}k=0$ and $L=Z(G)$). Therefore
$N_G(H)/H$ is extraspecial and $H\notin\CX$. This proves that $\CX$ is
empty.
\hfill\qed

\section{Two examples}

Theorem~\ref{Dadedetection} is not sufficient to determine $D_t(G)$ in
all cases when $p=2$. This seems to be in contrast to the case of an odd
prime, for which the solution of the detection conjecture for~$T(G)$ allows
for a complete description of~$D_t(G)$ (Theorem~\ref{odd-Dade-group}).

Our purpose is to illustrate the situation with the extraspecial groups of
type~2 and the almost extraspecial groups (type~3). For simplicity, we
shall only deal with the smallest of the groups, namely $D_8*Q_8$ and
$D_8*C_4$, but our results can easily be generalized to the other groups
of types 2 and~3.

If follows from Theorem~\ref{Dadedetection} that the product of all
restriction-deflation maps
 $$\prod_{H\in[\CX/G]} \Def_{N_G(H)/H}^{N_G(H)} \;
\Res_{N_G(H)}^G : D_t(G) \Rarr{}
\prod_{H\in[\CX/G]} D_t(N_G(H)/H)$$
 is injective. In the opposite direction, there is the sum of all maps
obtained by composing inflation maps $\Inf_{N_G(H)/H}^{N_G(H)}$ and
tensor induction $\Ten_{N_G(H)}^G$, namely
 $$
\sum_{H\in[\CX/G]} \Ten_{N_G(H)}^G \; \Inf_{N_G(H)/H}^{N_G(H)} :
\bigoplus_{H\in[\CX/G]} D_t(N_G(H)/H) \Rarr{} D_t(G) \,.
$$
We let $D_t^0(G)$ be the image of this map. The question of the
surjectivity of this map does not seem to be easy and this is why we have
to introduce the subgroup~$D_t^0(G)$. In similar situations for odd
primes, or for the Dade group tensored with~$\bQ$, we can prove the
surjectivity of the map (see Sections 4 and~6 of~\cite{Both}), so it seems
natural to conjecture that $D_t^0(G)=D_t(G)$. In our two examples, we
shall be able to compute $D_t^0(G)$ but it is not easy to know if $D_t(G)$
is larger or not.

In order to compute the image by restriction-deflation of elements of
$D_t^0(G)$, we need a technical formula which is derived from the
results of~\cite{Both}. There is a general formula describing the
restriction-deflation of an element of the form $\Ten_K^G \, \Inf_{K/H}^K
(x)$, but for simplicity we only consider two very special cases. The
Frobenius map $\lambda\mapsto \lambda^{p^n}$ is an endomorphism of~$k$ and
we let 
$$\gamma_{p^n}:D(G) \Rarr{} D(G)$$
be the group homomorphism induced by the Frobenius map, as defined in
Section~3 of~\cite{Both}.

\begin{lemma} \label{formula}
Let $G$ be a $p$-group and let $K$ and $H$ be subgroups of $G$
such that $H$ is a normal subgroup of~$K$.
\smallbreak
{\rm (a)} Let $P$ and $R$ be subgroups of $G$ such that $R$ is a normal subgroup
of~$P$.
Assume that $K$ and $P$ satisfy $KP=G$ \/{\rm (}\/a single double coset\/{\rm )}\/.
Assume further that the inclusions $P\cap K \to K$ and $P\cap K \to P$
induce isomorphisms
$$(P\cap K)/(R\cap H) \Rarr{\sim} K/H \qquad {\text {and}} \qquad
(P\cap K)/(R\cap H) \Rarr{\sim} P/R$$
 respectively. Then the following maps
from $D(K/H)$ to $D(P/R)$ are equal\/{\rm :}\/
$$ \Def_{P/R}^P \; \Res_P^G \;\Ten_K^G \; \Inf_{K/H}^K
= \gamma_{|R:R\cap H|} \; \Iso_{(P\cap K)/(R\cap H)}^{P/R} \;
(\Iso_{(P\cap K)/(R\cap H)}^{K/H})^{-1} \;,$$
 where the two latter maps are induced by the isomorphisms
$(P\cap K)/(R\cap H) \Rarr{\sim} P/R$ and
$(P\cap K)/(R\cap H) \Rarr{\sim} K/H$ respectively.

\smallbreak
{\rm (b)} Let $L$ be a normal subgroup of $K$. Then the following maps
from $D(K/H)$ to $D(K/L)$ are equal\/{\rm :}\/
$$\Def_{K/L}^K \; \Inf_{K/H}^K = \Inf_{K/HL}^{K/L} \; \Def_{K/HL}^{K/H}
\;.$$
\end{lemma}
\vglue8pt

\Proof 
(a) Since there is a single double coset, the Mackey formula implies that
$$ \Def_{P/R}^P \; \Res_P^G \;\Ten_K^G \; \Inf_{K/H}^K
= \Def_{P/R}^P \; \Ten_{P\cap K}^P \; \Res_{P\cap K}^K \; \Inf_{K/H}^K
\,.$$
Now Proposition 3.10 in~\cite{Both} asserts that
$$\Def_{P/R}^P \; \Ten_Q^P = \gamma_{|R:Q\cap R|} \;
\Ten_{QR/R}^{P/R} \; \Iso_{Q/Q\cap R}^{QR/R} \; \Def_{Q/Q\cap R}^Q
\;.$$
Applying this with $Q=P\cap K$, we have that $QR=P$ and $Q\cap R= R\cap H$,
because of the assumed isomorphism $(P\cap K)/(R\cap H) \Rarr{\sim} P/R$,
and therefore
$$ \Def_{P/R}^P \; \Ten_{P\cap K}^P = \gamma_{|R:R\cap H|} \;
 \Iso_{P\cap K/R\cap H}^{P/R} \; \Def_{P\cap K/R\cap H}^{P\cap K}
\;.$$
Composing on the right with $\Res_{P\cap K}^K \; \Inf_{K/H}^K$, it is
easy to see that
 $$\Def_{P\cap K/R\cap H}^{P\cap K} \; \Res_{P\cap K}^K \; \Inf_{K/H}^K
= (\Iso_{P\cap K/R\cap H}^{K/H})^{-1}\;,$$
using either the definitions of the maps or the methods of Corollary~3.9
in~\cite{Both}. It follows that
$$\Def_{P/R}^P \; \Ten_{P\cap K}^P \; \Res_{P\cap K}^K \; \Inf_{K/H}^K
= \gamma_{|R:R\cap H|} \; \Iso_{P\cap K/R\cap H}^{P/R} \;
(\Iso_{P\cap K/R\cap H}^{K/H})^{-1}\;,$$
and the result follows.
\smallbreak

(b) This follows either from the definitions of the maps or from the
methods of Corollary~3.9 in~\cite{Both}.
\Endproof\vskip4pt  

Now we can start with our first example $D_8*C_4$. Let $S_1,S_2,S_3$ be
representatives of the three  conjugacy classes of noncentral subgroups of
order~2 (the two classes in~$D_8$ and the product of a generator of~$C_4$
with an element of order~4 in~$D_8$).

\begin{prop} \label{aes16}
Let $G=D_8*C_4$ be the almost extraspecial group of order $16$. Then
$D_t^0(G)$ is cyclic of order~{\rm 2,} generated by the class of the module
$\Ten_{S_1\times C_4}^G\,\Inf_{S_1\times C_4/S_1}^{S_1\times C_4}
(\Omega_{S_1\times C_4/S_1}^1(k))$.
\end{prop}

\Proof 
We have that $N_G(S_i)=S_i\times C_4$ and so
$N_G(S_i)/S_i\cong C_4$ and $S_i$ is in the class~$\CX$ of
Theorem~\ref{Dadedetection}. These are the only subgroups in~$\CX$
(because every other nontrivial subgroup
$H$ contains the Frattini subgroup and $G/H$ is elementary abelian).
Therefore Theorem~\ref{Dadedetection} yields an injective map
$$\prod_{i=1}^3\; \Def_{S_i\times C_4/S_i}^{S_i\times C_4}\;
\Res_{S_i\times C_4}^G : D_t(G) \Rarr{}
\prod_{i=1}^3 D_t(S_i\times C_4/S_i) \cong (\bZ/2\bZ)^3\;,$$
each factor $D_t(S_i\times C_4/S_i)\cong D_t(C_4)$ being cyclic of
order~2 generated by the class of~$\Omega_{S_i\times C_4/S_i}^1(k)$.
Now by definition $D_t^0(G)$ is generated by the three elements
$$\Ten_{S_i\times C_4}^G\,\Inf_{S_i\times C_4/S_i}^{S_i\times C_4}
(\Omega_{S_i\times C_4/S_i}^1(k))\quad (1\leq i\leq3)\,.$$
We claim that they are all equal and have order~2. This will
complete the proof of the proposition.

In order to prove the claim, we show that the image of
any of these three elements by the injective map above is equal to
the ``diagonal element''
$$(\Omega_{S_1\times C_4/S_1}^1(k)\,,\,\Omega_{S_2\times C_4/S_2}^1(k)\,,\,
\Omega_{S_3\times C_4/S_3}^1(k)\,)\,.$$
This follows from a straightforward application of
Lemma~\ref{formula}. If $i\neq j$, we obtain
\begin{eqnarray*}
&&\hskip-.25in\Def_{S_j\times C_4/S_j}^{S_j\times C_4}\;
\Res_{S_j\times C_4}^G\; \Ten_{S_i\times C_4}^G\;
\Inf_{S_i\times C_4/S_i}^{S_i\times C_4}
(\Omega_{S_i\times C_4/S_i}^1(k)) \\
&&\qquad=  \gamma_{|S_j:1|} \; \Iso_{C_4/1}^{S_j\times C_4/S_j} \;
(\Iso_{C_4/1}^{S_i\times C_4/S_i})^{-1}
(\Omega_{S_i\times C_4/S_i}^1(k)) \\
&&\qquad=  \gamma_2(\Omega_{S_j\times C_4/S_j}^1(k))
=\Omega_{S_j\times C_4/S_j}^1(k) \;,
\end{eqnarray*}
using the fact that $\Omega^1(k)$ is defined over the prime
field~${\Bbb F}_2$ and hence is fixed by the Frobenius map~$\gamma_2$.
In the case $i=j$, we have for any $k[S_i\times C_4]$-module~$M$,
$$\Res_{S_i\times C_4}^G\, \Ten_{S_i\times C_4}^G(M)
=M\otimes{\,^gM}\,,$$
where $g$ is a representative of the nontrivial class of $G/S_i\times
C_4$ and ${\,^gM}$ denotes the conjugate module. Therefore, ignoring
inflation for simplicity, we obtain
\begin{eqnarray*}
&& 
\Def_{S_i\times C_4/S_i}^{S_i\times C_4}\;
\Res_{S_i\times C_4}^G\; \Ten_{S_i\times C_4}^G\;
\Inf_{S_i\times C_4/S_i}^{S_i\times C_4}
(\Omega_{S_i\times C_4/S_i}^1(k)) \\
&&\qquad= \Def_{S_i\times C_4/S_i}^{S_i\times C_4}\,
\big(\Omega_{S_i\times C_4/S_i}^1(k) \otimes
{\,^g(\Omega_{S_i\times C_4/S_i}^1(k))} \big) \\
&&\qquad= \Def_{S_i\times C_4/S_i}^{S_i\times C_4}\,
\big(\Omega_{S_i\times C_4/S_i}^1(k) \big) \otimes
\Def_{S_i\times C_4/S_i}^{S_i\times C_4}\,
\big(\Omega_{S_i\times C_4/gS_ig^{-1}}^1(k) \big) \\
&&\qquad= \Omega_{S_i\times C_4/S_i}^1(k) \otimes
\Def_{S_i\times C_4/S_i}^{S_i\times C_4}
(\Omega_{S_i\times C_4/gS_ig^{-1}}^1(k)) \,.
\end{eqnarray*}
But the second factor is trivial because, by part (b) of
Lemma~\ref{formula} with $K=S_i\times C_4$, we have
$$\Def_{K/S_i}^K\,\Inf_{K/gS_ig^{-1}}^K
= \Inf_{K/S_i(gS_ig^{-1})}^{K/S_i} \;
\Def_{K/S_i(gS_ig^{-1})}^{K/gS_ig^{-1}}
$$
and a deflation of the class of $\Omega^1(k)$ is trivial (see
Lemma~1.3 of~\cite{Both}).
\Endproof\vskip4pt  

In this example, we see that $D_t(G)$ embeds in three copies
of~$\bZ/2\bZ$ and that $D_t^0(G)\cong \bZ/2\bZ$. So in order to prove the
conjectural equality $D_t^0(G)=D_t(G)$, we would have to improve
Theorem~\ref{Dadedetection} by showing the injectivity of the
restriction-deflation map to a single section $S_i\times C_4/S_i$.
In this specific example, we have been able to do this by a rather
delicate argument  not given here.

The methods are similar with our second example $D_8*Q_8$, but another
complication occurs. Recall that $D_t(Q_8)$ is generated by the class
of~$\Omega_{Q_8}^1(k)$, which has order~4, and the class of a certain
5-dimensional module~$M$, which has order~2 (see~\cite[\S 6]{CT1}).
Moreover $M$ is defined over the field~${\Bbb F}_4$ (so we assume here
that $k$ contains~${\Bbb F}_4$) and $M$ is not invariant under the Galois
automorphism~$\gamma_2$. Actually $\gamma_2(M)\cong \Omega^2(M)$, another
5-dimensional module, and $\Omega^2(k)$, $M$, $\Omega^2(M)$ are the three 
elements of order~2 in $D_t(Q_8)\cong \bZ/4\bZ \oplus \bZ/2\bZ$.

Let $S_1,\dots,S_5$ be representatives of the five conjugacy classes of
noncentral subgroups of order~2 (the two classes in~$D_8$ and the product
of an element of order~4 in~$D_8$ with one of the three possible elements
of order~4 in~$Q_8$).

\begin{prop} \label{es32}
Let $G=D_8*Q_8$ be the extraspecial group of order \/{\rm 32 (}\/type~{\rm 2).} Then
$$D_t^0(G)\cong \bZ/4\bZ \oplus \bZ/2\bZ$$
generated by the class of
the module $$\Ten_{N_G(S_1)}^G\,\Inf_{N_G(S_1)/S_1}^{N_G(S_1)}
(\Omega_{N_G(S_1)/S_1}^1(k)) \quad\hbox{ \/{\rm (}\/order~{\rm 4)}}
$$
 and by the class
$$\Ten_{N_G(S_1)}^G\,\Inf_{N_G(S_1)/S_1}^{N_G(S_1)}
(M_{N_G(S_1)/S_1})\quad \hbox{ \/{\rm (}\/order~{\rm 2),}}
$$
 where $M_{N_G(S_1)/S_1}$ is
the module $M$ viewed as a module for the group $N_G(S_1)/S_1${\rm ,} which is
isomorphic to~$Q_8$.
\end{prop}

\Proof 
We have that $N_G(S_i)=S_i\times C$ (for some subgroup~$C$ isomorphic
to~$Q_8$)
and so
$N_G(S_i)/S_i\cong Q_8$ and $S_i$ is in the class~$\CX$ of
Theorem~\ref{Dadedetection}. These are the only subgroups in~$\CX$,
because every other nontrivial subgroup
$H$ contains the Frattini subgroup and $G/H$ is elementary abelian.
Therefore, by Theorem~\ref{Dadedetection}, the map
$$
\prod_{i=1}^5\; \Def_{N_G(S_i)/S_i}^{N_G(S_i)}\;
\Res_{N_G(S_i)}^G : D_t(G) \Rarr{}
\prod_{i=1}^5 D_t(N_G(S_i)/S_i) \cong (D_t(Q_8))^5
$$
is injective. 
Now by definition $D_t^0(G)$ is generated by the classes of the modules
$$\Ten_{N_G(S_i)}^G\,\Inf_{N_G(S_i)/S_i}^{N_G(S_i)}
(X)\quad (1\leq i\leq5)\,,$$
where $X$ is either $\Omega^1(k)$ or $M$ (viewed in
$D_t(N_G(S_i)/S_i)\,$).

If $X=\Omega^1(k)$, we always obtain the same element, independently
of~$i$, mapping to the diagonal element consisting of $\Omega^1(k)$
in each component under the injective map above. The proof of this
follows exactly the same argument as the one used in the proof of
Proposition~\ref{aes16}, with the following minor modification. For every
pair $S_i, S_j$ with $i\neq j$, the group generated by $S_i$ and $S_j$ is
isomorphic to~$D_8$.  Its centralizer~$C$ is isomorphic to~$Q_8$, and we
have $N_G(S_i)=S_i\times C$ and $N_G(S_j)=S_j\times C$. It follows that
we can use Lemma~\ref{formula} (with $P/Q=N_G(S_i)/S_i\,$,
$K/H=N_G(S_j)/S_j\,$, $P\cap K =C$). The rest of the argument is similar
to that used in Proposition~\ref{aes16}.

If now $X=M$, we  again use Lemma~\ref{formula}, but the computation
changes because of the presence of the Galois automorphism~$\gamma_2$
which does not fix the class of~$M$. Moreover, for each~$i$, we need to
fix a choice of isomorphism $N_G(S_i)/S_i \cong Q_8$ in order to be able
to make a consistent computation. We skip the details and only give the
result. It turns out that, under the injective map above, the image of
${\Ten_{N_G(S_i)}^G\,\Inf_{N_G(S_i)/S_i}^{N_G(S_i)}(M)}$
is the 5-tuple $(M,M,M,M,M)$, again independent  of~$i$. It follows that
we obtain just one extra generator of~$D_t^0(G)$, of order~2.
\Endproof\vskip4pt  

In this example, $D_t(G)$ is sandwiched between $D_t^0(G) \cong \bZ/4\bZ
\oplus \bZ/2\bZ$ and $(\bZ/4\bZ \oplus \bZ/2\bZ)^5$. The
question of the equality $D_t(G)=D_t^0(G)$ remains open.

\end{document}